\documentstyle[amscd,amssymb,diagrams,verbatim,12pt]{amsart}
\newcommand{\An}{\operatorname{An}}
\newcommand{\Per}{\operatorname{Per}}
\newcommand{\ab}{\operatorname{ab}}
\newcommand{\pa}{\partial}

\newcommand{\IT}{\operatorname{IT}}

\newcommand{\bH}{{\bf H}}
\newcommand{\bK}{{\bf K}}

\newcommand{\liminv}{\varprojlim}

\newcommand{\ch}{\operatorname{ch}}

\renewcommand{\mod}{\operatorname{mod}}

\newcommand{\gr}{\operatorname{gr}}

\newcommand{\Spf}{\operatorname{Spf}}
\newcommand{\und}{\underline}

\newcommand{\OO}{{\cal O}}
\newcommand{\Sym}{\operatorname{Sym}}

\newcommand{\lleft}{\operatorname{left}}
\newcommand{\DD}{{\cal D}}
\newcommand{\NN}{{\cal N}}

\newcommand{\bL}{{\bf L}}
\newcommand{\II}{{\cal I}}

\newcommand{\BB}{{\cal B}}

\newcommand{\hra}{\hookrightarrow}
\newcommand{\lan}{\langle}
\newcommand{\ran}{\rangle}
\newcommand{\Coh}{\operatorname{Coh}}
\newcommand{\GG}{{\cal G}}
\newcommand{\CC}{{\cal C}}

\newcommand{\Spec}{\operatorname{Spec}}

\newcommand{\Mat}{\operatorname{Mat}}

\renewcommand{\P}{{\Bbb P}}

\newcommand{\si}{\sigma}

\newcommand{\Pic}{\operatorname{Pic}}

\newcommand{\de}{\delta}

\renewcommand{\ker}{\operatorname{ker}}
\newcommand{\im}{\operatorname{im}}
\newcommand{\D}{{\cal D}}

\numberwithin{equation}{subsection}

\newcommand{\GL}{\operatorname{GL}}

\newtheorem{thm}{Theorem}[subsection]
\newtheorem{prop}[thm]{Proposition}

\newtheorem{lem}[thm]{Lemma}
\newtheorem{cor}[thm]{Corollary}
{  \theoremstyle{definition}
\newtheorem{defi}[thm]{Definition}
\newtheorem{ex}[thm]{Example}
\newtheorem{exs}[thm]{Examples}
\newtheorem{rem}[thm]{Remark}

}

\newcommand{\Pf}{\noindent {\it Proof}}
\newcommand{\id}{\operatorname{id}}
\newcommand{\Lie}{\operatorname{Lie}}
\newcommand{\LLie}{{\cal L}ie}
\newcommand{\ov}{\overline}
\newcommand{\we}{\wedge}

\newcommand{\ra}{\rightarrow}
\renewcommand{\AA}{{\cal A}}

\newcommand{\FF}{{\cal F}}
\newcommand{\EE}{{\cal E}}
\newcommand{\JJ}{{\cal J}}

\newcommand{\TT}{{\cal T}}

\newcommand{\PP}{{\cal P}}

\newcommand{\VV}{{\cal V}}

\newcommand{\LL}{{\cal L}}

\newcommand{\Om}{\Omega}
\newcommand{\Mor}{\operatorname{Mor}}
\newcommand{\bimod}{\operatorname{bimod}}

\newcommand{\dbar}{\overline{\partial}}

\newcommand{\Hom}{\operatorname{Hom}}

\newcommand{\Ext}{\operatorname{Ext}}
\newcommand{\End}{\operatorname{End}}

\newcommand{\Aut}{\operatorname{Aut}}

\renewcommand{\a}{\alpha}
\renewcommand{\b}{\beta}
\newcommand{\om}{\omega}
\newcommand{\De}{\Delta}
\newcommand{\la}{\lambda}

\newcommand{\C}{{\Bbb C}}
\newcommand{\N}{{\Bbb N}}
\newcommand{\R}{{\Bbb R}}
\newcommand{\Z}{{\Bbb Z}}

\newcommand{\La}{\Lambda}
\newcommand{\Ga}{\Gamma}

\newcommand{\wt}{\widetilde}
\newcommand{\ot}{\otimes}

\newcommand{\ad}{\operatorname{ad}}
\newcommand{\sub}{\subset}
\newcommand{\ed}{\qed\vspace{3mm}}

\newcommand{\Qcoh}{\operatorname{Qcoh}}

\newcommand{\dlb}{\langle\langle}
\newcommand{\drb}{\rangle\rangle}

\newsymbol\k 207C

\title{DG-resolutions of NC-smooth thickenings and NC-Fourier-Mukai transforms}
\author{Alexander Polishchuk and Junwu Tu}
\thanks{A.P. is partially supported by the NSF grant DMS-1001364}

\begin{document}
\begin{abstract} We give a construction of NC-smooth thickenings (a notion defined by Kapranov~\cite{Kapranov}) of a smooth variety equipped with a torsion free connection. We show that a twisted version
of this construction realizes all NC-smooth thickenings as $0$th cohomology of a differential graded
sheaf of algebras, similarly to Fedosov's construction in \cite{Fed}. We use this dg resolution to construct and study sheaves on NC-smooth thickenings. In particular, we construct an NC version of the Fourier-Mukai
transform from coherent sheaves on a (commutative) curve to perfect complexes on the canonical NC-smooth thickening of its Jacobian. We also define and study analytic NC-manifolds. We prove NC-versions
of some of GAGA theorems, and give a $C^\infty$-construction of analytic NC-thickenings that can
be used in particular for K\"ahler manifolds with constant holomorphic sectional curvature. Finally, we
describe an analytic NC-thickening of the Poincar\'e line bundle for the Jacobian of a curve,
and the corresponding Fourier-Mukai functor, in terms of $A_\infty$-structures.
\end{abstract}
\maketitle

\section{Introduction}

\subsection{Background.} In the remarkable paper \cite{Kapranov} Kapranov initiated the study of noncommutative algebras and
schemes which are complete with respect to the commutator filtration (called NC-complete algebras
and NC-schemes). 
The commutator filtration on a ring $R$ is
defined by
\begin{equation}\label{F-filtration-eq}
F^dR=\sum_{n_1+\ldots+n_m-m=d} R\cdot R^{\Lie}_{n_1}\cdot R\ldots\cdot R\cdot R^{\Lie}_{n_m}\cdot R,
\end{equation}
where $R^{\Lie}$ is $R$ viewed as a Lie algebra, $R^{\Lie}_n$ is the $n$th term of its lower central series.
In other words, $F^dR$ is the two-sided ideal spanned by all expression containing $d$ 
commutators (possibly nested). The ring $R$ is called
NC-complete if it is complete with respect to the filtration $(F^dR)$. Such a ring can be viewed as
an inverse limit of the {\it NC-nilpotent} rings $R/F^dR$ which are nilpotent extensions of
the abelianization of $R$, $R^{\ab}=R/F^1R$. {\it NC-schemes} are defined in \cite{Kapranov} as ringed spaces constructed by gluing formal spectra of NC-complete rings. Every NC-scheme $X$ has the abelianization $X^{\ab}$ which is a usual scheme, and we say that $X$ is an {\it NC-thickening} (or simply {\it thickening})
of $X^{\ab}$ (in particular, $X$ and $X^{\ab}$ have the same underlying topological spaces).
Intuitively one can think of NC-schemes as analogues of formal schemes where the formal
direction is allowed to be noncommutative. 

The important notion of {\it NC-smoothness} is introduced in \cite{Kapranov} by analogy with
the notion of smoothness of associative algebras due to Cuntz-Quillen \cite{CQ}. The only difference
is that in Kapranov's setup the lifting property is required to hold only in the class of NC-nilpotent algebras. 
One can think of an NC-smooth thickening of a smooth variety $X$ as an analog of deformation quantization,
where the algebra of functions on $X$ is replaced by its universal Poisson envelope. 
Kapranov showed in \cite{Kapranov} that every affine smooth variety $X$ admits
a unique NC-smooth thickening (up to a non-canonical isomorphism).
For non-affine smooth schemes there are in general obstructions to the existence of NC-smooth thickenings,
and they are not necessarily unique.  
Kapranov in \cite{Kapranov} constructs examples of NC-smooth thickenings for
projective spaces, Grassmannians and for some versal families of vector bundles.

There was a very limited development of the theory of NC-smooth thickenings after Kapranov
(see \cite{Cortinas-deRham}, \cite{Cortinas}, \cite{LeBruyn}).
In particular, before our work all constructions of smooth NC-thickenings followed
one of the two patterns: either step-by-step thickening using universal central extensions (this works
in the affine case), or via defining a functor on the category of NC-nilpotent algebras and proving its
representability.

\subsection{Summary of main results I: algebraic theory.} In the present paper we give a different construction of NC-smooth thickenings which is
analogous to Fedosov's construction of the deformation quantization of a symplectic manifold in \cite{Fed}. Our construction begins with the following definition.

\begin{defi}
\label{alge-nc-defi}
Let $X$ be a smooth algebraic variety. Let $\Omega^\bullet$ denote the algebraic de Rham complex of $X$, and $\hat{T}_\OO (\Om^1)$ the completed tensor algebra generated by one forms. An \emph{algebraic NC-connection} is a \emph{degree one, square zero derivation} $D$ of the graded algebra (graded by the de Rham degree)
\[ \AA_X:=\Om^\bullet\ot_{\OO} \hat{T}_\OO(\Om^1)\]
extending the de Rham differential on $\Om^\bullet$ and acting on $\Om^1\sub\hat{T}_\OO(\Om^1)$ by an operator of the form
\[ D(1\otimes \alpha)=\alpha\ot 1+\nabla_1(\alpha)+\nabla_2(\alpha)+\nabla_3(\alpha)+\cdots,\]
where $\nabla_i(\alpha)\in \Om^1\ot T^i(\Om^1)$.
\end{defi}

Note that $\nabla_1: \Om^1\to \Om^1\ot\Om^1$ is an algebraic connection on $\Om^1$, while for $i\geq 2$ 
the maps $\nabla_i$  are $\OO$-linear. The condition $D^2=0$ implies that
the $\nabla_1$ corresponds to a torsion free connection on $X$.
Conversely, we show that every torsion free connection extends to an
algebraic NC-connections and that the corresponding differential graded algebra is a resolution of 
an NC-smooth thickening of $X$.

\begin{thm}
\label{main-nc-thm}
Let $X$ be a smooth scheme endowed with an algebraic torsion free connection $\nabla$. Then
\begin{itemize}
\item[(1.)] there exists an algebraic NC-connection $D$ on $X$ with $\nabla_1=\nabla$;
\item[(2.)] one has $\underline{H}^i(\AA_X,D)=0$ for $i\geq 1$;
\item[(3.)] the sheaf $\underline{H}^0(\AA_X,D)$ 
with its induced algebra structure is an NC-smooth thickening of $X$.
\end{itemize}
Furthermore, such an NC-connection $D$ is determined uniquely up to conjugation by an automorphism
of $\hat{T}_\OO(\Om^1)$.
\end{thm}

Since every smooth affine variety admits a torsion-free connection, we get an explicit construction
of an NC-smooth thickening in this case, along with a dg-resolution.
In the non-affine case we prove that every NC-smooth thickening has a dg-resolution obtained by a twisted version of our construction, involving a sheaf of algebras locally isomorphic to $\hat{T}_\OO(\Om^1)$.

For an abelian variety $A$ our construction applied to the natural connection gives an NC-smooth
thickening that we call {\it standard}. 
If $J=J(C)$ is the Jacobian of a smooth projective curve $C$ then there is another NC-smooth thickening
of $J$ constructed by Kapranov in \cite{Kapranov} using the functor of noncommutative families of line
bundles on $C$. We prove that this modular thickening is isomorphic to the standard one for $J$
viewed as an abelian variety.

Another theme of this paper is the study of $\OO_X^{NC}$-modules for an NC-smooth thickening
$\OO_X^{NC}\to\OO_X$. We are partially guided in this study by the analogy with modules over
algebraic (or complex analytic) deformation quantization algebras that have been studied e.g. in \cite{KS}, \cite{BBP}, \cite{ABP}, \cite{BGP}, \cite{Pecharich}.
In the case when $\OO_X^{NC}$ is obtained from a torsion-free connection on $X$
we mimic our construction of the dg-resolution for modules.
Namely, starting with a quasicoherent sheaf $\FF$ on $X$ with a connection (not necessarily flat)
we construct a dg-module over $\AA_X$ that does not depend on a connection up to an isomorphism. 
Passing to the $0$th cohomology we get an $\OO_X^{NC}$-module, which is locally free over $\OO_X^{NC}$ in the case when $\FF$ is a vector bundle.
For an arbitrary NC-smooth thickening $\OO_X^{NC}$ a similar construction works if the connection on $\FF$
is flat. Furthermore, in this way we get an exact functor from 
$D$-modules on $X$ to $\OO_X^{NC}$-bimodules.

One can ask which vector bundles extend to locally free modules
(resp., bimodules) over a given NC-smooth thickening $\OO_X^{NC}\to\OO_X$. 
The first order obstruction to extend a vector bundle on $X$ to an $\OO_X^{NC}$-module
was computed by Kapranov in terms of the Atiyah classes (see \cite[Rem.\ (4.6.9)]{Kapranov}).
We extend this to the case of bimodules and also answer the question of extendability
completely in the case of the standard NC-smooth thickening of an abelian variety  $A$
and of line bundles on $A$.

In the case of abelian varieties it is natural to look for analogs of the Fourier-Mukai
transform involving their NC-smooth thickenings.
We show how to extend the Poincar\'e line bundle on $A\times\hat{A}$
to a line bundle (a bimodule) on $A^{NC}\times A^\natural$, where
$A^{NC}$ is the standard NC-smooth thickening of $A$ and
$A^\natural$ is the universal extension of $\hat{A}$ by a vector space
(the functions on $A^\natural$ are in the center). 
One can descend this NC-Poincar\'e line bundle to a twisted line bundle on $A^{NC}\times\hat{A}$
(since the projection $A^\natural\to \hat{A}$ splits locally).
If we restrict to a curve $C$ in $\hat{A}$ then we can untwist it and get a line bundle
on $A^{NC}\times C$. Using this line bundle we define an integral transform from
the derived category $D^b(\Coh(C))$ to the derived category of quasicoherent sheaves on $A^{NC}$.
We prove that objects in the image of this transform are perfect, i.e., can be represented by
complexes of vector bundles on $A^{NC}$. 

\subsection{Summary of main results II: analytic theory.}The second main topic of this paper may be 
called \emph{analytic NC geometry}. One advantage of the dg resolution $(\AA_X,D)$ of an NC-smooth scheme $X^{NC}$ considered above is that it enables us to define the \emph{analytification} of $X^{NC}$: simply take the analytification of $(\AA_X,D)$, and then apply $\underline{H}^0$. For a general complex analytic manifold $X$ (which might not admit any holomorphic torsion-free connection), we define an analytic NC-smooth thickening of $X$ to be a sheaf of algebras over $X$ that is locally in the analytic topology modeled by an algebra of the form $\underline{H}^0(\AA,D)$.

Similar to the case of complex geometry, in analytic NC geometry, it is often useful to work with $C^\infty$ bundles and their Dolbeault complexes. The following Definition generalizes Definition~\ref{alge-nc-defi}.

\medskip
\begin{defi}~\label{nc-connection-def}
Let $X$ be a complex manifold. Let $\Lambda_X$ denote the $C^\infty$ de Rham algebra of $X$, and $\Omega^{1,0}_X$ the space of $(1,0)$-type forms. An \emph{NC-connection} on $X$ is a \emph{degree one, square zero derivation} $\DD$ of the graded algebra (graded by the de Rham degree of $\Lambda_X$) 
\[\Lambda_X\otimes_{C^\infty_X} \hat{T}_{C_X^\infty}(\Omega_X^{1,0})\] 
extending the de Rham differential on $\Lambda_X$, and acting on 
$\Omega_X^{1,0}\sub\hat{T}_{C_X^\infty}(\Omega_X^{1,0})$ 
by an operator of the form 
\[\DD(1\ot \alpha)=\alpha\ot 1+\overline{\partial}+\nabla+A^\vee_2+A^\vee_3+\cdots\]
where $\nabla$ is a $(1,0)$-type connection on the bundle $\Omega_X^{1,0}$, and $A^\vee_k: \Omega_X^{1,0} \ra \Lambda^1_X\ot (\Omega_X^{1,0})^{\ot k}$ is a $C^\infty_X$-linear morphism for each $k\geq 2$. 
\end{defi}

Analogously to Theorem~\ref{main-nc-thm}, we prove the following results about NC connections, which demonstrates the important role they play in analytic NC geometry.

\medskip

\begin{thm}~\label{analytic-nc-thm}
Let $X$ be a complex manifold. Then 
\begin{itemize}
\item[(1.)] one has $\underline{H}^i(\Lambda_X\otimes_{C^\infty_X} \hat{T}_{C_X^\infty}(\Omega_X^{1,0}),\DD)=0$ for $i\geq 1$;
\item[(2.)] the sheaf $\underline{H}^0(\Lambda_X\otimes_{C^\infty_X} \hat{T}_{C_X^\infty}(\Omega_X^{1,0}),\DD)$ with its induced algebra structure is an analytic NC-smooth thickening of $X$;
\item[(3.)] there exists an NC-connection $\DD$ on $X$ if and only if $X$ admits an analytic NC-smooth thickening.
\end{itemize}
\end{thm}

\begin{thm}\label{examples-thm}
Let $(X,g)$ be a K\"ahler manifold with constant holomorphic sectional curvature. Then $X$ admits an NC-connection. Thus, by Theorem~\ref{analytic-nc-thm}, it also admits an analytic NC-smooth thickening.
\end{thm}

K\"ahler manifolds with constant sectional curvature are classified. There are three classes of such manifolds: quotients of $(\C^n,\omega_0=\sqrt{-1}/2\sum dz_id\overline{z}_j)$, complex projective spaces, and complex hyperbolic spaces. Of these three classes of examples, the first one admits holomorphic torsion-free connection, hence the existence of an NC-smooth thickening also follows from the holomorphic
version of Theorem~\ref{main-nc-thm}. 
An NC-smooth thickening of the projective space was constructed by Kapranov~\cite[Section 5.1]{Kapranov}.
Our construction for hyperbolic spaces seems to give new examples of NC-smooth thickenings.

\medskip

It is natural to ask whether analogs of GAGA results (see \cite{GAGA}) hold for algebraic NC geometry versus analytic NC geometry. 
We first deal with the algebraicity of analytic NC manifolds. Let $X$ be a smooth algebraic variety over $\C$. Following Kapranov~\cite{Kapranov}, we denote by $Th_X$ the groupoid whose objects are NC-smooth thickenings of $X$, and morphisms are morphisms of locally ringed spaces that are identity after taking Abelianization. In the analytic setting, we have the corresponding groupoid $Th_{X^{an}}$ consisting of analytic NC-smooth thickenings of $X^{an}$.

\begin{thm}
\label{main-an-thm}
Let $X$ be a smooth projective variety over $\C$. Then the analytification functor is an equivalence of categories $Th_X \ra Th_{X^{an}}$.
\end{thm}

Next we consider analytic NC vector bundles versus algebraic ones. Again let $X$ be a smooth projective variety over $\C$. Let $X^{NC}$ be an algebraic NC-smooth thickening of $X$, and $X^{nc}$ its analytification. Denote by $\Per^{gl}(X^{NC})$ the (global) perfect derived category of locally free $\OO_X^{NC}$-modules of finite rank, and by $\Per^{gl}(X^{NC})$ the similar perfect derived category for $\OO_{X^{an}}^{nc}$-modules.
Since the analytic topology is finer than the Zariski topology, and algebraic sections are holomorphic sections, there is a canonical functor
\[
\An: \Per^{gl}(X^{NC}) \ra \Per^{gl}(X^{nc}).
\]

\begin{thm}\label{nc-gaga-thm}
The functor $\An$ is an equivalence of categories.
\end{thm}

Thus, we can also use complex analytic methods to study NC vector bundles. As an example we present an analytic construction of the NC thickening of the Poincar\'e bundle $P$ over $C\times J$ studied earlier with algebraic methods. Using the analytic description we clarify the relationship between the NC Fourier transform of $P$ in a formal neighborhood of a point $L\in J$ and another such functor constructed by the first named author~\cite{P-BN} using different methods.

The paper is organized as follows. In Section \ref{connection-constr-sec}, after reviewing basic notions
on NC-thickenings, we present the Fedosov type
construction of the NC-connection associated with a torsion free connection. Then we show that
an NC-connection gives rise to an NC-smooth thickening $\OO_X^{NC}$
and give a proof of Theorem \ref{main-nc-thm}.
We also show in Theorem \ref{twisted-dg-thm}
that all NC-smooth thickenings can be described in terms of twisted NC-connections.

In Section \ref{module-sec} we consider the analog of the construction of NC-connections for $\OO_X$-modules.
In \ref{D-mod-sec} we construct an exact functor from $D$-modules on $X$ to $\OO_X^{NC}$-bimodules
and in \ref{extendability-sec} we study the question of extendability of vector bundles on $X$ to locally
free $\OO_X^{NC}$-modules (resp., bimodules). 

In Section \ref{Jac-FM-sec} we study NC-thickenings of the Jacobian $J=J(C)$ of a curve $C$. We show that
the canonical NC-smooth thickening $J^{NC}$ of $J$ constructed by Kapranov~\cite[Sec.\ 5.4]{Kapranov}
is isomorphic to the standard NC-thickening of $J$ as an algebraic variety. In
\ref{NC-Poincare-sec} we construct an NC-thickening of the Fourier-Mukai functor from $D^b(C)$ to the perfect
derived category of the NC-thickening of $J$.

In Section \ref{NC-analytic-sec} we define and study analytic NC manifolds. We prove Theorem 
\ref{main-an-thm}
and consider the relation between analytic NC-manifodls and $C^\infty$ NC-connections
(see Definition \ref{nc-connection-def}). We also construct such NC-connections on
K\"ahler manifolds with constant sectional curvature.

In Section \ref{analytic-vb-sec} we study the analytification functor for vector bundles over NC-smooth
thickenings and prove Theorem \ref{nc-gaga-thm}.

Finally, in Section \ref{ainf-sec} we use the analytic approach to describe the Poincar\'e line bundle
over $C\times J^{NC}$ and the corresponding Fourier-Mukai functor using $A_\infty$-structures.

\noindent
{\it Conventions.} We work over a base field $k$ of characteristic zero.
Starting from Section \ref{NC-analytic-sec} we assume that $k=\C$.
For a commutative ring $A$ we denote by $A\dlb x_1,\ldots,x_n\drb$ the ring
of noncommutative power series in independent variables $x_1,\ldots,x_n$ with coefficients in $A$.
For a smooth variety we set $D^b(X)=D^b(\Coh(X))$. By a curve we mean a smooth projective curve.
For an abelian group  (resp., a sheaf of abelian groups) $A$, $\Mat_n(A)$ denotes the group
(resp., sheaf) of $n\times n$-matrices with values in $A$. For a sheaf of rings $\OO$ we denote
by $\mod-\OO$ the category of right $\OO$-modules. For a ringed space $X$ we denote 
the derived category of right $\OO_X$-modules by $\DD(X)$.

\noindent
{\it Acknowledgments.} Part of this work was done while the first author was visiting the 
Mathematical Sciences Research Institute and the Institut des Hautes \'Etudes Scientifiques, 
and the second author was visiting the Chinese University of Hong Kong.
We would like to thank these institutions for excellent working conditions.
Also, we are grateful to Kevin Costello, Weiyong He and Dmitry Tamarkin for useful discussions.

\section{NC-thickenings from torsion free connections}\label{connection-constr-sec}

\subsection{NC-algebras and NC-schemes}

We start by recalling the main definitions from \cite{Kapranov}. Along the way we prove some
technical lemmas that will be needed later. 

\begin{defi} A $k$-algebra $R$ is called {\it NC-nilpotent of degree $d$} (resp., {\it NC-nilpotent}) 
if $F^{d+1}R=0$ (resp., $F^nR=0$ for $n\gg 0$).
The {\it NC-completion} of an algebra $R$ is 
$$R_{[[ab]]}=\liminv R/F^nR.$$ 
An algebra $R$ is {\it NC-complete} if it is complete with respect to the filtration $(F^nR)$,
i.e., the natural map $R\to R_{[[ab]]}$ is an isomorphism.
\end{defi}

For an NC-complete algebra $R$ the affine NC-scheme $X=\Spf(R)$ is defined
as the locally ringed space $(\Spec(R^{ab}),\OO_X)$, where $\OO_X$ is defined as the inverse
limit of sheaves of rings obtained from $R/F^dR$ by localization (similarly to the case of
usual affine schemes). The category of NC-schemes is the full subcategory of the category of
locally ringed spaces consisting of spaces locally isomorphic to affine NC-schemes.
For every NC-scheme $X$
the structure sheaf $\OO_X$ has a filtration by sheaves of ideals $F^nX$, and the quotients
$\OO_X/F^{n+1}\OO_X$ define NC-subschemes $X^{\le n}\sub X$
 with the same underlying topological space as $X$.
In particular, we set $X_{ab}=X^{\le 0}$.

\begin{lem}\label{NC-mor-lem}
To give a morphism $f:X\to Y$ of NC-schemes is the same as to give a collection of
compatible morphisms $f_n:X^{\le n}\to Y$.
\end{lem}

\Pf . Since the construction $X\mapsto X^{\le n}$ is compatible with passing to open subsets,
it is enough to consider the case when both $X$ and $Y$ are affine. In this case, by \cite[Prop.\ 2.2.3]{Kapranov},
our statement is that a homomorphism of NC-complete algebras $R\to S$ is determined by a
collection of compatible homomorphisms $R\to S/F^nS$, which follows immediately from
NC-completeness of $S$. 
\ed

\begin{defi}
(i) An NC-scheme $X$ is said to be of {\it finite type} if it is NC-nilpotent, i.e. $F^n\OO_X=0$ for $n\gg 0$,
$X_{ab}$ is a scheme of finite type over $k$ and the sheaves $F^n\OO_X/F^{n+1}\OO_X$ on $X_{ab}$
are coherent.

\noindent
(ii) An NC-scheme $X$ is called {\it $d$-smooth} if $F^{d+1}X=0$, $X$ is of finite type and for every extension
of NC-nilpotent algebras
\begin{equation}\label{nilp-extension-eq}
0\to I\to \La'\to \La\to 0
\end{equation}
where $I$ is a nilpotent ideal, the map $\Mor(\Spf(\La'),X)\to \Mor(\Spf(\La),X)$ is surjective.
An algebra $R$ is called $d$-smooth if $\Spf(R)$ is $d$-smooth.

\noindent
(iii) An NC-scheme $X$ is called {\it NC-smooth} if $X^{\le d}$ is $d$-smooth for every $d$.
An algebra $R$ is called NC-smooth if $R/F^{d+1}R$ is $d$-smooth for every $d$.
\end{defi}

It is easy to see that if an NC-scheme $X$ is $d$-smooth then $X^{\le d-1}$ is $(d-1)$-smooth.
In particular, $X_{ab}$ is smooth as a commutative scheme over $k$. If $X$ is a $d$-smooth (resp. NC-smooth)
NC-scheme then we say that $X$ is a {\it $d$-smooth (resp., NC-smooth) thickening} of the commutative smooth scheme $X_{ab}$.

Recall that an extension \eqref{nilp-extension-eq} is called central if $I^2=0$ and $I$ lies in the center of
$\La'$.
It is easy to see that any surjection $\La'\to \La$ of NC-nilpotent algebras whose kernel is a nilpotent 
ideal is a composition of central extensions (cf. \cite[Prop.\ 1.2.3]{Kapranov}. Thus,
in the definition $d$-smoothness it is enough to work with central extensions of NC-nilpotent algebras.

\begin{lem}\label{smooth-local-lem} Let $X$ be an NC-scheme. If $X$ is locally NC-smooth,
i.e., admits an open covering $(U_i)$ such that $U_i$ are NC-smooth, then $X$ is
NC-smooth.
\end{lem}

\Pf . Given a central extension $0\to I\to\La'\to\La\to 0$ of NC-nilpotent algebras we need to check
that any map $f:\Spec(\La)\to X$ extends to a map $\Spec(\La')\to X$. 
By assumption, we can assume that there is an affine open covering $(\Spec(\La_i))$
of $\Spec(\La)$ such that the corresponding restrictions $f_i:\Spec(\La_i)\to X$ extend to
$f'_i:\Spec(\La'_i)\to X$, where $0\to I_i\to \La'_i\to \La_i\to 0$ is the induced central extension
of $\La_i$. On intersections the liftings $f'_i$ and $f'_j$ differ by a derivation of $\OO_X$ with
values in $I$. Since $I$ is a central bimodule, such derivations correspond to homomorphisms
$f^*\Om_X^1\to I$. Thus, we get a Cech $1$-cocycle on $\Spec(\La^{ab})$ with values in
$(f^*\Om_X^1)^\vee\ot I$. Since $H^1(\Spec(\La^{ab},(f^*\Om_X^1)^\vee\ot I)=0$,
we can make the liftings $f'_i$ to be compatible on intersections.
\ed

Let $R$ be a $d$-smooth algebra. Then $R$ is a universal central extension of the $(d-1)$-smooth
algebra $R/F^dR$ (see \cite[Thm.\ 1.6.1, Prop.\ 1.6.2]{Kapranov}). For example, a $1$-smooth algebra
$R$ is a universal central extension of $R_{ab}$. The universality is equivalent to the property that
the commutator bracket $R_{ab}\times R_{ab}\to F^1R$ is the universal skew-symmetric biderivation,
so $F^1R\simeq\Om^2_{R_{ab}}$ (see \cite[Sec.\ 1.3]{Kapranov}). This analogous property holds
in the nonaffine setting. Namely, a central extension of sheaves of $k$-algebras
$$0\to \II\to \wt{\OO}\to\OO_Y\to 0$$
on a smooth commutative scheme $Y$ is $1$-smooth if and only the commutator pairing
$$(f,g)\mapsto [\wt{f},\wt{g}]\in \II,$$
(where $\wt{f},\wt{g}\in\wt{\OO}$ are any liftings of $f$ and $g$) viewed as a
biderivation with values in $\II$ induces an isomorphism $\Om^2_X\to \II$.

\begin{defi}\label{standard-1-thickening}(cf. \cite[Ex.\ (1.3.9)]{Kapranov}) Let $X$ be a smooth scheme. The {\it standard $1$-smooth thickening} of $X$
is the sheaf of $k$-vector spaces $\wt{\OO}=\Om^2_X\oplus \OO_X$ equipped with the product
$$(\a,f)\cdot (\a',f')=(f'\a+f\a'+df\we dg, ff').$$
\end{defi}

Thus, for any NC-smooth scheme $X$ one has a natural isomorphism 
$$F^1\OO_X/F^2\OO_X\simeq\Om^2_{X_{ab}}.$$
The description of the quotients $F^n\OO_X/F^{n+1}\OO_X$ is more complicated. Kapranov proves that 
$$F^n\OO_X/F^{n+1}\OO_X\simeq P_n(\OO_{X_{ab}}),$$
where $P(R)=\bigoplus_{n\ge 0}P_n(A)$ denotes the Poisson envelope of a commutative
algebra $A$. However, the description of $P_n(A)$ for a smooth algebra $A$ in terms of some polynomial
functors proposed in \cite[Thm.\ 4.1.3]{Kapranov} is only known to hold locally (see 
\cite[Thm.\ 1.4]{Cortinas}; the proof of the general case in \cite{Kapranov} has a gap).
In this work we will consider another algebra filtration on the structure sheaf of an NC-smooth scheme
for which the graded quotients will have a simple descriptions as polynomial functors.

For any ring (or a sheaf of rings) $R$ and $n\ge 1$ we define the two-sided ideal 
\begin{equation}\label{I-n-def-eq}
I_n(R)=\sum_{i_1\ge 2,\ldots,i_m\ge 2,i_1+\ldots+i_m\ge n} R\cdot R_{i_1}^{\Lie}\cdot R\cdot\ldots\cdot R\cdot R^{\Lie}_{i_m}\cdot R.
\end{equation}
Note that $R_i^{\Lie}\sub F^{i-1}R$, so $I_n(R)\sub F^{[n/2]}R$.
On the other hand, it is clear from the definition that $F^{n-1}R\sub I_n(R)$ for $n\ge 2$. Thus,
$(I_n(R))$ define the same topology on $R$ as $(F^nR)$.

For any ring $R$ the filtration $I_n(R)$ is an algebra filtration, so we can
consider the associated graded algebra 
$$\gr^\bullet_I(R)=\bigoplus_{n\ge 0} I_n(R)/I_{n+1}(R),$$
where we set $I_0(R)=R$. Note that $\gr^0_I(R)=R^{ab}$.
Note also that $[R,I_n(R)]\sub I_{n+1}(R)$, so $R^{ab}$ is in the center
of $\gr^\bullet_I(R)$.

For a vector space $V$ we denote by $\LLie(V)=\bigoplus_{n\ge 1}\LLie_n(V)$ the free Lie algebra
generated by $V$ and we set 
\begin{equation}\label{Lie+eq}
\LLie_+(V)=\bigoplus_{n\ge 2}\LLie_n(V).
\end{equation}
There is a natural $\GL(V)$-action on $\LLie_n(V)$, so we can define the functor $\LLie_n(?)$
on the category of vector bundles over a scheme.
We consider the grading on $\LLie_+(V)$ such that $\LLie_n(V)$ lives in degree $n$.\footnote{This
grading is different from the one considered in \cite{Kapranov}, cf. Remark \ref{filtrations-rem}.}
We equip the universal enveloping algebra $U(\LLie_+(V))$ with the induced grading, so that
$$U(\LLie_+(V))_n\sub U(\LLie(V))_n=V^{\ot n}.$$

\begin{lem}\label{homomorphism-lem} 
Assume that $R^{ab}$ is smooth. There is a natural surjective homomorphism
of graded algebras
\begin{equation}\label{Lie-map-eq}
U(\LLie_+^{R^{ab}}(\Om^1_{R^{ab}}))\to gr^\bullet_I(R)
\end{equation}
sending $f_0[df_1,[df_2,\ldots,[df_{n-1},df_n]\ldots]]$ to 
$\wt{f_0}[\wt{f_1},[\wt{f_2},\ldots,[\wt{f_{n-1}},\wt{f_n}]\ldots]]$,
where $\wt{f_i}\in R$ is a lifting of $f_i\in R^{ab}$.
\end{lem}

\Pf . First, one can easily check that the map
\begin{equation}\label{tensor-lie-map}
T^n(\Om^1_{R^{ab}})\to gr^n_I(R):
f_0df_1\ot df_2\ot\ldots\ot df_n\mapsto
\wt{f_0}[\wt{f_1},[\wt{f_2},\ldots,[\wt{f_{n-1}},\wt{f_n}]\ldots]]
\end{equation}
is well defined for $n\ge 2$. We can also define similar maps for iterated brackets in a different order. More
precisely, let $L=\bigoplus_n L(n)$ be the Lie operad. Then we get a natural map
$$L(n)\ot_{k[S_n]} T^n(\Om^1_{R^{ab}})\to gr^n_I(R).$$
Recall that $\LLie_n(V)=L(n)\ot_{S_n} T^n(V)$ (see \cite[Sec.\ 3]{KM}). Thus, we obtain a map of Lie
algebras
$$\LLie_+^{R^{ab}}(\Om^1_{R^{ab}})\to gr^\bullet_I(R)$$
which extends to \eqref{Lie-map-eq}.
To check surjectivity of \eqref{Lie-map-eq} it is enough to see that $\gr^\bullet_I(R)$ is generated
as $R^{ab}$-algebra by the images of $R_n^{\Lie}\sub I_n(R)$.
But this follows from the inclusion
$$R\cdot R_{i_1}^{\Lie}\cdot R\cdot R_{i_2}^{\Lie}\cdot \ldots\cdot R\cdot R^{\Lie}_{i_m}\cdot R\sub
R\cdot R_{i_1}^{\Lie}\cdot R_{i_2}^{\Lie}\cdot \ldots\cdot R^{\Lie}_{i_m}+I_{i_1+\ldots+i_m+1}(R),$$
which holds since commuting factors from $R$ to the left creates additional commutators.
\ed

We will show in Corollary \ref{I-filtration-gen-cor}
that the map \eqref{Lie-map-eq} is an isomorphism if and only if $R$ is NC-smooth.


\subsection{NC-thickening of the affine space}\label{NC-affine-space-sec}

Let us recall Kapranov's description of the completion $k\lan e_1,\ldots,e_n\ran_{[[ab]]}$ of
the free algebra $k\lan e_1,\ldots,e_n\ran$ with respect to the commutator filtration (see \cite[Sec.\ 3]{Kapranov}).
For a polynomial $P(y_1,\ldots,y_n)\in k[y_1,\ldots,y_n]$ we define
$[[P(e)]]\in k\lan e_1,\ldots,e_n\ran$ by replacing every monomial $cy_1^{i_1}\ldots y_n^{i_n}$ in $P$
by the ordered monomial $ce_1^{i_1}\ldots e_n^{i_n}$.
Let $V$ be the $k$-vector space with the basis $e_1,\ldots,e_n$, so we can identify
$k\lan e_1,\ldots,e_n\ran$ with the tensor algebra $T(V)$.
Let $\BB=\{\b_1,\b_2,\ldots\}$ be a $k$-basis of $\LLie_+(V)$ (see \eqref{Lie+eq}), ordered in such a
way that elements of $\LLie_n(V)$ precede elements of $\LLie_{n+1}(V)$.
For any function $\la:\BB\to\Z_{\ge 0}$ with finite support we consider the monomial 
$$M_\la=\b_1^{\la(\b_1)}\b_2^{\la(\b_2)}\ldots \in T(V)$$
(where we view $\LLie_+(V)$ as a subspace of $T(V)$).
The definition of the commutator filtration implies that
every infinite series
$$\sum_{\la:\BB\to\Z_{\ge 0}} [[a_\la(e)]] M_\la,$$
where $a_\la(y)\in k[y_1,\ldots,y_n]$, converges in $T(V)_{[[ab]]}$.
In fact, every element of $T(V)_{[[ab]]}$ can be written uniquely as such a series and 
a multiplication rule is given by \cite[Prop.\ 3.4.3]{Kapranov}.

We can rewrite the above description by realizing $T[V]_{[[ab]]}$ as a subalgebra in the algebra
of formal noncommutative power series in $e_1,\ldots,e_n$ over the commutative ring
$A=k[x_1,\ldots,x_n]$. Namely,
consider the completed tensor algebra of the free $A$-module $V\ot A$,
$$\hat{T}_A(V\ot A):=\prod_{n\ge 0}T^n_A(V\ot A)\simeq A\dlb{e_1,\ldots,e_n}\drb.$$
Let $\Om^{\bullet}_A=\Om^\bullet_{A/k}$ be the de Rham complex of $A$. 
We equip the algebra $\Om^\bullet_A\ot_A \hat{T}_A(V\ot A)$
with the differential $D$ which is a unique derivation extending the
de Rham differential on $\Om^\bullet_A$ such that
$$D(e_i)=dx_i \text{ for } i=1,\ldots,n.$$
Let $\pi:\hat{T}_A(V\ot A)\to A$ be the natural projection (sending all tensors of degree $\ge 1$ to zero).
Consider the $k$-algebra homomorphism
$$\de: T(V)\to \hat{T}_A(V\ot A): e_i\mapsto e_i-x_i.$$
Clearly the image is contained in the kernel of $D$, and $\pi\circ\de:T(V)\to A$ 
is the homomorphism sending $e_i$ to $-x_i$.

\begin{thm}~\label{local-constr-thm}
The homomorphism $\de$ extends to an isomorphism of the NC-completion 
$T(V)_{[[ab]]}$ with the subalgebra $\ker(D)\sub \hat{T}_A(V\ot A)$.
\end{thm}


\Pf . Note that any Lie word $w\in \LLie_+(V)$ in generators $e_i$ is mapped by $\de$
to itself viewed as an element in $T(V)\sub T_A(V\ot A)$. Hence, the same is true for
the monomials $M_\la$, where $\la:\BB\to Z_{\ge 0}$.
Thus, the homomorphism $\de$ extends to a homomorphism
$$\hat{\de}: T(V)_{[[ab]]}\to \hat{T}_A(V\ot A).$$
mapping the infinite series $\sum_{\la} [[a_\la]] M_\la$ 
to $\sum_{\la} [[a_\la(e_i-x_i)]] M_\la$ which converges in $\hat{T}_A(V\ot A)$
and still belongs to $\ker(D)$.
We need to check that the image of $\hat{\de}$ is the entire $\ker(D)$.
Note that any element $x\in\hat{T}_A(V\ot A)$ can be written as
$x=\sum_\la [[f_\la(e)]] M_\la$, where $f_\la(y)$ are some formal power series in commuting variables
$y_1,\ldots,y_n$ with coefficients in $A$. The condition $D(x)=0$ means that each $f_{\la}(y)$ is
of the form $f_{\la}(y)=a_{\la}(y_i-x_i)$ for some power series $a_{\la}$ with coefficients in $\C$. 
Note that $a_{\la}(-x)=f_\la(0)$ is the constant coefficient of $f_{\la}(y)$, 
hence $a_{\la}\in A$.
\ed


\subsection{Global construction}\label{NC-constr-sec}

Let $(V, \nabla)$ be a vector bundle with a connection over a smooth scheme $X$, 
and let $\iota:T\to V$ be a homomoprhism
of vector bundles. We denote still by $\nabla$ the induced connection on the dual bundle $V^*$.
Then we can 
consider the completed tensor algebra $\hat{T}_{\OO}(V^*)$ and define the differential $D=D_{\nabla,\iota}$ on 
the algebra $\Om^\bullet\ot_{\OO}\hat{T}_{\OO}(V^*)$ as a unique derivation,
extending the de Rham differential on $\Om^\bullet$ and satisfying
$$D(\alpha)=\nabla(\alpha)+\iota^*(\alpha)\ot 1 \text{ for } \alpha\in V^*.$$

\begin{lem}\label{torsion-free-lem} 
One has $D^2=0$ if and only if $\nabla$ is flat and
$\iota$ is torsion free, i.e., satisfies the equation
\begin{equation}\label{torsion-free-eq}
\iota([t_1,t_2])=\nabla_{t_1}(\iota(t_2))-\nabla_{t_2}(\iota(t_1)).
\end{equation}
\end{lem}

\Pf . Let $\wt{\nabla}:\Om^1\ot V^*\to\Om^2\ot V$ and $\wt{\iota^*}:\Om^1\ot V\to \Om^2$ be the maps given
by
$$\wt{\nabla}(\om\ot\alpha)=d\om\ot\alpha-\om\we\nabla(\alpha),\ \ \wt{\iota^*}(\om\ot\a)=
\om\we\iota*(\a).$$
Then 
$$D^2(\a)=\wt{\nabla}\nabla(a)+[d\iota^*(\a)-\wt{\iota^*}\nabla(\a)]\in \Om^2\ot V\oplus\Om^2.$$
Hence $D^2=0$ if and only if 
$\wt{\nabla}\nabla=0$ and 
$$d\circ\iota^*=\wt{\iota^*}\circ\nabla.$$
The former condition just means that $\nabla$ is flat. To show that the latter condition is equivalent
to \eqref{torsion-free-eq} we 
note that for a pair of local vector fields $t_1$, $t_2$ and for $\a\in V^*$ one has
$$d\iota^*(\a)(t_1, t_2)=t_1\cdot \a(\iota(t_2))-t_2\cdot\a(\iota(t_1))-\a(\iota([t_1,t_2])),$$
$$\wt{\iota^*}(\nabla(\a))(t_1,t_2)=\nabla_{t_1}(\a)(\iota(t_2))-\nabla_{t_2}(\a)(\iota(t_1)).$$
It remains to use the formulas 
$t_1\cdot \a(\iota(t_2))=\nabla_{t_1}(\a)(\iota(t_2))+\a(\nabla_{t_1}(\iota(t_2)))$ and a similar
formula for $t_2\cdot \a(\iota(t_2))$.
\ed

Now let us specialize to the case $V=T$ and $\iota=\id$, so that $\nabla$ is a a connection on the
tangent bundle. Assume that $\nabla$ is torsion free and flat.
Then the previous construction gives an algebraic NC-connection on $X$ in the sense of Definition
\ref{alge-nc-defi}.

Let $D$ be an arbitrary algebraic NC-connection on $X$, i.e., a differential  of degree $1$ on the algebra
$$\AA_X:=\Om^\bullet\ot_{\OO}\hat{T}_{\OO}(\Om^1),$$ 
where the grading on $\AA_X$ is given by the natural grading on $\Om^\bullet$,
such that
\[ D(1\otimes \alpha)=\alpha\ot 1+\nabla_1(\alpha)+\nabla_2(\alpha)+\nabla_3(\alpha)+\cdots\]
 for $\a\in\Om^1$.
The same calculation as in Lemma \ref{torsion-free-lem} shows that 
$\nabla_1$ corresponds to a torsion free connection on $X$.
In Theorem \ref{dg-constr-thm} below we show that conversely, any torsion free connection $\nabla$
(not necessarily flat) gives rise to an NC-connection with $\nabla_1=\nabla$.
The idea is to add higher terms to the differential that would kill the effect of the nonzero curvature.

Below we view the universal enveloping algebra $U(\LLie_+(V))$ as a graded subalgebra of 
$T(V)=U(\LLie(V))$. 

\begin{thm}\label{dg-constr-thm} 
(i) Let $\nabla$ be a torsion free connection on the tangent bundle on a smooth scheme $X$. Then there exists
a natural construction of an algebraic NC-connection $D$ with $\nabla_1=\nabla$, the induced connection
on $\Om^1$.

\noindent
(ii) For any two algebraic NC-connections $D,D'$ on $X$
(possibly for different torsion free connections)
there exists an algebra automorphism $\phi$ of $\hat{T}_\OO(\Om^1)$, such that
for $\a\in\Om^1$, $\phi(\a)=\a+\phi_2(\a)+\phi_3(\a)+\ldots$
with $\phi_i(\a)\in T^i(\Om^1)$, such that $D'=(\id\ot\phi) D (\id\ot\phi)^{-1}$,
where $\id\ot\phi$ is the induced automorphism of $\AA_X$.

\noindent
(iii) Let $G=G_D$ be the sheaf of groups of automorphisms $\phi$ of $\hat{T}_\OO(\Om^1)$, such 
for $\a\in\Om^1$, $\phi(\a)=\a+\phi_2(\a)+\phi_3(\a)+\ldots$ with $\phi_i(\a)\in T^i(\Om^1)$,
and such that $D=(\id\ot\phi) D (\id\ot\phi)^{-1}$. The group $G$ has a decreasing filtration by normal
subgroups 
$$G_n=\{\phi\in G\ |\ \phi_i=0 \text{ for } 2\le i\le n\}, $$
where $n\ge 1$ (so $G_1=G$).
Then 
$$G_n/G_{n+1}\simeq \underline{\Hom}(\Om^1,U(\LLie_+(\Om^1))_{n+1})$$ 
as a sheaf of (abelian) groups.
\end{thm}
 
First, we need to describe a certain version of de Rham complex for the free algebra (it is different
from the well known noncommutative forms considered in \cite{CQ}).

\begin{prop}\label{nc-dR-prop}
Let $V$ be a finite-dimensional vector space. Let us equip the algebra 
$$\AA_V={\bigwedge}^\bullet(V)\ot T(V)$$
with a dg-algebra structure, where
the grading is given by the natural grading of the exterior algebra (so $T(V)$ is in degree $0$)
and the differential $\tau$ is the unique derivation, trivial on $\bigwedge^\bullet(V)$ and
satisfying $\tau(1\ot v)=v\ot 1$.
Then we have $H^i(\AA_V)=0$ for $i>0$ and
$$H^0(\AA_V)=U(\LLie_+(V))\sub T(V).$$
\end{prop}

\Pf . It is easy to check that for any $v\in V$ and $x\in T(V)$ one has
$$\tau(vx-xv)=v\tau(x)-\tau(x)v.$$
This implies that $\tau(\LLie_+(V))=0$ and hence $U(\LLie_+(V))\sub \ker(\tau)$.
Let $e_1,\ldots,e_n$ be a basis of $V$. 
By the PBW-theorem, we can write any element of $T(V)$ uniquely as a finite sum
\begin{equation}\label{x-M-la-eq}
x=\sum_{\la:\BB\to\Z_{\ge 0}} [[f_{\la}(e)]] M_\la
\end{equation}
with $f_\la\in S(V)$, where we use the notation of Section \ref{NC-affine-space-sec}.
Since $M_\la$ are elements in $U(\LLie_+(V))$, for $x$ given by \eqref{x-M-la-eq} we have
$$\tau(x)=\sum_{i=1}^n e_i\ot \sum_{\la}[[\pa_if_\la(e)]] M_\la.$$
Let $\ov{M}_\la$ be the monomials $M_\la$ viewed as elements of $S(\LLie_+(V))$.
Then we have a natural isomorphism
$$\Om^\bullet_{S(V)}\ot_k S(\LLie_+(V))\to \AA_V: f\cdot de_{i_1}\we\ldots\we de_{i_r}\ot \ov{M}_\la\mapsto
de_{i_1}\we\ldots\we de_{i_r}\ot [[f(e)]] M_\la,$$
where $f\in S(V)$. Under this isomorphism the differential $\tau$ corresponds to the differential $d_{dR}\ot\id$ on
$\Om^\bullet_{S(V)}\ot_k S(\LLie_+(V))$, where $d_{dR}$ is the de Rham differential.
This immediately implies our assertion, since the $k$-linear span of the monomials $M_\la$ is contained
in the subalgebra $U(\LLie_+(V))\sub T(V)$.
\ed 

\begin{cor}\label{homotopy-cor} 
In the above situation there exists a $\GL(V)$-equivariant homotopy operators
$h=(h_n)$, where 
$$h_n: {\bigwedge}^n(V)\ot T^i(V)\to {\bigwedge}^{n-1}(V)\ot T^{i+1}(V)$$
satisfying $h_{n+1}\tau+\tau h_n=\id$ on $\bigwedge^n(V)\ot T(V)$ for $n\ge 1$.
Hence, for any vector bundle $\VV$ on a scheme $X$ there exist similar operators
on $\bigwedge^\bullet(\VV)\ot T(\VV)$.
\end{cor}

\Pf . This immediately follows from the vanishing of $H^{>0}(\AA_V)$ and from the reductivity of $GL(V)$.
\ed
 
We will need the following simple calculation with graded Lie commutators.

\begin{lem}\label{graded-Lie-lem} 
Let $\bL$ be the graded Lie algebra over $\Z[\frac{1}{2}]$
with generators $D_1,\ldots,D_n$, $F_1,\ldots,F_n$
of degree $1$ and relations
$$\sum_{i=0}^m [D_i,D_{m-i}]=\sum_{i=0}^m [D_i,F_{m-i}]=0$$
for $m=0,\ldots.n$.
Then one has
$$\sum_{i=1}^n [D_0,[D_i,D_{n+1-i}]]=\sum_{i=1}^n [D_0,[D_i,F_{n+1-i}]]=0$$
in $\bL$.
\end{lem}

\Pf . Note that $\bL$ has a Lie endomorphism, identical on $D_i$ and sending $F_i$ to $D_i$.
Thus, it is enough to prove the identity
\begin{equation}\label{formal-Lie-identity}
\sum_{i=1}^n [D_0,[D_i,F_{n+1-i}]]=0.
\end{equation}
Applying Jacobi identity we can rewrite the left-hand side as follows:
$$\sum_{i=1}^n [D_0,[D_i,F_{n+1-i}]]=
\sum_{i=1}^n [[D_0,D_i],F_{n+1-i}]+
\sum_{i=1}^n [[D_0,F_{n+1-i}],D_i].$$
Next, applying the relations in $\bL$ we get
$$\sum_{i=1}^n[[D_0,D_i],F_{n+1-i}]=-\frac{1}{2}\sum_{i\ge 1,j\ge 1}[[D_i,D_j],F_{n+1-i-j}],$$
$$\sum_{i=1}^n [[D_0,F_{n+1-i}],D_i]=-\sum_{i\ge 1,j\ge 1}[[D_j,F_{n+1-i-j}],D_i]$$
with the convention $F_m=0$ for $m<0$.
Now applying the Jacobi identity we get
$$\sum_{i\ge 1,j\ge 1}[[D_i,D_j],F_{n+1-i-j}]=2\sum_{i\ge 1,j\ge 1}[D_i,[D_j,F_{n+1-i-j}]]$$
and the assertion follows.
\ed
 
\noindent
{\it Proof of Theorem \ref{dg-constr-thm}.}
(i) Let us denote by $D_1$ the unique derivation of $\AA_X$ extending the de Rham differential on 
$\Om^\bullet$ and such that $D_1(1\ot \a)=\nabla(\a)\in \Om^1\ot\Om^1$ for $\a\in\Om^1$ (the existence of
$D_1$ follows from the Leibnitz identity for $\nabla$). We also set
$D_0=\tau$, which is the differential on $\AA_X$, zero on $\Om^\bullet$, and satisfying
$\tau(1\ot\a)=\a\ot 1$ for $\a\in\Om^1$ (see Proposition \ref{nc-dR-prop}). 
Note that $D_0^2=0$ and the condition that our connection is
torsion-free is equivalent  (cf.\ the proof of Lemma \ref{torsion-free-lem}) to the identity
$$[D_0,D_1]=D_0D_1+D_1D_0=0,$$
where we use $[\cdot,\cdot]$ to denote the supercommutator.
This implies that $[D_0,D_1^2]=0$, hence
$$D_0(D_1^2(1\ot\a))=0$$
for any $\a\in\Om^1$. Note that the derivation $D_1^2$ is zero on $\Om^\bullet$, so it is $\OO$-linear.
By Corollary \ref{homotopy-cor}, this implies that the $\OO$-linear map 
$$\nabla_2(\a):=-h_2(D_1^2(1\ot\a))$$
satisfies $D_0(\nabla_2(\a))=-D_1^2(1\ot\a)$.
Extending $\nabla_2$ to a degree $1$ derivation $D_2$ of $\AA_X$, zero on $\Om^\bullet$, this can be rewritten
as 
$$[D_0,D_2]+D_1^2=0.$$
Assume that $\OO$-linear 
maps $\nabla_i:\Om^1\to \Om^1\ot T^i(\Om^1)$ are defined for $2\le i\le n$ and that the corresponding degree $1$ derivations
$D_i$ of $\AA_X$ (zero on $\Om^\bullet$) satisfy
\begin{equation}\label{D-m-eq}
\sum_{i=0}^m [D_i,D_{m-i}]=0
\end{equation}
for all $m\le n$. We are going to construct an $\OO$-linear map
$\nabla_{n+1}$ such that the corresponding degree $1$ derivation
$D_{n+1}$ satisfies the equation \eqref{D-m-eq} for $m=n+1$, which can be rewritten as
\begin{equation}\label{D-n+1-eq}
[D_0,D_{n+1}]=-\frac{1}{2}\sum_{i=1}^n [D_i,D_{n+1-i}].
\end{equation}
By Lemma \ref{graded-Lie-lem}, we have
$$\sum_{i=1}^n [D_0,[D_i,D_{n+1-i}]]=0.$$
Note that the derivation $\sum_{i=1}^n[D_i,D_{n+1-i}]$ acts by zero on $\Om^\bullet$, so it is $\OO$-liinear.
Hence, by Corollary \ref{homotopy-cor}, we can set
$$\nabla_{n+1}(\a)=-\frac{1}{2}\sum_{i=1}^n h_{n+1}([D_i,D_{n+1-i}](1\ot\a)).$$

The above recursive procedure gives a degree $1$ derivation $D=D_0+D_1+D_2+D_3+\ldots$
of $\AA_X$ and the equations \eqref{D-m-eq} mean that $D^2=0$.

\noindent
(ii) Let us write $D=D_0+D_1+D_2+\ldots$, $D'=D'_0+D'_1+D'_2+\ldots$,  where $D_0=D'_0=\tau$,
$D_1$ and $D'_1$ correspond to some torsion-free connections, and $D_i$ (resp., $D'_i$)
are derivations extending the $\OO$-linear operators $\nabla_i$ (resp., $\nabla'_i$).
We claim that if $D_i=D'_i$ for $i<n$ then there exists an automorphism $\phi$ of $\hat{T}_\OO(\Om^1)$
with
$$\phi(\a)=\a+\phi_{n+1}(\a),$$
where $\phi_{n+1}(\a)\in T^{n+1}(\Om^1)$, such that
$$(\id\ot\phi) D(\id\ot\phi)^{-1})_i=D'_i \text{ for } i\le n.$$
Indeed, for any $\phi$ as above we have $(\id\ot\phi) D (\id\ot\phi)^{-1})_i=D_i$ for $i<n$, while
$$(\id\ot\phi) D (\id\ot\phi)^{-1})_n(1\ot \a)=D_n(1\ot\a)-D_0(1\ot \phi_{n+1}(\a)).$$
Thus, we need to find $\phi_{n+1}$ such that
$$D_n(1\ot \a)-D'_n(1\ot\a)=D_0(1\ot\phi_{n+1}(\a)).$$
By Corollary \ref{homotopy-cor}, it is enough to check that
$$D_0(D_n-D'_n)(1\ot\a)=0.$$
Since $D_i=D'_i$ for $i<n$ the equations \eqref{D-m-eq} for $D_i$ and $D'_i$ imply that
$$[D_0,D_n-D'_n]=0.$$
Hence, 
$$D_0(D_n-D'_n)(1\ot\a)=-(D_n-D'_n)(\a\ot 1)=0,$$
since $D_n$ and $D'_n$ have the same restriction to $\Om^\bullet$ (zero if $n\neq 1$ and the
de Rham differential for $n=1$).

\noindent
(iii) The condition that $\phi\in G_{n-1}$ preserves $D$ implies (by looking at the components
in $\Om^1\ot T^{n-1}(\Om^1)$) that $\tau\phi_n=0$, i.e., $\phi_n$ factors through
$U(\LLie_+(\Om^1))_n\sub T^n(\Om^1)$. Thus, we have a natural homomorphism 
$$G_{n-1}\to \underline{\Hom}(\Om^1,U(\LLie_+(\Om^1))_n): \phi\mapsto \phi_n$$
with the kernel $G_n$. It remains to prove the surjectivity of this map.
In other words, for any given $\phi_n:\Om^1\to T^n(\Om^1)$ such that $\tau\phi_n=0$, we have
to construct an element $\phi\in G_{n-1}$ such that $\phi(\a)=\a+\phi_n(\a)+\phi_{n+1}(\a)+\ldots$
for $\a\in\Om^1$. Let us first consider an automorphism $\wt{\phi}$ of $\hat{T}_\OO(\Om^1)$ defined by 
$\wt{\phi}(\a)=\a+\phi_n(\a)$. Then the differential
$\wt{D}=\wt{\phi}D\wt{\phi}^{-1}$ has components $(\wt{D})_i=D_i$ for $i<n$.
Now the proof of (ii) shows that there exists an automorphism $\psi$ of $\hat{T}_\OO(\Om^1)$
with $\psi(\a)=\a+\psi_{n+1}(\a)+\psi_{n+2}(\a)+\ldots$, such that $\wt{D}=\psi D\psi^{-1}$.
Then $\phi=\psi^{-1}\wt{\phi}$ is the required element of $G_{n-1}$.
\ed
 
\begin{rem} Note that $D_1^2(1\ot\a)=R(\a)\in \Om^2\ot \Om^1$, where $R$ is the curvature tensor.
Hence, $\nabla_2(1\ot\a)=h_2(R(\a))$.
It is easy to see that the identities $\tau(R(\a))=D_0(D_1^2(1\ot\a))=0$ and 
$\tau([D_1,D_2](1\ot\a))=\tau(\nabla(\nabla_2)(\a))=0$
correspond to the first and second Bianchi identities (in a particular case of a torsion-free connection).
\end{rem} 
 

\begin{lem}\label{D-homotopy-lem} Let $(C^\bullet,d)$ be a complex in an abelian category, such that each
$C^n$ is complete with respect to a decreasing filtration $(F^\bullet C^n)$.
Assume that $d=d_0+d_{\ge 1}$, where $d_0$ is also a differential on $C^\bullet$, i.e. $d_0^2=0$, and 
$$d_{\ge 1}(F^i C^n)\sub F^i C^{n+1}.$$
Let also $h$ be an operator of degree $-1$ on $C^\bullet$ such that 
$$h(F^i C^n)\sub F^{i+1} C^{n-1}.$$
Finally assume that the operator $p_0=\id-hd_0-d_0h$ satisfies $p_0 d_{\ge 1}=0$.
Then $\id+hd_{\ge 1}$ is an invertible operator on $C^\bullet$ and
$$d=(\id+hd_{\ge 1})^{-1}d_0(\id+hd_{\ge 1}).$$
\end{lem}

\Pf . Since $d_{\ge 1}$ preserves the filtration $F^\bullet C$ and $h$ raises the filtration index by $1$,
we immediately see that $\id+hd_{\ge 1}$ is invertible (since $C^n$ is complete with respect to the
filtration).
Thus, it is enough to check the identity
$$(\id+hd_{\ge 1})d=d_0(\id+hd_{\ge 1}).$$
We have
\begin{align*}
&(\id+hd_{\ge 1})d=d+hd_{\ge 1}d=d-hd_0d=d_0+d_{\ge 1}-hd_0d_{\ge 1}=\\
&d_0+(\id-hd_0)d_{\ge 1}=d_0+(p_0+d_0h)d_{\ge 1}=d_0+d_0hd_{\ge 1},
\end{align*}
as required.
\ed

\begin{prop}\label{D-homotopy-prop} Let $D=D_0+D_1+D_2+\ldots$ 
be a differential on $\AA_X$ as in Theorem \ref{dg-constr-thm}, and let $h=(h_n)$ be the homotopy
from Corollary \ref{homotopy-cor}. Set $D_{\ge 1}=D_1+D_2+\ldots$.
Then $\id+hD_{\ge 1}$ is an invertible operator on $\AA_X$ and 
$$D=(\id+hD_{\ge 1})^{-1}D_0(\id+hD_{\ge 1}).$$
Hence, the operator 
$$h_D=(\id+hD_{\ge 1})^{-1}h(\id+hD_{\ge 1}):\AA_X\to\Om^\bullet\ot\hat{T}^{\ge 1}(\Om^1)$$
of degree $-1$ satisfies $h_DD+Dh_D=\id$ on $\AA_X^i$ for $i>0$.
\end{prop}

\Pf . We apply Lemma \ref{D-homotopy-lem} to $d_0=D_0$, $d_{\ge 1}=D_{\ge 1}$,
the filtration $\Om^\bullet \ot \hat{T}^{\ge i}(\Om^1)$ on $\AA_X$, and the homotopy operator
$h$ from Corollary \ref{homotopy-cor}. Note that $p_0=\id-hD_0-D_0h$ satisfies
$p_0(\AA_X^{\ge 1})=0$ which implies that $p_0 D_{\ge 1}=0$.
\ed

\begin{cor}\label{coh-vanishing}  
Let $D$ be an algebraic NC-connection on a smooth scheme $X$.

\noindent
(i) The complex of sheaves of abelian groups $(\AA_X, D)$ is homotopy equivalent to its $0$th cohomology. 
In particular, $\underline{H}^i(\AA_X,D)=0$ for $i>0$. 

\noindent
(ii) We have an isomorphism of sheaves of abelian groups
\begin{equation}\label{ker-D-D0-eq}
(\id+hD_{\ge 1}):\ker(D)\rTo{\sim} \ker(D_0)\simeq \hat{U}(\LLie_+(\Om^1)),
\end{equation}
where $\hat{U}$ denotes the completion of the universal enveloping algebra
with respect to the augmentation ideal corresponding to the natural grading on
$\LLie_+(\Om^1)$.

\noindent
(iii) For any $d\ge 1$ consider the map 
$$D^{(d)}: T(\Om^1)/T^{\ge d+1}(\Om^1)\to \Om^1\ot (T(\Om^1)/T^{\ge d}(\Om^1))$$
induced by $D$. Then the natural map
$$K/(K\cap T^{\ge d+1}(\Om^1))\to \ker(D^{(d)}),$$
where $K=\ker(D:\AA^0_X\to \AA^1_X)$, is an isomorphism.
\end{cor}

\Pf . Parts (i) and (ii) follow directly from Proposition \ref{D-homotopy-prop}. For part (iii) we need
to check that for any $x\in \hat{T}(\Om^1)$ such that $D(x)\in \Om^1\ot T^{\ge d}(\Om^1)$, one
has $x\in K+T^{\ge d+1}(\Om^1)$. For this we use the operator $h_D$ from Proposition
\ref{D-homotopy-prop}. We have
$$D(h_DDx-x)=(Dh_D+h_DD)Dx-Dx=0,$$
so $h_DDx-x\in K$. It remains to note that the homotopy $h$ from Corollary \ref{homotopy-cor}
maps $\Om^\bullet\ot T^{\ge d}(\Om^1)$ to $\Om^\bullet\ot T^{\ge d+1}(\Om^1)$, while
$\id+hD_{\ge 1}$ preserves the filtration $(\Om^\bullet\ot T^{\ge i}(\Om^1))$.
Hence, 
$$h_DDx\in h_D(\Om^1\ot T^{\ge d}(\Om^1))\sub T^{\ge d+1}(\Om^1),$$
and the decomposition $x=(x-h_DDx)+h_DDx$ has the required properties.
\ed

Let us fix a smooth scheme $X$ with an algebraic NC-connection $D$. Consider the sheaf of rings
\begin{equation}\label{O-NC-def-eq}
\OO^{NC}=\und{H}^0(\AA_X,D)=\ker(D:\AA^0_X\to \AA^1_X).
\end{equation}
We have a natural projection $\pi:\AA^0_X=\hat{T}_\OO(\Om^1)\to \OO$ sending tensors of degree $\ge 1$ to zero,
and hence the induced homomorphism
$\pi:\OO^{NC}\to \OO$. Note that by Theorem \ref{dg-constr-thm}(ii),
the data $(\OO^{NC},\pi)$ does not depend on a choice of $D$ up to an
isomorphism.

\begin{prop}\label{section-prop} 
The homomorphism $\pi$ is surjective. Furthermore, it has a splitting 
$\si:\OO\to\OO^{NC}$ as a homomorphism of sheaves (not compatible with multiplication).
\end{prop}

\Pf . For $f\in\OO$ we have $D(f)=df$. By Proposition \ref{D-homotopy-prop}, 
$$\si(f)=f-h_D(df)\in\OO^{NC}$$ 
gives the required splitting.
\ed

We want to prove that $\OO^{NC}$ is a smooth NC-thickening of $\OO$. First, we need to check
that it is NC-complete. Consider the following
filtration by two-sided ideals on $\OO^{NC}$:
$$\OO^{NC}_{\ge d}:=\OO^{NC}\cap \hat{T}^{\ge d}_\OO(\Om^1) \text{ for } d\ge 1.$$

\begin{prop}\label{filtration-prop} 
(i) One has $\OO^{NC}_{\ge 1}=\OO^{NC}_{\ge 2}=\ker(\pi)$.

\noindent
(ii) One has $\OO^{NC}_{\ge n}=I_n(\OO^{NC})$, where the filtration $I_n(?)$ is
given by \eqref{I-n-def-eq}.

\noindent
(iii) The projection $\OO^{NC}_{\ge n}\to T^n(\Om^1)$ induces an isomorphism
$\OO^{NC}_{\ge n}/\OO^{NC}_{\ge n+1}\simeq U(\LLie_+(\Om^1))_n$.
\end{prop}

\Pf . (i) This follows from the fact that the component of the differential 
$$\Om^1\to \Om^1\ot \hat{T}(\Om^1)\to \Om^1$$
is injective.

\noindent
(ii) For any local section $a\in\OO^{NC}$ we have $\ad(a)\OO^{NC}\sub\OO^{NC}_{\ge 2}$
and $\ad(a)\OO^{NC}_{\ge n}\sub\OO^{NC}_{\ge n+1}$. 
Hence, $(\OO^{NC})_i^{\Lie}\sub\OO^{NC}_{\ge i}$. This implies the inclusion
$$I_n\OO^{NC}\sub \OO^{NC}_{\ge n}.$$
Since $\OO^{NC}_{\ge n}$ is complete with respect to the filtration by $\OO^{NC}_{\ge i}$ ($i>n$), 
it remains to check the surjectivity of the map
$$I_n(\OO^{NC})\to\OO^{NC}_{\ge n}/\OO^{NC}_{\ge n+1}.$$
It follows from Proposition \ref{nc-dR-prop} that
the latter quotient is a subsheaf of $U(\LLie_+(\Om^1))_n\sub T^n(\Om^1)$.
Thus, it is enough to check that $U(\LLie_+(\Om^1))_n$ is locally generated as an $\OO$-module by
elements of the form 
\begin{align}\label{gen-elements-eq}
&[\si(f_1),[\ldots,[\si(f_{i_1-1}),\si(f_{i_1})]\ldots]]\cdot 
[\si(f_{i_1+1}),[\ldots,[\si(f_{i_1+i_2-1}),\si(f_{i_1+i_2})]\ldots]]\cdot\ldots\cdot  \nonumber \\
&[\si(f_{i_1+\ldots+i_{k-1}+1}),[\ldots,[\si(f_{n-1}),\si(f_n)]\ldots]]\in I_n(\OO^{NC}),
\end{align}
where $f_1,\ldots,f_n\in\OO$ and $\si:\OO\to\OO^{NC}$ is the section defined in Proposition
\ref{section-prop}.
Since $\si(f)=f-df+\ldots$, for any $\a=\a_i+\a_{i+1}+\ldots\in \OO^{NC}_{\ge i}$ we have
$$[\si(f),\a]\equiv -df\ot\a_i+\a_i\ot df\mod\OO^{NC}_{\ge i+2}.$$
Hence, the expression \eqref{gen-elements-eq} is in $\OO^{NC}_{\ge n}$ with the
leading term 
\begin{align*}
&[df_1,[\ldots,[df_{i_1-1},df_{i_1}]\ldots]]\cdot
[df_{i_1+1},[\ldots,[df_{i_1+i_2-1},df_{i_1+i_2}]\ldots]]\cdot\ldots\cdot\\
&[df_{i_1+\ldots+i_{k-1}+1},[\ldots,[df_{n-1},df_n]\ldots]]\in U(\LLie_+(\Om^1))_n.
\end{align*}
By definition, such elements generate $U(\LLie_+(\Om^1))_n$.

\noindent
(iii) This follows from the proof in (ii).
\ed

\begin{cor} $\OO^{NC}$ is NC-complete.
\end{cor}

\begin{cor}\label{I-filtration-cor} 
The map \eqref{Lie-map-eq} is an isomoprhism for $\OO^{NC}$.
\end{cor}

\Pf . The proof of Proposition \ref{filtration-prop}(ii) shows that the composition of the
map \eqref{Lie-map-eq} with the isomorphism
$$\gr^n_I(\OO^{NC})\simeq \OO^{NC}_{\ge n}/\OO^{NC}_{\ge n+1}\rTo{\sim}U(\LLie_+(\Om^1))_n$$
is the identity.
\ed

Now we can prove that $\OO^{NC}$ is NC-smooth, thus finishing the proof of Theorem
\ref{main-nc-thm}.

\begin{thm}\label{smoothness-thm} The NC-thickening
$\OO^{NC}\to\OO$ defined by \eqref{O-NC-def-eq} is NC-smooth.
\end{thm}

\Pf . By Lemma \ref{smooth-local-lem}, we can assume $X$ to be affine, so there exists an NC-smooth
thickening $\OO^{univ}\to\OO$. Since $\OO^{NC}$ is NC-complete, by universality we get a homomorphism
of algebras $f:\OO^{univ}\to\OO^{NC}$ compatible with projections to $\OO$ 
(cf. \cite[(1.3.8), (1.6.2)]{Kapranov}). Since both algebras
are complete with respect to the $I$-filtration, to check that $f$ is an isomorphism 
it is enough to check that the induced maps
$$f_n: I_n\OO^{univ}/I_{n+1}\OO^{univ}\to I_n\OO^{NC}/I_{n+1}\OO^{NC}$$
are isomorphisms. By Lemma \ref{homomorphism-lem}, we have a commutative diagram
\begin{diagram}
U(\LLie_+(\Om^1))_n&\rTo{\id}&U(\LLie_+(\Om^1))_n\\
\dTo{}&&\dTo{}\\
\gr_I^n\OO^{univ}&\rTo{f_n}& \gr_I^n\OO^{NC}
\end{diagram}
in which vertical arrows are surjective. 
By Corollary \ref{I-filtration-cor}, the right vertical arrow is an isomorphism.
Hence, $f_n$ is an isomorphism as well.
\ed

\begin{cor}\label{I-filtration-gen-cor} 
Let $\wt{\OO}\to\OO_X$ be an NC-thickening of a smooth commutative scheme $X$. Then the map
\eqref{Lie-map-eq} is an isomorphism for $\wt{\OO}$ if and only if $(X,\wt{\OO})$ is NC-smooth.
\end{cor}

\Pf . ``If". The question is local, so we can assume that $X$ is affine and admits a torsion-free connection.
Then we have an isomorphism $\wt{\OO}\simeq\OO^{NC}$, since any two NC-smooth thickenings are
isomorphic in the affine case (see \cite[Thm.\ 1.6.1]{Kapranov}). Now we can use Corollary
\ref{I-filtration-cor}. 

\noindent ``Only if". By Lemma \ref{smooth-local-lem}, the question is local, so again we can assume that
$X$ is affine and admits a torsion-free connection. 
Let $\OO^{NC}\to\OO$ be the NC-smooth thickening constructed above. By universality
of the NC-smooth thickening in the affine case ((cf. \cite[(1.3.8), (1.6.2)]{Kapranov}) there exists
a homomorphism $f: \OO^{NC}\to\wt{\OO}$ of NC-thickenings. By assumption, $f$ induces
an isomorphism 
$$\gr^n_I\OO^{NC}\to \gr^n_I\wt{\OO}$$
for each $n$. Hence, $f$ is an isomorphism.
\ed

\begin{prop}\label{Aut-O-NC-prop} 
For the NC-thickening $\OO^{NC}\to \OO$ defined by \eqref{O-NC-def-eq}
let $\underline{\Aut}(\OO^{NC})$ denote the sheaf of groups of automorphisms of the extension
$\OO^{NC}\to\OO$. Let also be $G$, $G_n$ be the sheaves of groups defined in 
Theorem \ref{dg-constr-thm}(iii).
Then the natural maps 
$$G\to \underline{\Aut}(\OO^{NC}),$$
$$G/G_n\to \underline{\Aut}(\OO^{NC}/I_{n+1}\OO^{NC}),$$
where $n\ge 1$, are isomorphisms. 
\end{prop}

\Pf .
Recall that $G_n\sub G$ consists of automorphisms of $(\hat{T}_{\OO}(\Om^1),D)$ which induce the identity on
$\hat{T}(\Om^1)/\hat{T}^{\ge n+1}(\Om^1)$ (so $G=G_1$), and that we have natural identifications
\begin{equation}\label{G-quot-eq}
G_{n-1}/G_n\rTo{\sim}\Hom(\Om^1,U(\LLie_+(\Om^1)_n))
\end{equation} 
(see Theorem \ref{dg-constr-thm}(iii)).
Similarly, we can define for each $n\ge 1$ the subgroup $\und{\Aut}_n(\OO^{NC})\sub\und{\Aut}(\OO^{NC})$ consisting
of automorphisms that induce the identity on $\OO^{NC}/I_{n+1}\OO^{NC}$.
By Proposition \ref{filtration-prop}(ii), 
the homomorphism $G\to\und{\Aut}(\OO^{NC})$ sends $G_n$ to $\und{\Aut}_n(\OO^{NC})$.
We observe that we have a natural exact sequence of groups
$$1\to \und{\Aut}_{n-1}(\OO^{NC})\to\und{\Aut}(\OO^{NC})\to\und{\Aut}(\OO^{NC}/I_n\OO^{NC})\to 1.$$
Indeed, surjectivity here follows from the NC-smoothness of $\OO^{NC}$.
Note that 
$$\OO^{NC}/I_{n+1}\OO^{NC}\to\OO^{NC}/I_n\OO^{NC}$$ 
is a central extension
with the kernel $I_n\OO^{NC}/I_{n+1}\OO^{NC}$ which is generated by expressions 
\eqref{gen-elements-eq}.
Therefore, any automorphism of $\OO^{NC}/I_{n+1}\OO^{NC}$, identical modulo $I_n\OO^{NC}$,
acts on $I_n\OO^{NC}/I_{n+1}\OO^{NC}$ as the identity.
Hence, by \cite[Prop.\ 1.2.5]{Kapranov} we get an exact sequence
$$1\to \Hom(\Om^1,I_n\OO^{NC}/I_{n+1}\OO^{NC})\to
\und{\Aut}(\OO^{NC}/I_{n+1}\OO^{NC})\to\und{\Aut}(\OO^{NC}/I_n\OO^{NC})\to 1$$
(surjectivity again follows from the NC-smoothness of $\OO^{NC}$).
Hence, using Corollary \ref{I-filtration-cor} we obtain a natural identification 
\begin{equation}\label{Aut-quot-eq}
\und{\Aut}_{n-1}(\OO^{NC})/\und{\Aut}_n(\OO^{NC})\rTo{\sim} \Hom(\Om^1,U(\LLie_+(\Om^1)_n)).
\end{equation}
We claim that the homomorphism $G_{n-1}/G_n\to \und{\Aut}_{n-1}(\OO^{NC})/\und{\Aut}_n(\OO^{NC})$
is compatible (up to a sign) with the identifications \eqref{G-quot-eq} and \eqref{Aut-quot-eq} (and
hence, is an isomorphism). Note that this will imply that the homomorphism $G\to \und{\Aut}(\OO^{NC})$
is an isomorphism. Given an element $\phi\in G_{n-1}$ of $\hat{T}(\Om^1)$, for $\a\in\Om^1$ we have
$\phi(\a)=\a+\phi_n(\a)+\ldots$, where $\phi_n:\Om^1\to U(\LLie_+(\Om^1)_n)$.)
We have to restrict $\phi$ to $\OO^{NC}\sub\hat{T}(\Om^1)$ and see how it acts modulo
$I_{n+1}\OO^{NC}=\OO^{NC}_{\ge n+1}$. By definition, for a lifting $\si(f)=f-df+\ldots\in\OO^{NC}$ of $f\in\OO$
we have $\phi(\si(f))=\si(f)-\phi_n(df)\mod \OO^{NC}_{\ge n+1}$, which implies our claim.
\ed

Theorems \ref{dg-constr-thm} and \ref{smoothness-thm}
give a construction of an NC-smooth thickening $\OO^{NC}\to\OO$ starting from a torsion-free
connection on $X$. Such a connection always exists in the affine case as the following Lemma shows.

\begin{lem}\label{affine-torsion-free-lem} 
Let $X$ be an affine smooth scheme. Then there exists a torsion-free connection on $X$.
\end{lem}

\Pf . Locally we can always find a basis of $\Om^1$ consisting of closed forms (e.g., using local parameters),
hence a torsion-free connection. Covering $X$ with open subsets and choosing torsion-free connections
on them we will get a class in $H^1(X,\underline{\Hom}(\Om^1,S^2\Om^1))$ which is an obstacle to the
existence of a torsion-free connection on $X$. When $X$ is affine this class automatically vanishes.
\ed

Thus, we recover Kapranov's result on existence of NC-smooth thickenings for smooth affine varieties.
In the non-affine case the NC-smooth thickenings obtained from torsion-free connections are
rather special. The following result points out two simple restrictions on NC-thickenings that we get.

\begin{prop}\label{opposite-prop} Let $\OO^{NC}$ be an NC-smooth thickening
associated with an algebraic NC-connection via \eqref{O-NC-def-eq}.

\noindent
(i) The NC-thickening $\OO^{NC}\to \OO$ is isomorphic to its opposite.

\noindent
(ii) The $1$-smooth thickening $\OO^{NC}/F^2\OO^{NC}\to \OO$ is isomorphic to the standard one (see 
Definition \ref{standard-1-thickening})
\end{prop}

\Pf . (i) We have a natural anti-involution $\iota$ of $\AA_X$ which acts as the identity on the generators
$\Om^1\sub\Om^\bullet\sub\AA_X$ and $\Om^1\sub\hat{T}_\OO(\Om^1)$.
We have some differential $D=\tau+D_1+\ldots$ on $\AA_X$ such that $\OO^{NC}=\ker(D)$.
Now set $D'=\iota D\iota$. This is again a differential of the same type since
$\iota\tau\iota=\tau$. By Theorem \ref{dg-constr-thm}(ii), 
there is an isomorphism of algebras $\ker(D)\simeq\ker(D')$.
But $\iota$ induces also an isomorphism of $\ker(D')$ with the opposite algebra to $\OO^{NC}$.

\noindent
(ii) We have
$$\OO^{NC}/F^2\OO^{NC}=\OO^{NC}/I_3\OO^{NC}=\OO^{NC}/\OO^{NC}_{\ge 3}.$$
Set
$$\wt{\OO}=\{(f,x,y)\in \OO\oplus \Om^1\oplus T^2(\Om^1) \ |\ df+\tau(x)=0, \nabla(x)+\tau(y)=0\}.$$
Corollary \ref{coh-vanishing}(iii) implies that the natural projection
$$\OO^{NC}/\OO^{NC}_{\ge 3}\to \wt{\OO}$$
is an isomorphism. Note that the map 
$$\tau|_{S^2(\Om^1)}: S^2(\Om^1)\to \im(\tau)\sub \Om^1\ot\Om^1$$
is the multiplication by $2$. Thus, we have a natural section
$$\si:\OO\to\wt{\OO}: f\mapsto (f,-df,\frac{1}{2}\nabla(df)).$$
Now an easy computation shows that
$$\si(f)\si(g)=\si(fg)+\frac{1}{2}(df\ot dg-dg\ot df),$$
so $\wt{\OO}$ is isomorphic to the standard $1$-smooth thickening of $\OO_X$.
\ed

\begin{defi}\label{standard-thickening-def}
Let $X$ be a commutative algebraic group over a field of characteristic zero.
Then $X$ has a canonical torsion-free flat connection, such that flat sections are exactly translation-invariant vector fields (the vanishing of torsion follows from the fact that translation-invariant forms are closed).
We call the associated NC-smooth thickening of $X$ the {\it standard NC-thickening of $X$}. 
Note that since in this case the connection is flat, the corresponding NC-connection
has just two terms $D=D_0+D_1$, where $D_0$ is the $\OO$-linear derivation $\tau$ and $D_1$
is induced by the connection.
\end{defi}



\begin{prop}\label{NC-homomorphism-prop} 
Let $f:X\to Y$ be a homomorphism of commutative algebraic groups over a field of characteristic zero. Then it extends naturally to a morphism between the standard NC-smooth thickenings of $X$ and $Y$.
\end{prop}

\Pf . We have to extend the usual homomorphism $f^{-1}\OO_Y\to \OO_X$ to a homomorphism
of the standard NC-smooth thickenings $f^{-1}\OO_Y^{NC}\to \OO_X^{NC}$.
Such a homomorphism is obtained by passing to $0$th cohomology from the natural homomorphism
of dg-algebras
\begin{equation}\label{NC-homomorphism-eq}
f^{-1}\AA_Y^{NC}\to \AA_X^{NC}
\end{equation}
induced by the pull-back maps on functions and forms.
We only have to check that \eqref{NC-homomorphism-eq} is compatible with the 
differentials $D=D_0+D_1$. 
For $D_0$-parts this is clear, so it is enough to prove the commutativity of the diagram
\begin{diagram}
f^{-1}\Om^1_Y &\rTo{\nabla_Y}& f^{-1}(\Om^1_Y\ot_{\OO_Y}\Om^1_Y)\\
\dTo{}&&\dTo{}\\
\Om^1_X &\rTo{\nabla_X}&\Om^1_X\ot_{\OO_X}\Om^1_X
\end{diagram}
But this follows immediately from the fact that the pull-back of translation-invariant $1$-forms on $Y$
(which are $\nabla_Y$-horizontal) are translation-invariant $1$-forms on $X$.
\ed

\begin{prop}
Let $\OO^{NC}$ is the standard NC-thickening of a commutative algebraic group $X$.
Then for each $i$ we have 
$$H^i(X,\OO^{NC})\simeq H^i(X,U(\LLie_+(\Om^1)))\simeq H^i(X,\OO)\ot \hat{U}(\LLie_+(\Om^1|_e)),$$
where $\Om^1|_e$ is the fiber of $\Om^1$ at the neutral element $e\in X$.
In particular, if $X$ is an abelian variety then the algebra of global sections of $\OO^{NC}$ is 
isomorphic to $\hat{U}(\LLie_+(\Om^1|_e))$.
\end{prop}

\Pf . This follows from \eqref{ker-D-D0-eq}.
\ed

Next, we are going to show that all NC-smooth thickenings are obtained by a twisted version of our construction.

\begin{defi}\label{twisting-def}
For a smooth scheme $X$ we define the category $\CC_X$ as follows. An object of $\CC_X$,
called a \emph{twisted algebraic NC-connection}, consists of the data $(\TT,\JJ,\varphi,D)$, 
where $\TT$ is a sheaf of $\OO_X$-algebras equipped with a two-sided ideal $\JJ$,
$$\varphi:\bigoplus_{n\ge 0}\JJ^n/\JJ^{n+1}\to T(\Om^1)$$
is an isomorphism of graded algebras, and $D$ is a differential of degree $1$ on
$$\AA_\TT:=\Om^\bullet\ot_\OO \TT$$ 
such that $D$ restricts to the de Rham differential on $\Om^\bullet$,
and such that the composition
$$\JJ\hra\TT\rTo{D} \Om^1\ot_\OO \TT\to \Om^1\ot_\OO \TT/\JJ\simeq \Om^1$$
is equal to the composition of the projection $\JJ\to\JJ/\JJ^2$ with the isomorphism
$\varphi_1:\JJ/\JJ^2\to\Om^1$.
In addition, we require that $(\TT,\JJ)$ is locally isomorphic to 
$(\hat{T}_\OO(\Om^1),\hat{T}^{\ge 1}_\OO(\Om^1))$ (in particular,
$\TT$ is complete with respect to the filtration $(\JJ^n)$).
Morphisms between two such data are given by homomorphisms $f:\TT\to \TT'$ of $\OO$-algebras,
sending $\JJ\sub\TT$ to $\JJ'\sub\TT'$, induced the identity on $T(\Om^1)$ via
the isomorphisms $\varphi$ and $\varphi'$, and such that the induced homomorphism
$\id\ot f:\Om^\bullet\ot\OO \TT\to \Om^\bullet\ot_\OO \TT'$ is compatible with the differentials $D$.
\end{defi}

\begin{thm}\label{twisted-dg-thm} For any twisted algebraic NC-connection
one has $\und{H}^i(\AA_\TT,D)=0$ for $i>0$, and $\und{H}^0(\AA_\TT,D)$ is an NC-smooth
thickening of $X$ (with respect to its natural projection to $\TT/\JJ\simeq\OO$).
The functor $(\TT,\JJ,\varphi,D)\mapsto \und{H}^0(\AA_\TT)$ 
gives an equivalence of categories 
$$\CC_X\rTo{\sim} Th_X,$$ 
where $Th_X$ is the category of NC-smooth thickenings of $X$.
\end{thm}

\Pf . The first two assertion are local, so they follow from Corollary \ref{coh-vanishing}(i)
and Theorem \ref{smoothness-thm}.
Since $U\to \CC_U$ and $U\to Th_U$ are both gerbes (for $Th_U$ this follows from 
\cite[Thm.\ 1.6.1]{Kapranov}), it is enough to prove that
our functor induces an isomorphism on sheaves on automorphisms. Thus, for a data
$(\TT,\JJ,\varphi,D)$ and the corresponding NC-smooth thickening 
$$\OO^{NC}=\ker(D:\TT\to\Om^1\ot_\OO \TT)$$ 
we have to check that the natural map
$$G\to \underline{\Aut}(\OO^{NC})$$
is an isomorphism, where $G$ is the sheaf of groups of automorphisms of
$(\TT,\JJ,\varphi,D)$ and $\underline{\Aut}(\OO^{NC})$ is the sheaf of groups of automorphisms of the extension
$\OO^{NC}\to\OO$. This is a local question, so we can assume that
$\TT=\hat{T}(\Om^1)$, in which case $D$ is one of the differentials constructed in Theorem 
\ref{dg-constr-thm}. It remains to apply Proposition \ref{Aut-O-NC-prop}.
\ed

The following technical notion will be useful for us.

\begin{defi}\label{algebraic-d-def}
An NC-thickening $\OO_X^{(d)}\to\OO_X$ of a smooth variety $X$
is called \emph{$(d)_I$-smooth} 
if $I_{d+1}\OO_X^{(d)}=0$ and the map \eqref{Lie-map-eq} for
$\OO_X^{(d)}$ is an isomorphism in degrees $\le d$.
\end{defi}

There is also a suitable truncated notion of a twisted NC-connection.

\begin{defi} Let $d\ge 1$.
A \emph{$d$-truncated twisted algebraic NC-connection} is given by the data
$(\TT^{(d)},\JJ,\varphi,D^{(d)})$, 
where $\TT^{(d)}$ is a sheaf of $\OO_X$-algebras $\TT$ equipped with a two-sided ideal $\JJ$, such
that $\JJ^{d+1}=0$,
$$\varphi:\bigoplus_{n=0}^d\JJ^n/\JJ^{n+1}\to T(\Om^1)/T^{\ge d+1}(\Om^1)$$
is an isomorphism of graded algebras, and $D$ is a derivation
$$D^{(d)}:\TT^{(d)}\to\Om^1\ot_\OO (\TT^{(d)}/\JJ^d)$$
with values in a $\TT^{(d)}$-bimodule, 
such that $D^{(d)}(f)=df\ot 1$ for $f\in \OO$
and such that the composition
$$\JJ\hra\TT^{(d)}\rTo{D} \Om^1\ot_\OO (\TT^{(d)}/\JJ^d)\to \Om^1\ot_\OO (\TT/\JJ)\simeq \Om^1$$
is equal to the composition of the projection $\JJ\to\JJ/\JJ^2$ with the isomorphism
$\varphi_1:\JJ/\JJ^2\to\Om^1$.
In addition, we require that our data is locally isomorphic to the truncation of an algebraic NC-connection $D$,
i.e., $\hat{\TT}^{(d)}\simeq T(\Om^1)/T^{\ge d+1}(\Om^1)$, $\JJ$ is the image of $T^{\ge 1}(\Om^1)$
and $D^{(d)}$ is induced by $D$. 
Morphisms between two twisted NC-connections 
are given by homomorphisms $f:\TT^{(d)}_1\to \TT^{(d)}_2$ of $\OO$-algebras compatible with all the data.
\end{defi}

\begin{prop}\label{trunc-twisted-prop} For any $d$-truncated twisted algebraic NC-connection
the sheaf of rings $\ker(D^{(d)})$, with its natural projection to $\TT^{(d)_I}/\JJ\simeq\OO$, is
a $(d)_I$-smooth thickening of $X$. This gives an equivalence
\begin{equation}\label{trunc-thickenings-class-eq}
\CC^{(d)}_X\rTo{\sim} Th^{(d)_I}_X,
\end{equation}
from the category of $d$-truncated twisted algebraic NC-connections on $X$
with the category  $Th_X^{(d)_I}$ of $(d)_I$-smooth thickenings of $X$.
\end{prop}

\Pf . Note that both these categories are groupoids: for $Th_X^{(d)_I}$ this is deduced easily 
by induction in $d$ from the fact that the map \eqref{Lie-map-eq} for an $(d)_I$-smooth
NC-thickening $\OO_X^{(d)}$ is an isomorphism in degrees $\le d$. 
Furthermore, we claim that
$$U\mapsto \CC^{(d)}_U, \ Th^{(d)_I}_U$$
are gerbes on $X$. Indeed, for $\CC^{(d)}_U$ this follows from the definition.
For $(d)_I$-smooth NC-thickenings it is enough to check that every $(d)_I$-smooth 
thickening $\OO_U^{(d)})$
over a smooth affine $U$ is of the form $\OO_U^{NC}/I_{d+1}\OO_U^{NC}$
for the (unique) NC-smooth thickening $\OO_U^{NC}$ of $U$. By the universality of $\OO_U^{NC}$ 
(see \cite[Prop.\ 1.6.2]{Kapranov}) we can
construct a homomorphism $\OO_U^{NC}\to \OO_U^{(d)}$, identical on the abelianizations.
It factors through a homomorphism $\OO_U^{NC}/I_{d+1}\OO_U^{NC}\to \OO_U^{(d)}$.
The fact that it is an isomorphism follows from Corollary \ref{I-filtration-gen-cor} and the definition 
of $(d)_I$-smoothness. It remains to check that the functor
\eqref{trunc-thickenings-class-eq} induces an isomorphism of the local groups of isomorphisms.
To this end we recall that we have a natural isomorphism
$$G/G_d\rTo{\sim} \und{\Aut}(\OO_X^{NC}/I^{d+1}\OO_X),$$
where $G,G_d$ are sheaves of groups introduced in Theorem \ref{dg-constr-thm}(iii)
(see Proposition \ref{Aut-O-NC-prop}).
It remains to check that for the $d$-truncation $D^{(d)}$ of an algebraic NC-connection $D$
the natural map
$$G/G_d\to \und{\Aut}(T(\Om^1)/T^{\ge d+1}(\Om^1),\JJ,\varphi,D^{(d)})$$
is an isomorphism. The injectivity follows from the definition of $G_d$. The induction in $d$ reduces
the proof of surjectivity to checking that all automorphisms $\phi$ of $T(\Om^1)/T^{\ge d+1}(\Om^1)$
with $\phi(\a)=\a+\phi_d(\a)$ (where $\phi_d(\a)\in T^d(\Om^1)$), compatible with $D^{(d)}$,
are in the image of $G_{d-1}/G_d\sub G/G_d$ under the above map.
But one immediately checks that the compatibility with $D^{(d)}$ implies that $\tau\phi_d(\a)=0$
in $T^{d-1}(\Om^1)$, i.e., $\phi_d$ factors through a map $\Om^1\to U(\LLie_+(\Om^1))_d$,
so the assertion follows from Theorem \ref{dg-constr-thm}(iii).
\ed


\begin{rem}\label{filtrations-rem}
It is easy to see from the proof of the above Proposition that 
extensions of a fixed $(d)_I$-smooth thickening to a $(d+1)_I$-smooth thickening form
a gerbe banded by $\und{Hom}(\Om^1_X,(U\LLie_+(\Om^1_X))_{d+1})$
(see Lemma \ref{classification-lem} below). This is analogous to Kapranov's \cite[Prop.\ 4.3.1]{Kapranov}
which deals with $(d+1)$-smooth thickenings extending a given $d$-smooth thickening.
We'd like to point out that the results of sections (4.4)--(4.6) in \cite{Kapranov}
hold only if one replaces the commutator filtration with our $I$-filtration (see \eqref{I-n-def-eq}).
Indeed, for example, in \cite[Prop.\ 4.4.2]{Kapranov} one needs to use an embedding of 
the $d$-th associated graded quotient of the filtration on $\OO_X^{NC}$ into 
$(\Om^1_X)^{\ot d}$.
This does not work for the commutator filtration (e.g., the $1$st graded quotient is $\Om^2_X$),
however, we do have such embeddings for the $I$-filtration (see Lemma \ref{homomorphism-lem}).
In particular, \cite[Thm.\ 4.6.5]{Kapranov} holds with $d$-smooth replaced by $(d)_I$-smooth.
In \cite[Sec.\ 4.6]{Kapranov} second order thickenings should be replaced by $(3)_I$-smooth thickenings
($1$-smooth is the same as $(2)_I$-smooth).
\end{rem}

\section{Modules over NC-smooth thickenings}\label{module-sec}

\subsection{Modules from connections}
\label{module-conn-sec}

The construction of Theorem \ref{dg-constr-thm} has a natural analog for modules.

\begin{thm}\label{module-thm}
(i) Let $D=D_0+D_1+\ldots$ be an NC-connection on $X$. Let also
$(\FF,\nabla^{\FF})$ be a quasicoherent sheaf on $X$ with a connection (not necessarily integrable). 
Then there exists a natural
construction of a differential $D^\FF$ of degree $1$ on the $\AA_X$-module 
$$\FF\hat{\ot}_{\OO}\AA_X:=\liminv \FF\ot_\OO \Om^\bullet\ot_\OO (T(\Om^1)/T^{\ge n}(\Om^1)),$$
giving it
a structure of dg-module over $(\AA_X,D)$ such that for $s\in\FF$,
$$D^\FF(s\ot 1)=\nabla^\FF(s) \mod \FF\ot\Om^1\ot \hat{T}^{\ge 1}(\Om^1).$$
We refer to such a differential $D^\FF$ as a \emph{module NC-connection} (\emph{mNC-connection}).

\noindent
(ii) For any two mNC-connections $D^\FF$ and $\wt{D}^\FF$ on $\FF\hat{\ot}_{\OO}\AA_X$ 
there exists an $\AA_X$-linear automorphism $\Phi$ of $\FF\hat{\ot}_\OO\AA_X$ of degree $0$ 
inducing identity modulo $\FF\hat{\ot}_\OO\Om^\bullet\ot_\OO\hat{T}^{\ge 1}(\Om^1)$,
such that 
$$\wt{D}^\FF=\Phi D^\FF\Phi^{-1}.$$

\noindent
(iii) For a pair of sheaves with connection $(\FF,\nabla^\FF)$, $(\GG,\nabla^\GG)$, let
$D^\FF$ and $D^\GG$ be mNC-connections compatible with $\nabla^\FF$ and $\nabla^\GG$. 
Consider the sheaf 
$$\bH(\FF,\GG)=\bH((\FF,D^\FF),(\GG,D^\GG)):=
\Hom_{\AA_X}((\FF\hat{\ot}_\OO\AA_X,D^\FF),(\GG\hat{\ot}_\OO\AA_X,D^\GG))$$
consisting of homomorphisms $F$ of right $\AA_X$-modules of degree zero such that $D^\GG F=F D^\FF$
(i.e., of closed homomorphisms of dg-modules over $\AA_X$).
It has a natural decreasing filtration with
$$\bH(\FF,\GG)_n=\{F\in\bH(\FF,\GG) \ |\ F(s\ot 1)\in \GG\ot\hat{T}^{\ge n}(\Om^1) \text{ for }s\in\FF\}.$$
Then for each $n\ge 0$ there is a natural isomorphism
\begin{equation}\label{H-graded-pieces-eq}
\bH(\FF,\GG)_n/\bH(\FF,\GG)_{n+1}\simeq \Hom_\OO(\FF,\GG\ot U(\LLie_+(\Om^1))_n).
\end{equation}

\noindent
(iv) One has $\underline{H}^i(\FF\hat{\ot}_\OO\AA_X, D^\FF)=0$ for $i>0$.
\end{thm}

\Pf . (i) Let us set $\nabla^\FF_1=\nabla^\FF$. We are looking for a collection of $\OO$-linear maps 
$$\nabla^\FF_i:\FF\to \FF\ot\Om^1\ot T^{i-1}(\Om^1),$$ 
$i\ge 2$, such that the map 
$$D^\FF(s\ot 1)=\nabla^\FF_1(s)+\nabla^\FF_2(s)+\nabla^\FF_3(s)+\ldots,$$
extended to an endomorphism of $\FF\hat{\ot}_{\OO}\AA_X$ by the Leibnitz rule, satisfies $(D^\FF)^2=0$.
We can write this extension in the form
$$D^\FF=D^\FF_0+D^\FF_1+D^\FF_2+D^\FF_3+\ldots,$$
where $D^\FF_0=\id_\FF\ot D_0$ and 
for $i\ge 1$, $D^\FF_i$ are endomorphisms of $\FF\hat{\ot}_\OO\AA_X$ of degree $1$, 
extending $\nabla^\FF_i$ and satisfying the following versions
of the Leibnitz rule:
$$D^\FF_i(xa)=D^\FF_i(x)a+(-1)^{\deg(x)}x D_i(a).$$
Thus, the equation $(D^\FF)^2=0$ becomes a system of equations
\begin{equation}\label{V-D-eq}
[D^\FF_0,D^\FF_m]+\frac{1}{2}\sum_{i=1}^{m-1}[D^\FF_i,D^\FF_{m-i}]=0
\end{equation}
for $m\ge 1$. For example, the first equation is
$$[D^\FF_0, D^\FF_1]=0,$$
which follows from the fact that $[D^\FF_0, D^\FF_1]$ is $\AA_X$-linear (since $[D_0,D_1]=0$) and vanishes on the elements $s\ot 1$, where $s\in \FF$. 
Assume that $\nabla^\FF_i$ with $i\le n$ are already constructed (where $n\ge 1$), so that the equations
are \eqref{V-D-eq} are satisfied for $m\le n$. By Lemma \ref{graded-Lie-lem}, we have
$$\sum_{i=1}^n[D^\FF_0,[D^\FF_i, D^\FF_{n+1-i}]]=0.$$
Hence, $(D^\FF_0)\left(\sum_{i=1}^n [D^\FF_i, D^\FF_{n+1-i}](s\ot 1)\right)=0$, so setting
$$\nabla^\FF_{n+1}(s)=-\frac{1}{2}(\id_\FF\ot h_{n+1})\left(\sum_{i=1}^n [D^\FF_i, D^\FF_{n+1-i}](s\ot 1)\right),$$
where $(h_i)$ are the homotopy maps from Corollary \ref{homotopy-cor}, we will get that 
$$[D^\FF_0, D^\FF_{n+1}](s\ot 1)=(D^\FF_0)\nabla^\FF_{n+1}(s)=-\frac{1}{2}\left(\sum_{i=1}^n [D^\FF_i, D^\FF_{n+1-i}](s\ot 1)\right).$$
To see that this implies the equality
$$[D^\FF_0, D^\FF_{n+1}]=-\frac{1}{2}\sum_{i=1}^n [D^\FF_i, D^\FF_{n+1-i}]$$
we observe that the left-hand side satisfies the version of the Leibnitz rule with respect to the derivation
$[D_0,D_{n+1}]$ of $\AA_X$, while the right-hand side satisfies the similar rule with respect to
$-\frac{1}{2}\sum_{i=1}^n [D_i, D_{n+1-i}]$. It remains to use the equality \eqref{D-n+1-eq}.

\noindent
(ii) 
Assume we are given two mNC-connections $D^\FF$ and $\wt{D}^\FF$ as in (i)
with $D^\FF_i=\wt{D}^\FF_i$ for $i<n$ (where $D^\FF_0=\wt{D}^\FF_0=\id\ot D_0$). 
Then the equations \eqref{V-D-eq} imply that
$[\id\ot D_0,D^\FF_n-\wt{D}^\FF_n]=0$,
hence $(\id\ot D_0)(D^\FF_n-\wt{D}^\FF_n)(s\ot 1)=0$. Therefore, we can find
$\phi_n:\FF\to \FF\ot T^n(\Om^1)$ 
such that $(D^\FF_n-\wt{D}^\FF_n)(s\ot 1)=(\id\ot D_0)\phi_n(s)$ for $s\in \FF$. 
Consider the automorphism $\Phi$ of $\FF\hat{\ot}_\OO\AA_X$ as a right $\AA_X$-module given by
$s\ot 1\mapsto s\ot 1+\phi_n(s)$. Then it easy to check that
$(\Phi D^\FF\Phi^{-1})_i=\wt{D}^\FF_i$ for $i\le n$.
It remains to note that any infinite composition of automorphisms
$$\ldots\Phi_3\circ\Phi_2\circ\Phi_1$$ 
with $\Phi_i(s\ot 1)=s\ot 1\mod \FF\hat{\ot}_\OO \hat{T}^{\ge i}(\Om^1)$ converges
to an automorphism of $\FF\hat{\ot}_\OO\AA_X$.

\noindent
(iii) For $F\in\bH(\FF,\GG)_n$ we have $F(s\ot 1)=F_n(s)\mod \GG\hat{\ot} \hat{T}^{\ge n+1}(\Om^1)$
for some $\OO$-linear map $F_n:\FF\to \GG\ot T^n(\Om^1)$. It is easy to see that the condition
that $F$ commutes with the differential implies that $(\id\ot D_0)(F_n(s))=0$, so 
$F_n$ lands in $U(\LLie_+(\Om^1))_n\sub T^n(\Om^1)$ (see Proposition \ref{nc-dR-prop}). 
By definition $F\in \bH(\FF,\GG)_{n+1}$
if and only if $F_n=0$. Thus, it is enough to check that the map
$$\bH(\FF,\GG)_n\to \Hom(\FF,\GG\ot U(\LLie_+(\Om^1))_n): F\mapsto F_n$$
is surjective. In other words, we have to find $F_i:\FF\to \GG\ot T^i(\Om^1)$ with
$i\ge n+1$ such that $F(s\ot 1)=F_n(s)+F_{n+1}(s)+\ldots$ defines an
$\AA_X$-linear homomorphism in $\bH(\FF,\GG)$. Let us extend each $F_i$ to an $\AA_X$-linear homomorphism from $\FF\hat{\ot}_\OO\AA_X$ to $\GG\hat{\ot}_\OO\AA_X$.
Then we can write the equations on $(F_i)$ in the form
$$\sum_{i=0}^m[D_i,F_{m-i}]=0,$$
where $D_i$ appearing to the left (resp., to the right) of $F_j$ are $D^\GG_i$ (resp., $D^\FF_i$).
If $F_i$ with $i\le n$ are already found then we can construct $F_{n+1}$ using Lemma 
\ref{graded-Lie-lem}.

\noindent
(iv) The proof is entirely parallel to that of Proposition \ref{D-homotopy-prop}: we apply
Lemma \ref{D-homotopy-lem} to the decomposition $D^\FF=D^\FF_0+D^\FF_{\ge 1}$,
where $D^\FF_{\ge 1}=D^\FF_1+D^\FF_2+\ldots$ and use the homotopy operator
$\id_\FF\ot h$ for $D^\FF_0=\id_\FF\ot D_0$, where $h$ is the homotopy from Corollary
\ref{homotopy-cor}.
\ed

Theorem \ref{module-thm} gives a construction of right $\OO^{NC}$-modules: namely, starting with
any sheaf with connection we get such a module by taking 
$$\und{H}^0(\FF\hat{\ot}_\OO \AA_X,D^\FF)=\ker(D^\FF)\cap \FF\hat{\ot}_\OO\hat{T}(\Om^1)$$
equipped with a natural right action of $\OO^{NC}=\ker(D)\cap\hat{T}(\Om^1)$.
Note that $D^\FF$ is continuous with respect to the topology given by the filtration
$(\FF\ot_\OO\Om^\bullet\hat{\ot}_\OO\hat{T}^{\ge n}(\Om^1))$, so 
$\und{H}^0(\FF\hat{\ot}_\OO \AA_X,D^\FF)$ is complete with respect to the induced topology. 

\begin{cor}\label{module-cor}
(i) In the situation of Theorem \ref{module-thm}(i) set $\wt{\FF}=\und{H}^0(\FF\hat{\ot}_\OO \AA_X,D^\FF)$.
Then the natural map $\wt{\FF}\to\FF$
induced by the projection $\FF\hat{\ot}\hat{T}(\Om^1)\to\FF$, is surjective. Furthermore,
it induces an isomorphism
$$\wt{\FF}/\overline{\wt{\FF}\OO^{NC}_{\ge 1}}\rTo{\sim} \FF,$$
where the bar means taking the closure with respect to the natural topology on $\wt{\FF}$.

\noindent
(ii) Assume in addition that $\FF$ is a vector bundle of rank $r$.
Then $\wt{\FF}$ is locally free of rank $r$ as a
 right $\OO^{NC}$-module and 
 $$\wt{\FF}\ot_{\OO^{NC}}\OO\simeq \FF.$$

\noindent
(iii) In the situation of Theorem \ref{module-thm}(iii) assume that $\FF$ is a vector bundle.
Then the natural map
$$\bH((\FF,D^\FF),(\GG,D^\GG))\to \Hom_{\OO^{NC}}(\und{H}^0(\FF\hat{\ot}_\OO \AA_X,D^\FF),
\und{H}^0(\GG\hat{\ot}_\OO \AA_X,D^\GG))$$
is an isomorphism.
\end{cor}

\Pf . (i) Let us consider the natural decreasing filtration on $\wt{\FF}$ given by
$$\wt{\FF}_{\ge n}=\ker(D^\FF)\cap \FF\hat{\ot}_\OO\hat{T}^{\ge n}(\Om^1).$$
We have
\begin{equation}\label{K-H-eq}
\wt{\FF}=\underline{\bH}((\OO,D),(\FF,D^\FF)),
\end{equation}
where $\underline{\bH}(\cdot,\cdot)$ is the sheafified version of $\bH(\cdot,\cdot)$ (see Theorem
\ref{module-thm}(iii)).
Hence, from \eqref{H-graded-pieces-eq} we get natural isomorphisms
$$\wt{\FF}_{\ge n}/\wt{\FF}_{\ge n+1}\rTo{\sim} \FF\ot_\OO U(\LLie_+(\Om^1))_n.$$
In the case $n=0$ we get the required surjectivity.
Next, we observe that the natural map 
$$\wt{\FF}\ot_{\OO^{NC}}\OO^{NC}_{\ge 1}\to \wt{\FF}_{\ge 1}$$
is compatible with filtrations: it sends
$\wt{\FF}\ot_{\OO^{NC}}\OO^{NC}_{\ge n}$ to $\wt{\FF}_{\ge n}$.
Since the induced map on graded quotients
$$\wt{\FF}\ot_{\OO^{NC}}\OO^{NC}_{\ge n}/\OO^{NC}_{\ge n+1}\to
\wt{\FF}_{\ge n}/\wt{\FF}_{\ge n+1}$$
is surjective, we deduce the surjectivity of the map
$$\wt{\FF}\ot_{\OO^{NC}}\OO^{NC}_{\ge 1}\to \wt{\FF}_{\ge 1}/\wt{\FF}_{\ge n}$$
for $n\ge 1$.
In other words, we get that
$$\wt{\FF}_{\ge 1}=\wt{\FF}\ot_{\OO^{NC}}\OO^{NC}_{\ge 1}+\wt{\FF}_{\ge n}$$
for any $n\ge 1$, so $\wt{\FF}_{\ge 1}$ is the closure of $\wt{\FF}\ot_{\OO^{NC}}\OO^{NC}_{\ge 1}$.

\noindent
(ii) This follows from the fact that $\wt{\FF}\ot_{\OO^{NC}}\OO$ 
does not depend on a choice of connection (up to a canonical isomorphism). Hence, 
localizing we reduce to the case $\FF=\OO^{\oplus n}$.

\noindent
(iii) It is enough to prove that the same map is an isomorphism locally, 
so we can assume that $\FF=\OO$ with trivial connection. 
Then $\wt{\FF}=\OO^{NC}$, so the assertion follows from \eqref{K-H-eq}.
\ed


\subsection{The functor from the category of $D$-modules}\label{D-mod-sec}

Let $X$ be a smooth variety with an NC-smooth thickening $\OO_X^{NC}\to \OO_X$ obtained
from a twisted NC-connection $(\TT,\JJ,\varphi,D)$ as in Theorem \ref{twisted-dg-thm}, so that
$\OO_X^{NC}$ is the $0$th cohomology of the sheaf of dg-algebras
$$\AA_{\TT,D}:=(\AA_\TT=\Om^\bullet_X\ot\TT,D).$$

For any scheme $S$ over $k$ we consider the induced relative NC-thickening 
$$\OO_{S\times X/S}^{NC}\to \OO_{S\times X}$$ 
obtained as the tensor product $p_S^{-1}\OO_S\ot p_X^{-1}\OO_X^{NC}$,
where $p_S$ and $p_X$ are the projections of the product $S\times X$ onto its factors.
(thus, $\OO_S$ is a central subring in $\OO_{S\times X/S}^{NC}$).
Note that $\OO_{S\times X/S}$ has as dg-resolution the sheaf of dg-algebras
$p_X^*\AA_{\TT,D}$ on $S\times X$.

\begin{defi}\label{D-mod-defi}
Let $\FF$ be a quasicoherent sheaf on $S\times X$ equipped with a relative flat connection $\nabla$ in the
direction of $X$ (so $\nabla$ is $\OO_S$-linear). We associate with $\FF$ a right dg-module over
$p_X^*\AA_{\TT,D}$ by setting
\begin{equation}\label{K-F-nabla-eq}
K(\FF,\nabla):=(\FF\ot_\OO p_X^*\Om^\bullet)\hat{\ot}_{p_X^*\Om^\bullet} p_X^*\AA_{\TT,D},
\end{equation}
where $\FF\ot_\OO p_X^*\Om^\bullet$ is the de Rham complex of $\nabla$, which we view as a dg-module
over $p_X^*\Om^\bullet$. In other words,
$$K(\FF,\nabla)=\FF\hat{\ot}_\OO p_X^*\AA_{\TT,D}$$ 
as a right $p_X^*\AA_{\TT,D}$-module and the differential $D^\FF$ on it is uniquely determined by
\begin{equation}\label{flat-D-F-eq}
D^\FF(s\ot 1)=\nabla(s)\in \FF\ot_\OO \Om^1\sub K(\FF,\nabla),
\end{equation}
for $s\in\FF$. 
Since $\Om^\bullet$ is in the graded center of $\Om^\bullet\ot\TT$, we also have
a left action of $p_X^*\AA_{\TT,D}$ on $K(\FF,\nabla)$, so 
we can view it as a dg-bimodule over $p_X^*\AA_{\TT,D}$. Set 
$$\bK(\FF,\nabla)=\und{H}^0K(\FF,\nabla).$$
This is bimodule over $\OO^{NC}_{S\times X/S}$, so we get a functor
\begin{equation}\label{D-mod-functor-eq}
\bK=\bK_{S,X^{NC}}:D_{S\times X/S}-\mod_{qc}\to\OO^{NC}_{S\times X/S}-\bimod_{/S},
\end{equation}
where $\OO^{NC}_{S\times X/S}-\bimod_{/S}$ is the category of $\OO^{NC}_{S\times X/S}$-bimodules
with central action of $p_S^{-1}\OO_S$.
\end{defi}

\begin{rem}\label{D-module-rem} 
The results of Section \ref{module-conn-sec} (which are valid for NC-thickenings constructed
from global torsion-free connections) also have natural relative analogs.
For example, in Theorem \ref{module-thm} we can start with a quasicoherent sheaf $\FF$ on $S\times X$ with
a relative connection $\nabla$ (over $S$) and construct a dg-module over $p_X^*\AA_X$, which does not
depend on a connection up to a (noncanonical) isomorphism.
Note that in the untwisted case $\TT=\hat{T}(\Om^1)$ the dg-module \eqref{K-F-nabla-eq} obtained
from a flat relative connection $\nabla$ on $\FF$ fits into the framework
of the relative version of Theorem \ref{module-thm} (in the case of flat connection we don't need to
add higher corrections to the formula \eqref{flat-D-F-eq}).
\end{rem}

\begin{prop}\label{D-module-prop} Let $X$ be a smooth variety with a torsion free connection, 
and let $D$ be the corresponding NC-connection
giving rise to a relative NC-smooth thickening
$\OO^{NC}_{S\times X/S}=\und{H}^0(p_X^*\AA_X,D)$ over a scheme $S$.

\noindent
(i) Then for a quasicoherent sheaf $\FF$ on $S\times X$ with a relative connection $\nabla$ (over $S$) 
the dg-module $K(\FF,\nabla)$ and the $\OO^{NC}_{S\times X/S}$-module $\bK(\FF,\nabla)$
do not depend on a choice of a relative flat connection $\nabla$ on $\FF$
up to an isomorphism. 

\noindent
(ii) If $\FF$ is a vector bundle of rank $r$ on $S\times X$
with a relative flat connection then $\bK(\FF,\nabla)$ is locally free 
of rank $r$ as a right $\OO^{NC}_{S\times X/S}$-module.

\noindent
(iii) If $L$ is a line bundle with a flat connection $\nabla$ then 
$\bK(L,\nabla)$ is an invertible $\OO^{NC}$-bimodule. 
\end{prop}

\Pf . (i) By Remark \ref{D-module-rem} this follows from Theorem \ref{module-thm}(ii).

\noindent
(ii) This follows from the relative analog of Corollary \ref{module-cor}(ii).

\noindent
(iii) This can be checked locally using Corollary \ref{module-cor}(iii).
\ed

\begin{rem} The assertion of Proposition \ref{D-module-prop}(i) does not hold 
for the functor $\bK$ associated with an
arbitrary twisted NC-connection 
$(\TT,\JJ,\varphi,D)$: the module $\bK(\FF,\nabla)$ does depend on $\nabla$
in general (see Remark \ref{twisted-ab-var-conn} below).
\end{rem}

Note that for every $D$-module $(\FF,\nabla)$ we have a natural map 
\begin{equation}\label{bimodule-map-eq}
\bK(E,\nabla)\ot_{\OO^{NC}} \OO\to E
\end{equation}
induced by the projection $K(E,\nabla)\to E$. Furthermore, this map
is compatible with the $\OO^{NC}$-bimodule structure on $\bK(E,\nabla)$ and the
central $\OO$-bimodule structure on $E$ via the homomorphism $\OO^{NC}\to \OO$.

\begin{prop}\label{D-mod-prop} (i) The 
functor \eqref{D-mod-functor-eq} is exact.

\noindent (ii) If $(E,\nabla)$ is a vector bundle with a relative flat connection on $S\times X$ then the map
\eqref{bimodule-map-eq} is an isomorphism.
\end{prop}

\Pf . (i) Since by the relative version of
Theorem \ref{module-thm}(iv), $K(\FF,\nabla^\FF)$ has no higher cohomology,
it suffices to prove that the functor $\FF\mapsto \FF\hat{\ot}_\OO \TT$ is exact
(on quasicoherent sheaves). But this is clear, since over a sufficiently small open affine $U$
we can identify $\TT$ with $\hat{T}(\Om^1)$, and so
$$(\FF\hat{\ot}_\OO \TT)(U)\simeq \prod_n (\FF\ot_\OO T^n(\Om^1))(U).$$

\noindent
(ii) It suffices to check this locally, so we can use Corollary \ref{module-cor}(ii).
\ed


\begin{lem}\label{rel-conn-cd-lem} Let $f:S\to S'$ be a morphism of schemes of finite type over $k$.

\noindent
(i) Assume $f$ is an affine morphism. Then we have a commutative diagram
\begin{equation}\label{NC-push-forward-cd}
\begin{diagram}
D_{S\times X/S}-\mod_{qc} &\rTo{\bK_{S,X^{NC}}}& \OO_{S\times X/S}^{NC}-\bimod_{/S} \\
\dTo{(f\times\id)_*}&&\dTo{(f\times\id)_*}\\
D_{S'\times X/S'}-\mod_{qc} &\rTo{\bK_{S',X^{NC}}}& \OO_{S'\times X/S'}^{NC}-\bimod_{/S'} 
\end{diagram}
\end{equation}
where the vertical arrows are the natural push-forward functors.
If $S\to S'$ is an arbitrary morphism then there is a similar commutative diagram
for the derived categories and derived push-forward functors.

\noindent
(ii) We have a commutative diagram
\begin{equation}\label{NC-pull-back-cd}
\begin{diagram}
D_{S\times X/S}-\mod^{lff} &\rTo{\bK_{S,X^{NC}}}& \OO_{S\times X/S}^{NC}-\bimod_{/S}^{lff} \\
\uTo{f^*}&&\uTo{f^*}\\
D_{S'\times X/S'}-\mod^{lff}_{qc} &\rTo{\bK_{S',X^{NC}}}& \OO_{S'\times X/S'}^{NC}-\bimod_{/S'}^{lff} 
\end{diagram}
\end{equation}
where the superscript $lff$ means ``locally free of finite rank" (as right modules) and the vertical arrows are the pull-back functors given by 
$$f^*:\FF\mapsto p_S^{-1}\OO_S\ot_{p_S^{-1}f^{-1}\OO_{S'}}(f\times\id_X)^{-1}\FF.$$
\end{lem}

\Pf . (i) The case of affine morphism is straightforward. The case of general morphism
reduces to it by using Cech resolution to calculate the push-forwards.

\noindent
(ii) We have a natural isomorphism of
dg-bimodules over $p_X^*\AA_{\TT,D}$
\begin{equation}\label{K-F-F'-eq}
p_S^{-1}\OO_S\ot_{p_S^{-1}f^{-1}\OO_{S'}}(f\times\id_X)^{-1}K(\FF,\nabla)\rTo{\sim} K(\FF',f^*\nabla).
\end{equation}
Passing to the $0$th cohomology we obtain the required isomorphism of 
$\OO^{NC}_{S\times X/S}$-bimodules.
\ed

\subsection{Extendability of vector bundles to an NC-smooth thickening}\label{extendability-sec}

\begin{defi}\label{extension-defi}
Let $\wt{\OO}\to \OO_X$ be an NC-thickening, and let $E$ be a vector bundle on $X$. 
We say that a locally free right $\wt{\OO}$-module $\wt{E}$ is an {\it extension of $E$ to $\wt{\OO}$} if
$\wt{E}\ot_{\wt{\OO}}\OO\simeq E$. 
We say that $\wt{E}$ is a {\it bimodule extension of $E$}
if in addition $\wt{E}$ has a compatible left $\wt{\OO}$-module
structure making it an $\wt{\OO}$-bimodule (with $k$ in the center), and the map
$$\wt{\OO}\to\und{\End}(\wt{E})\to \und{\End}(E)$$
induced by the left action of $\wt{\OO}$ on $\wt{E}$ factors through the standard embedding
$\OO\to\und{\End}(E)$. 
\end{defi}

Given a vector bundle $E$ on $X$ and an NC-smooth thickening $\OO^{NC}\to \OO$
one can ask what are the obstructions for existence of an extension (resp., bimodule extension) 
of $E$ to $\OO^{NC}$. Our construction of $\OO^{NC}$-modules from connections gives
such extensions in the case of the NC-smooth thickening constructed from a torsion free connection on $X$
(see Corollary \ref{module-cor}(ii)). Also, it gives the following result for arbitrary NC-smooth
thickenings.

\begin{prop} Let $\wt{\OO}\to\OO$ be an NC-smooth thickening of $X$. Then for any vector bundle $E$ of rank
$n$ on $X$ admitting a flat connection there exists a
$\wt{\OO}$-bimodule $\wt{E}$, locally free of rank $n$ as a left (resp., right) $\wt{\OO}$-module
such that $\wt{E}\ot_{\wt{\OO}}\OO\simeq E$.
\end{prop}

\Pf . By Theorem \ref{twisted-dg-thm}, we can assume that $\wt{\OO}$ arises from a twisted NC-connection
$(\TT,\JJ,\varphi,D)$. Hence, the assertion follows from Proposition \ref{D-mod-prop}(ii).
\ed

Below we will give an obstruction for extending to $1$-smooth thickenings in terms
of Atiyah classes. Recall (see \cite[]{Kapranov})
that $1$-smooth thickenings $\wt{\OO}\to X$ form a gerbe banded
by $T\ot\Om^2$, where $T=T_X$, and that there is a standard $1$-smooth thickening (see
Definition \ref{standard-1-thickening}).
Thus, isomorphism classes of $1$-smooth thickenings are classified by $H^1(X,T\ot\Om^2)$.

Recall (see \cite{Atiyah})
that for a vector bundle $E$ on $X$ the Atiyah class $\a_E\in H^1(X,\Om^1\ot\underline{\End}(E))$
is the class of the extension
$$0\to E\ot\Om^1\to J^{\le 1}(E)\to E\to 0$$
where $J^{\le 1}(E)$ is the bundle of 1st jets in $E$. Equivalently, if $\nabla_i$ are local connections
on $E$ defined over an open covering $(U_i)$ then 
$(\nabla_j-\nabla_i\in \Ga(U_i\cap U_j,\Om^1\ot\underline{\End}(E)))$ is a \v Cech representative
of $\a_E$. 

\begin{prop}\label{obs-prop} 
Let $\wt{\OO}\to\OO$ be a $1$-smooth thickening of $X$ corresponding to the class
$c(\wt{\OO})\in H^1(X,T\ot\Om^2)$, and let $E$ be a vector bundle on $X$. 

\noindent
(i) $E$ extends to a right locally free $\wt{\OO}$-module if and only its Atiyah class $\a_E$ satisfies 
\begin{equation}\label{alpha-c-alpha-eq}
\a_E\cup \a_E-\lan c(\wt{\OO}),\a_E\ran=0
\end{equation}
in $H^2(X,\Om^2\ot\und{\End}(E))$.

\noindent
(ii) Let $\wt{\OO'}\to\OO$ be another $1$-smooth thickening of $X$.
Then $E$ extends to a $\wt{\OO'}-\wt{\OO}$-bimodule, locally free over $\wt{\OO}$,
and such that the induced $\OO-\OO$-bimodule structure on $E$ is standard, if and only if
\begin{equation}\label{kappa-alpha-c-c-eq}
2(\kappa\ot \id_{\und{\End}(E)})(\a_E)=\left(c(\wt{\OO})-c(\wt{\OO'})\right)\cdot \id_E
\end{equation}
in $H^1(X,T\ot\Om^2\ot\und{\End}(E))$,
where 
\begin{equation}\label{kappa-map-eq}
\kappa:\Om^1\to \und{\Hom}(\Om^1,\Om^2)\simeq T\ot\Om^2:x\mapsto (y\mapsto x\we y).
\end{equation}

\noindent
(iii) Assume $\dim X>1$. Then $E$ extends to an $\wt{\OO}$-bimodule if and only if $\a_E=0$.
\end{prop}

\Pf . (i) 
We have an exact sequence of sheaves of groups on $X$
$$0\to \Mat_n(\Om^2)\to \GL_n(\wt{\OO})\to \GL_n(\OO)\to 1$$
where $\Mat_n(\Om^2)$ embeds as an abelian normal subgroup in $\GL_n(\wt{\OO})$ via $M\mapsto 1+M$.
Let $(g_{ij}=\varphi_i^{-1}\varphi_j)$ be transition functions of $E$ with respect to some 
trivializations $\varphi_i:\OO_{U_i}^n\to E|_{U_i}$ over open covering $(U_i)$.
We can assume that $U_i$ are such that $\wt{\OO}|_{U_i}$ is isomorphic to the standard $1$-smooth
thickening of $U_i$, so over each $U_i$ we have a section
$\si_i:\OO_{U_i}\to\wt{\OO}_{U_i}$ such that
$$\si_i(f)\si_i(g)=\si_i(fg)+df\we dg.$$
Note that 
$$\si_j(f)-\si_i(f)=\lan v_{ij},df\ran,$$
for some Cech $1$-cocycle $(v_{ij})$ with values in $T\ot\Om^2$, such that $(v_{ij})$ represents the class
$c(\wt{\OO})$.
The bundle $E$ extends to $\wt{\OO}$ if and only if for some open covering we can lift some 
transition functions $(g_{ij})$ for $E$ to a Cech $1$-cocycle $(\wt{g}_{ij})$ with values
in $\GL_n(\wt{\OO})$.  Let us write
$$\wt{g}_{ij}=\si_i(g_{ij})+\om_{ij},$$
where $\om_{ij}\in\Mat_n(\Om^2)(U_{ij})$.
Then the $1$-cocycle condition for $\wt{g}_{ij}$ is
$$(\si_i(g_{ij})+\om_{ij})(\si_j(g_{jk})+\om_{jk})=\si_i(g_{ik})+\om_{ik}$$
over $U_{ijk}$.
Since $\si_j(g_{jk})=\si_i(g_{jk})+\lan v_{ij},dg_{jk}\ran$ and $g_{ij}g_{jk}=g_{ik}$, this is
equivalent to
$$dg_{ij}\we dg_{jk}+g_{ij}\lan v_{ij},dg_{jk}\ran+\om_{ij}g_{jk}+g_{ij}\om_{jk}=\om_{ik}.$$
Multiplying with $\varphi_i$ on the left and $\varphi_k^{-1}$ on the right we obtain the
following equation in $\Om^2\ot\und{\End}(E)|_{U_{ij}}$:
\begin{equation}\label{dg-cocycle-eq}
\varphi_i dg_{ij}\we dg_{jk} \varphi_k^{-1}+\lan v_{ij},\varphi_j dg_{jk}\varphi_k^{-1}\ran=
\eta_{ik}-\eta_{ij}-\eta_{jk},
\end{equation}
where $\eta_{ij}=\varphi_i\om_{ij}\varphi_j^{-1}$.
On the other hand, setting 
$$\nabla_i=\varphi_i\circ d\circ \varphi_i^{-1},$$
where $d$ is the natural connection on $\OO^n$, we get connections on $E|_{U_i}$.
Thus, the Atiyah class $\a_E$ is represented by the Cech $1$-cocycle
$$\a_{ij}=\nabla_j-\nabla_i=\varphi_i[g_{ij}\circ d\circ g_{ij}^{-1}-d]\varphi_i^{-1}=
-\varphi_i(dg_{ij})g_{ij}^{-1}\varphi_i^{-1}=-\varphi_i dg_{ij} \varphi_j^{-1}$$
Hence, the class $\a_E\cup \a_E$ is represented by the Cech $2$-cocycle
$$\a_{ij}\we\a_{jk}=\varphi_i dg_{ij} \varphi_j^{-1}\we\varphi_j dg_{jk} g_{jk}^{-1}\varphi_k^{-1}
=\varphi_i dg_{ij}\we dg_{jk} \varphi_k^{-1},$$
so the left-hand side of \eqref{dg-cocycle-eq} represents the class
$\a_E\cup \a_E-\lan c(\wt{\OO}),\a_E\ran$.

\noindent
(ii) We use the notation of part (i). Assume that $E$ is lifted to a right $\wt{\OO}$-module,
as in part (i).
The additional structure that makes an extension of $E$ into
a bimodule should be given by a collection of homomorphisms of sheaves of algebras
$$A_i:\wt{\OO'}|_{U_i}\to \Mat_n(\wt{\OO})|_{U_i}$$
such that
\begin{equation}\label{A-ij-bimodule-str-eq}
A_i(\wt{f})\wt{g}_{ij}=\wt{g}_{ij}A_j(\wt{f})
\end{equation}
for $\wt{f}\in \wt{\OO'}_{U_{ij}}$. 
The compatibility with the central bimodule structure on $E$ in Definition \ref{extension-defi}
is that $A_i$ is compatible with the standard diagonal embedding 
$\OO_{U_i}\to\Mat_n(\OO_{U_i}):f \mapsto f I$.
This implies (by considering commutators) that the restriction of $A_i$ induces the similar diagonal embedding 
$\Om^2\to\Mat_n(\Om^2)$. 
We can assume that the $1$-smooth thickenings $\wt{\OO'}|_{U_i}$ are standard, so there exist
sections $\si'_i:\OO\to\wt{\OO'}|_{U_i}$, such that 
$$\si'_i(f)\si'_i(g)=\si'_i(fg)+df\we dg,$$
$$\si'_j(f)-\si'_i(f)=\lan v'_{ij},df\ran,$$
where $(v'_{ij})$ is the $1$-cocycle with values in $T\ot\Om^2$ representing $c(\wt{\OO'})$.
Let us write for $f\in\OO_{U_i}$
$$A_i(\si'_i(f))=\si_i(f)I+\eta_i(f),$$
where $\eta_i(f)\in\Mat_n(\Om^2_{U_i})$. Then the condition that $A_i$ is a homomorphism is equivalent to
$\eta_i$ being a derivation with values in $\Mat_n(\Om^2)|_{U_i}$, so we can write 
$\eta_i(f)=\lan w_i,df\ran$, where $w_i\in T\ot\Mat_n(\Om^2)|_{U_i}$.
Note that 
$$A_j(\si'_i(f))=A_j(\si'_j(f)-\lan v'_{ij},df\ran)=\si_j(f)I+\eta_j(f)-\lan v'_{ij},df\ran I=
[\si_i(f)+\lan v_{ij}-v'_{ij},df\ran]\cdot I+\eta_j(f).$$
Thus, the equation \eqref{A-ij-bimodule-str-eq}
becomes
$$(\si_i(f)I+\eta_i(f))(\si_i(g_{ij})+\om_{ij})=(\si_i(g_{ij})+\om_{ij})([\si_i(f)+\lan v_{ij}-v'_{ij},df\ran]\cdot I+\eta_j(f))$$
which is equivalent to
$$2df\we dg_{ij}=g_{ij}\lan w_j+(v_{ij}-v'_{ij})I,df\ran-\lan w_i,df\ran g_{ij},$$
We can rewrite this as an equality in $\Om^2\ot \und{\End}(E)|_{U_{ij}}$:
$$2df\we \varphi_i dg_{ij}\varphi_j^{-1}=\lan u_j-u_i+(v_{ij}-v'_{ij})\id_E, df\ran,$$
where $u_i=\varphi_i w_i \varphi_i^{-1}$.
This equality means that the $1$-cocycles representing $2(\kappa\ot\id)(\a_E)$ and
$c(\wt{\OO})-c(\wt{\OO'})$ differ by the coboundary $(u_j-u_i)$, 
and our assertion follows.

\noindent
(iii) For $\wt{\OO'}=\wt{\OO}$ the equation \eqref{kappa-alpha-c-c-eq} becomes
$$2(\kappa\ot\id)(\a_E)=0.$$
It remains to observe that if
$\dim X>1$ then $\kappa$ is an embedding of a direct summand, so $\kappa\ot\id$ 
induces an embedding on $H^1$. Thus, the above equation is equivalent to
$\a_E=0$. 
\ed

\begin{rem} 
It is easy to see that $c(\wt{\OO}^{op})=-c(\wt{\OO})$, so if we are considering extension of $E$ to
left $\wt{\OO}$-modules then the equation \eqref{alpha-c-alpha-eq} should be replaced by
$$\a_E\cup \a_E+\lan c(\wt{\OO}),\a_E\ran=0.$$
This condition differs slightly from the condition in \cite[Rem.\ (4.6.9)]{Kapranov}. 
The equivalence of the two conditions
follows from the Bianchi identity for Atiyah classes proved in \cite[Prop.\ (1.2.2)]{Kapranov2}.
\end{rem}

\begin{cor}\label{van-cor}
Let $X$ be a smooth variety. 
If a vector bundle $E$ on $X$ extends to the standard $1$-smooth thickening of $X$
then $\ch_i(E)=0$ in $H^i(X,\Om^i)$ for $i>1$.
\end{cor}

\Pf . For the standard thickening $c(\wt{\OO})=0$, so the condition of Proposition \ref{obs-prop}(i)
becomes $\a_E\cup \a_E=0$. Now we
use the fact that the component $\ch_i(E)$ of the Chern character, up to a sign, is obtained from the 
Atiyah class $\a_E\in H^1(\Om^1\ot\underline{\End}(E))$ as the trace of
$\frac{1}{i!}\a_E^{\cup i}\in H^i((\Om^i\ot\underline{\End}(E))$  (see \cite[Sec.\ 5]{Atiyah}).
\ed

\begin{cor}\label{ample-bun-cor}
Let $X$ be a smooth proejctive variety of dimension $\ge 1$, and let $L$ be a line bundle on $X$,
such that either $L$ or $L^{-1}$ is ample. Let $\wt{\OO}\to\OO_X$ be a $1$-smooth thickening of
$X$ such that $c(\wt{\OO})=\kappa(\xi)$ for some $\xi\in H^1(X,\Om^1)$,
where $\kappa$ is the map \eqref{kappa-map-eq}.
Then $L$ extends to a right $\wt{\OO}$-module if and only if $\xi=c_1(L)$.
\end{cor}

\Pf . Since 
$$\lan\kappa(\a), c_1(L)\ran=\a\cup c_1(L)\in H^2(X,\Om^2),$$
the equation \eqref{alpha-c-alpha-eq} gives
$$(c_1(L)-\xi)\cup c_1(L)=0.$$
By Lefschetz theorem, this implies that $c_1(L)-\xi=0$.
\ed

Note that the assumption $c(\wt{\OO})=\kappa(\a)$ in the above Corollary is always satisfied for
surfaces, since in this case the map \eqref{kappa-map-eq} is an isomorphism.


\begin{ex} Kapranov constructed an NC-smooth thickening $\OO^{NC}_{\lleft}$
of the projective space $\P^n$, such that the line bundle $\OO(1)$ extends to a right $\OO^{NC}_{\lleft}$-module.
One can check that the map 
$$\kappa:H^1(\P^n,\Om^1)\to H^1(\P^n,T\ot\Om^2)$$
is an isomorphism.
Hence, by Corollary \ref{ample-bun-cor}, we obtain
$$c(\OO^{NC,\le 1}_{\lleft})=\kappa(c_1(\OO(1))).$$
\end{ex} 

In the case of line bundles on abelian varieties and the standard NC-smooth thickening
we can give a complete answer to the extension problem.

\begin{thm}\label{line-bundle-thm} Let $L$ be a line bundle on an abelian variety $A$, where $\dim A>1$.

\noindent
(i) $L$ extends to a (right) line bundle on the standard NC-smooth thickening $A^{NC}$
if and only if $c_1(L)^2=0$ in $H^2(A,\Om^2)$.

\noindent
(ii) $L$ extends to an invertible $\OO^{NC}_A$-bimodule
if and only if $c_1(L)=0$ in $H^1(A,\Om^1)$.
\end{thm}

\Pf . (i) The ``only if" part follows from Corollary \ref{van-cor}. Conversely, assume
$c_1(L)^2=0$. This means that the map $\eta: T_{A,0}\to H^1(A,\OO)$ given by $c_1(L)$
satisfies $\wedge^2\eta=0$, i.e., the image of $\eta$ is at most $1$-dimensional.
But $\eta$ is the tangent map to the homomorphism $\phi_L:A\to\hat{A}$ associated with $L$,
Hence, either $\phi_L=0$ or $\phi_L$ factors through an elliptic curve. In former case
$L$ is in $\Pic^0(A)$, so it extends to the NC-thickening. In the latter case we have a codimension $1$
abelian subvariety $B\sub A$ such that $L|_B$ is in $\Pic^0(B)$. Thus, for some $L_0\in\Pic^0(A)$
the restriction of $L\ot L_0^{-1}$ to $B$ is trivial. This implies that $L\ot L_0^{-1}=f^*M$,
where $f:A\to A/B$ is the corresponding morphism to an elliptic curve
and $M$ is a line bundle on $A/B$. Since any line bundle in $\Pic^0(A)$ admits a flat connection, it
lift to an invertible $\OO^{NC}$-bimodule. Thus, it is enough to
deal with pull-backs of line bundles from elliptic curves. In fact, pull-back of any vector bundle
from elliptic curve can be extended to the thickening because in this situation the commutative sheaf of functions on the elliptic curves embeds into $\OO^{NC}$ on the abelian variety, since all homomorphisms between
abelian varieties always extend to NC-thickenings---see Proposition \ref{NC-homomorphism-prop}.

\noindent
(ii) The ``if" part follows from Proposition \ref{D-mod-prop}(ii) since any line bundle in $\Pic^0(A)$
admits a flat connection. 
The ``only if" part follows from Proposition \ref{obs-prop}(ii).
\ed

\begin{cor}\label{nonstandard-NC-cor} 
Let $f:A\to E$ be a surjective morphism from an abelian variety to an elliptic curve,
and let $L_E$ be a line bundle on $E$ of nonzero degree.
Then $f^*L_E$ extends to a right line bundle $\LL$ on the standard NC-smooth thickening of $A$
and $\und{\End}(\LL)$ is an NC-smooth thickening of $A$, which is not isomorphic to the standard one.
\end{cor}

\Pf . The algebra $\und{\End}(\LL)$ is locally isomorphic to $\OO^{NC}_A$ (the standard NC-smooth thickening
 of $\OO_A$), so it is an NC-smooth thickening of $\OO_A$. If it were isomorphic to $\OO^{NC}_A$,
 we would get a $\OO^{NC}_A$-bimodule structure on $\LL$. But by Theorem \ref{line-bundle-thm}(ii),
 this is impossible since $c_1(f^*L_E)\neq 0$.
 \ed

\begin{rem} Let $\OO^{NC}\to \OO$ be an NC-smooth thickening constructed using a
torsion-free connection on $X$ (see Sec.\ \ref{NC-constr-sec}).
By Proposition \ref{obs-prop}(ii), if a vector bundle $E$ on $X$ admits a bimodule extension to $\OO^{NC}$ 
then its Atiyah class $\a_E$ vanishes. In the other direction, we know by Corollary \ref{module-cor}(ii), that if
$\a_E=0$ i.e., $E$ has a connection, then $E$ admits an extension to $\OO^{NC}$, but not necessarily
as a bimodule. Also, we know that bundles with flat connection admit bimodule extensions 
(see Proposition \ref{obs-prop}(ii)). It seems plausible that 
extendability of $E$ as a bimodule to the $2$-smooth extension 
$\OO^{NC}/F^3\OO^{NC}$ should imply the existence of a flat connection on $E$.
\end{rem}




\section{NC-Jacobian and the NC-Fourier-Mukai transform}\label{Jac-FM-sec}

\subsection{NC-thickening of a versal family of line bundles}\label{moduli-sec}

Let $Z$ be a complete variety, $L$ a line bundle on $Z\times B$, where 
$B$ is a smooth algebraic variety. Following \cite[Sec.\ 5.4]{Kapranov} let us make the following assumptions:

\noindent
(a) the Kodaira-Spencer map $\kappa:\TT_B\to H^1(Z,\OO)\ot\OO_B$ is an isomorphism,

\noindent
(b) $H^2(Z,\OO)=0$.

In addition let us fix a point $p_0\in Z$ and make one more assumption:

\noindent
(c) the line bundle $L|_{p_0\times B}$ is trivial. 

Let us fix a trivialization $s_0:\OO_B\to L|_{p_0\times B}$ and 
consider the functor $h^{rig}=h^{rig}_B$ on the category of NC-nilpotent algebras $\NN$, where
$h^{rig}(\La)$ is the set of isomorphism classes of the following data:

\noindent
(1) a line bundle $\LL$ on $Z\times\Spec(\La)$ together with a trivialization
$s:\OO\to \LL|_{p_0\times\Spec(\La)}$,

\noindent
(2) a morphism $f:\Spec(\La^0_{ab})\to B$, where $\La^0_{ab}$ is the quotient of $\La_{ab}$ by
the ideal of nilpotent elements,

\noindent
(3) an isomorphism $\varphi:\OO_{C\times\Spec(\La^0_{ab})}\otimes\LL\stackrel{\sim}{\to} (\id\times f)^*L$,
compatible with the trivializations over $p_0$.

Here by a line bundle on an NC-scheme we mean a sheaf of left $\OO$-modules, locally isomorphic
to $\OO$.

The following Lemma implies that an isomorphism of two data $(\LL,s)$ in (1) is unique 
(and also an isomorphism in (3) is unique). 

\begin{lem}\label{rig-lem} 
Let $\LL$ be a line bundle on $Z\times\Spec(\La)$.
Then any endomorphism of $\LL$ on $Z\times\Spec(\La)$ is uniquely determined
by its restriction to $p_0\times\Spec(\La)$.
\end{lem}

\Pf . 
Let $0=J_n\sub J_{n-1}\sub\ldots\sub J_1\sub J_0=\La$ be a filtration
by ideals such that $\La/J_1=\La^0_{ab}$ and $J_i/J_{i+1}$ are central bimodules over $\La/J_1$. 
Consider the induced filtration of $\LL$ by $J_i\LL$ and set $\LL^0=\LL/J_1\LL$ Note that 
the natural map 
$$J_i/J_{i+1}\ot_{\La^0_{ab}} \LL^0\to J_i\LL/J_{i+1}\LL$$
is an isomorphism (by local triviality of $\LL$).
Now suppose we have an endomorphism $f$ of $\LL$ on  $C\times\Spec(\La)$ which is zero over $p_0$.
Since $f$ is compatible with the filtration $(J_i\LL)$, it induces an endomorphism of
$\LL^0$. But $\End(\LL^0)=\La^0_{ab}$ and the condition that $f$ vanishes over $p_0$ implies
that $f$ induces zero endomorphism of $\LL^0$. Hence, $f$ factors through $J_1\LL\sub\LL$.
Next, assume $f$ factors through $J_i\LL$ for some $i\ge 1$ and consider the induced
map 
$$\LL^0=\LL/J_1\LL\to J_i\LL/J_{i+1}\LL\simeq (J_i/J_{i+1})\ot_{\La^0_{ab}} \LL^0.$$
As we have seen above such maps correspond to elements of 
$J_i/J_{i+1}$, so again the condition that restriction of $f$ to $p_0$ 
vanishes implies that $f$ factors through $J_{i+1}\LL$.
\ed

Let $h(\Lambda)$ be the set of isomorphism classes of the same data as for $h^{rig}$,
except for a trivialization over $p_0$ (this is the functor considered in \cite{Kapranov}). 
We have a natural map 
\begin{equation}\label{forget-rig-map}
h^{rig}(\Lambda)\to h(\Lambda).
\end{equation}

\begin{lem}\label{h-rig-h-lem} The map \eqref{forget-rig-map} is surjective. 
It is an isomorphism for commutative $\Lambda$.
\end{lem}

\Pf . First, let us check that any vector bundle over $\Spec(\La)$,
which becomes trivial on $\Spec(\La^0_{ab})$, is itself trivial. Since $\La$ is obtained
from $\La^0_{ab}$ by successive central extensions, it suffices to check that if a finitely generated
projective module $P$ over
$\La'$, where $\La'$ is a central extension of $\La$ by $I$, is such that $P/IP$ is a free $\La$-module 
then $P$ itself is free. Lifting generators of $P/IP$ to $P$ we get a morphism $(\La')^n\to P$
Using the condition $I^2=0$, we derive that it is an isomorphism.
Thus, if we have a line bundle $\LL$ over $Z\times\Spec(\La)$ together with an isomorphism
$\varphi:\LL|_{Z\times\Spec(\La^0_{ab})}\to (\id\times f)^*L$ then the restriction 
$\LL|_{p_0\times\Spec(\La)}$ is a trivial bundle over $\Spec(\La)$.
Let us choose its trivialization. It induces a trivialization of 
$\LL_{p_0\times\Spec(\La^0_{ab})}$, which differs from the one coming from that of $L|_{p_0\times B}$
by some invertible element of $\La^0_{ab}$. Since such an element can be lifted to an invertible element
of $\La$, we can choose a trivialization of $\LL|_{p_0\times\Spec(\La)}$ inducing a given trivialization
on $p_0\times\Spec(\La^0_{ab})$.

Now suppose that $\Lambda$ is commutative. Any two trivializations of $\LL|_{p_0\times \Spec(\La)}$
differ by an invertible element of $\La$. Since such an element gives an automorphism of $\LL$,
we see that different choices of trivialization lead to isomorphic data.
\ed

\begin{thm}\label{h-rig-repr-thm}
The functor $h^{rig}$ is represented by an NC-smooth thickening $B^{NC}$ of $B$.
Furthermore, there is a universal line bundle $L^{NC}$ on $Z\times B^{NC}$, trivialized over 
$p_0\times B^{NC}$ and extending the line bundle $L$ on $Z\times B$.
\end{thm}

\Pf . First, let us prove formal smoothness of the functors $h$ and $h^{rig}$.
Let 
$$0\to I\to\La'\to \La\to 0$$ 
be a central extension of NC-nilpotent algebras. 
Then we have an exact sequence 
$$1\to \II \to \OO_{Z\times\Spec(\La')}^*\to \OO_{Z\times\Spec(\La)}^*\to 1$$
of sheaves of groups on $Z\times\Spec(\La')$, where
we view $Z\times\Spec(\La)$ as a closed NC-subscheme of $Z\times\Spec(\La')$ with
the corresponding sheaf of ideals $\II\simeq\wt{I}\boxtimes\OO_Z$, where
$\wt{I}$ is the quasicoherent sheaf on $\Spec(\La_{ab})$ associated with the $\La_{ab}$-module $I$.
We have 
$$H^2(Z\times\Spec(\La_{ab}),\II)\simeq I\ot H^2(Z,\OO)=0$$ 
since $H^2(Z,\OO)=0$ by assumption. This implies surjectivity
of the map of nonabelian cohomology
$$H^1(Z,(\La'\ot\OO_Z)^*)\to H^1(Z,(\La\ot\OO_Z)^*),$$
hence the map $h(\La')\to h(\La)$ is surjective, i.e., $h$ is formally smooth. 
An element $x\in h^{rig}(\La)$ induces an
element $y\in h(\La)$ which can be lifted to $y'\in h(\La')$ and then to $x'\in h^{rig}(\La')$. However,
the element $x'$ does not necessarily lift $x$: the induced rigidification over $\La$ may differ from
the original one by some invertible element of $\La$. Lifting it to an invertible element of $\La'$
and changing $x'$ accordingly, we get a lifting of $x$ to $h^{rig}(\La')$.

The fact that on commutative algebras $h$ (and hence, by Lemma \ref{h-rig-h-lem}, $h^{rig}$) is isomorphic to $h_B$ was proved in \cite[Prop.\ 5.4.3(a)]{Kapranov}.

Finally, representability for $h^{rig}$ is proved along the lines of \cite[Sec.\ 5.4]{Kapranov} (with a slight modification, see Remark \ref{repr-rem} below). By \cite[Thm.\ 2.3.5]{Kapranov}, we have to check that for any
Cartesian diagram of algebras in $\NN$ 
\begin{diagram}
\La_{12} &\rTo{q_1}& \La_1\\
\dTo{q_2}&&\dTo{p_1}\\
\La_2 &\rTo{p_2}&\La
\end{diagram}
such that $p_1$ is a central extension, the natural map
\begin{equation}\label{NC-descent-map}
h^{rig}(\La_{12})\to h^{rig}(\La_1)\times_{h^{rig}(\La)} h(\La_2)
\end{equation}
is an isomorphism.
To this end we are going to construct the inverse map. Namely, assume we are
given for $i=1,2$ the data $(\LL_i, f_i,\varphi_i)\in h^{rig}(\La_i)$ that induce the same
element of $h^{rig}(\La)$. Let us view $\LL_i$, for $i=1,2$, as locally free modules
of rank $1$ over $\OO_Z\otimes_k\La_i$. Then by Lemma \ref{rig-lem}, we have a unique isomorphism 
$$\LL_1\ot_{\La_1}\La\simeq \LL_2\ot_{\La_2}\La$$
of $\OO_Z\otimes_k\La$-modules, compatible with trivializations over $p_0$.
Thus, setting $\LL=\LL_i\ot_{\La_i}\La$, we can define a
locally free $\OO_Z\otimes\La_{12}$-module of rank $1$
$$\LL_{12}:=\LL_1\times_{\LL} \LL_2.$$

\ed

\begin{rem}\label{repr-rem}
In \cite[Sec.\ 5.4]{Kapranov} it is claimed that the functor $h$ is representable. However,
the proof of \cite[Prop.\ 5.4.3]{Kapranov} has a gap: a choice in the gluing procedure makes
the proposed inverse map to \eqref{NC-descent-map} ill-defined. This problem does not
appear for $h^{rig}$ (due to Lemma \ref{rig-lem}).
Note that since $h$ is formally smooth and $h^{rig}$ and $h$ restrict to the same functor on
commutative schemes, it follows that 
$h$ is representable if and only if the map $h^{rig}\to h$ is an isomorphism.
This is equivalent to the statement that for an NC-nilpotent algebra $\La$, and 
a line bundle $\LL$ over $C\times\Spec(\La)$, such that the restriction $\LL|_{C\times\Spec(\La^0_{ab})}$
is induced via some map $f:\Spec(\La^0_{ab})\to B$ by our family of line bundles over $B$, the map
$$\Aut(\LL)\to \Aut(\LL_{p_0\times\Spec(\La)})$$
is an isomorphism (note that by Lemma \ref{rig-lem}, this map is injective).
In the case of commutative algebras this statement
follows immediately from the identification of endomorphisms of a line bundle over 
$C\times\Spec(\La)$ with $\La$. 
\end{rem}

\subsection{NC-thickening of the Jacobian}

Let $C$ be a smooth projective curve, $p_0\in C$ a point.
Let $\PP_C$ be the Poincar\'e line bundle on $C\times J$, where $J=J(C)$ is the Jacobian of $C$,
trivialized over $p_0\times J$. This family of line bundles over $J$ 
satisfies the assumptions (a),(b),(c) from Section \ref{moduli-sec},
so we have the corresponding functor $h^{rig}_J$ of NC-families on the category $\NN$ 
of NC-nilpotent algebras.

On the other hand, we 
define the functor $H^{rig}$ on NC-complete schemes by letting $H^{rig}(S)$ be the set of isomorphism
classes of pairs $(\LL,s)$, where $\LL$ is a line
bundle on $C\times S$ together with a trivialization
$s:\OO\to \LL|_{p_0\times S}$. 

\begin{lem}\label{h-rig-H-lem} 
The natural map of functors $h^{rig}_J\to H^{rig}|_{\NN}$ is an isomorphism.
\end{lem}

\Pf . This follows immediately from the universality of the Poincar\'e line bundle: for any NC-nilpotent
algebra $\La$ and $(\LL,s)\in H^{rig}(\Spec(\La))$ the line bundle
$\OO_{C\times\Spec(\La_{ab})}\ot\LL$ on $C\times\Spec(\La_{ab})$ corresponds to a morphism
$f:\Spec(\La_{ab})\to J$, which allows to construct an element of $h^{rig}_J(\Spec(\La))$.
\ed

\begin{cor}\label{Jac-mod-cor}
The functor $H^{rig}$ is represented by an NC-smooth thickening $J^{NC,mod}$ of $J$.
We denote by $\PP^{NC,mod}$ 
the corresponding universal family, which is a line bundle on $C\times J^{NC,mod}$ trivialized over 
$p_0\times J^{NC,mod}$.
\end{cor}

\Pf . 
First, by Theorem \ref{h-rig-repr-thm} and by Lemma \ref{h-rig-H-lem} we know the representability of
$H^{rig}|_\NN$, so
$$H^{rig}|_{\NN}\simeq h_{J^{NC,mod}}|_\NN$$
for some NC-smooth NC-thickening $J^{NC,mod}$ of $J$.
We claim that this isomorphism of functors extends uniqely to an isomorphism of functors
$H^{rig}$ and $h_{J^{NC,mod}}$ on the category of all NC-schemes.
Indeed, this easily reduces to the case of affine NC-schemes.
By Lemma \ref{NC-mor-lem}, it is enough to check that the natural map
$$H^{rig}(R)\to \liminv H^{rig} (R/F^nR)$$
is an isomorphism for any NC-complete algebra $R$. 
For this it suffices to prove that for any NC-scheme $X$ the category of line bundles on $X$
is equivalent to the category of collections $(L_n,f_n)$, where $L_n$ is a line bundle on $X^{\le n}$
and $f_n$ is an isomorphism $L_{n+1}|_{X^{\le n}}\simeq L_n$. Since
$\OO_X^*=\liminv \OO_{X^{\le n}}^*$, the natural restriction functor is fully faithful, so it suffices
to check that any collection $(L_n,f_n)$ is locally trivial. Indeed, assume that $X$ is affine
and $L_0\simeq\OO_X$. Note that for each $n$, we have an exact sequence 
$$0\to (F^{n+1}\OO_X)L_{n+1}\to L_{n+1}\to L_n\to 0$$
of sheaves on $X$, where 
$$(F^{n+1}\OO_X)L_{n+1}\simeq (F^{n+1}\OO_X/F^{n+2}\OO_X)\ot_{\OO_{X^{\le 0}}}L_0$$
is a quasicoherent sheaf of $\OO_{X^{\le 0}}$-modules. Hence, the map $H^0(X,L_{n+1})\to H^0(X,L_n)$
is surjective. Using this we can lift a trivialization $\OO_{X^{\le 0}}\rTo{\sim} L_0$ to a compatible
system of maps $\OO_{X^{\le n}}\to L_n$, which will automatically be isomorphisms.
\ed

\begin{defi} We call the NC-smooth thickening $J^{NC,mod}$ of the Jacobian $J$, constructed above,
the {\it modular NC-thickening of $J$}, and call the universal line bundle $\PP^{NC,mod}$ on
$C\times J^{NC,mod}$ the {\it NC-Poincar\'e bundle}.
\end{defi}

Below we will show that the modular NC-thickening of $J$
is isomorphic to the standard NC-thickening of $J$ as an abelian variety
(see Corollary \ref{J-identification-cor}).

\subsection{NC-Poincar\'e bundle and NC-Fourier-Mukai transform}\label{NC-Poincare-sec}

Recall that for an abelian variety $A$ there is a universal
extension $A^\natural\to \hat{A}$ of $\hat{A}$ by a vector space, namely, $A^\natural$ is the moduli space
of line bundles with flat connection on $A$. Then the line bundle $\PP^\natural$ on $A^\natural\times A$
obtained as the pull-back of $\PP$,
has a relative flat connection $\nabla^\natural$ (in the direction of $A$, i.e., relative over $A^\natural$). 
This is used in \cite{Laumon} to construct the Fourier-Mukai transform (in fact, an equivalence)
$$F^D: D(\Qcoh(A^\natural))\to D(D_A-\mod_{qc}):E\mapsto Rp_{2*}(\PP^\natural\ot p_1^*E).$$
Applying the above functor to $\PP^\natural$ we obtain a right 
$\OO^{NC}_{A^\natural\times A/A^\natural}$-module 
$\bK_{A^\natural, A}(\PP^\natural,\nabla^\natural)$, locally free of rank $1$.
Hence, we can define the corresponding Fourier-Mukai type functor
\begin{equation}\label{FNC-natural-eq}
F^{NC,\natural}:D(\Qcoh(A^\natural))\to D(\mod-\OO^{NC}_A): 
E\mapsto Rp_{2*}(\bK_{A^\natural,A}(\PP^\natural,\nabla^\natural)\ot_{p_1^{-1}\OO_{A^\natural}} p_1^{-1}E).
\end{equation}

\begin{prop}\label{F-NC-natural-prop}
One has a commutative triangle of functors
\begin{diagram} 
D(\Qcoh(A^\natural))&\rTo{F^D}& D(D_A-\mod_{qc})\\
&\rdTo{F^{NC,\natural}}&\dTo_{\bK_{pt,A}}\\
&&D(\mod-\OO^{NC}_A)
\end{diagram}
\end{prop}

\Pf . For $E\in D(\Qcoh(A^\natural))$ we consider 
$\PP^\natural\ot_{p_1^{-1}\OO_{A^\natural}} p_1^{-1}E$ equipped with the relative connection over $A^\natural$.
Then we can proceed in two ways: either apply the functor $\bK_{A^\natural,A}$ and
then take push-forward to $A$, or first take push-forward and then apply the functor
$\bK_{pt,A}$. The first procedure gives $F^{NC,\natural}(E)$, while the second gives
$\bK_{pt,A}(F^D(E))$. It remains to use Lemma \ref{rel-conn-cd-lem}(i)
for the projection $A^\natural\to pt$. 
\ed

\begin{rem}
If $U\sub \hat{A}$ is an open affine
then there exists a section $U\to A^\natural$, hence a 
flat connection on $\PP|_{U\times A}$, relative to $U$. Applying the functor $\bK_{U,A}$ we get a line bundle
on $U\times A^{NC}$. Thus, for an affine open covering $(U_i)$ of $\hat{A}$
we obtain line bundles 
 $\wt{\PP}_{U_i}$ on $U_i\times A^{NC}$ extending $\PP|_{U_i\times A}$.
On double intersections we will have isomorphisms of line bundles (by Proposition \ref{D-module-prop}(i)) 
but they will not be compatible on triple intersections.
Thus, we obtain a twisted line bundle on the NC-scheme
$\hat{A}\times A^{NC}$ (in particular, we get a canonical gerbe on this NC-scheme).
\end{rem}

\noindent
{\bf Construction}. Assume that we have a map $f:C\to\hat{A}$ where $C$ is a curve.
Pulling back the Poincar\'e line bundle 
from $\hat{A}\times A$ we obtain a line bundle $\PP(f)=(f\times\id)^*\PP$ on $C\times A$.
We can cover $C$ by two open affine
subsets $U_1$ and $U_2$, and choose liftings $\si_1:U_1\to A^\natural, \si_2:U_2\to A^\natural$.
Hence, we get the induced relative connections $\nabla_1, \nabla_2$ on 
$\PP(f)|{U_1\times A}$ and $\PP(f)|_{U_2\times A}$, respectively. By Proposition \ref{D-module-prop}(i),
the corresponding line bundles
over the relative NC-thickenings $U_1\times A^{NC}$ and $U_2\times A^{NC}$
are isomorphic over the intersection $U_1\cap U_2$. Hence, 
choosing such an isomorphism $\varphi$ and gluing them we obtain
a line bundle $\PP^{NC}_{\si_1,\si_2,\varphi}$ on the relative NC-thickening $C\times A^{NC}$,
which extends the line bundle $\PP(f)$. 

Let us assume in addition that $f(p_0)=0$ for some fixed point $p_0\in C$.
We can choose $U_1\sub C$ to be an open neighborhood of $p_0$ and require that
a lifting $\si_1:U_1\to A^\natural$ satisfies $\si_1(p_0)=0$,
where we use the vector space structure on the fiber of $A^\natural\to \hat{A}$ over zero.
Let us also set $U_2=C\setminus p_0$ and choose any lifting $\si_2:U_2\to A^\natural$. 
Then the above construction gives an extension $\PP(f)^{NC}_{\si_1,\si_2,\varphi}$ of $\PP(f)$
to the relative NC-thickening $C\times A^{NC}$, equipped with a trivialization
over $p_0\times A^{NC}$.

\begin{prop}\label{J-morphism-prop} 
There exist a morphism of NC-schemes $\phi:A^{NC}\to J^{NC,mod}$,
where $J^{NC,mod}$ is the modular NC-thickening of the Jacobian $J=J(C)$
and an isomorphism 
$$\PP(f)^{NC}_{\si_1,\si_2,\varphi}\simeq (\id_C\times\phi)^*\PP^{NC,mod}$$ 
of line bundles on
$C\times A^{NC}$, compatible with trivializations over $p_0\times J^{NC}$ and with the
morphism of abelian varieties $A\to J$ associated with $f$.
\end{prop}

\Pf . This follows immediately from Corollary \ref{Jac-mod-cor}.
\ed

\begin{cor}\label{J-identification-cor} Let $f: C\to J$ be the standard embedding of a curve into its Jacobian,
such that $f(p_0)=0$. Applying the above construction to some liftings $U_i\to J^\natural$, where 
$C=U_1\cup U_2$ is an open covering, we get an extension $\wt{\PP^{NC}}$ of the Poincar\'e line bundle
$\PP$ on $C\times J$ to $C\times J^{NC}$ (trivialized over $p_0$), where $J^{NC}$ is the standard
NC-smooth thickening of $J$ as an abelian variety. There exists an isomorphism
$$\phi:J^{NC}\rTo{\sim} J^{NC,mod}$$ 
of NC-thickenings of $J$ and an isomorphism
$$\wt{\PP}^{NC}\simeq (\id_C\times\phi)^*\PP^{NC,\mod}$$
of line bundles on $C\times J^{NC}$, compatible with trivializations over $p_0\times J^{NC}$ and
inducing the identity isomorphism of $\PP$ over $C\times J$.
\end{cor}

\Pf . This follows from Proposition \ref{J-morphism-prop} together with the fact that any morphism of NC-smooth thickenings
of a given variety is an isomorphism (see \cite[Prop.\ 4.3.1(a)]{Kapranov}).
\ed

\begin{rem}\label{twisted-ab-var-conn} 
Suppose a curve $C$ admits a nonconstant map to an elliptic curve $E$.
Then we have the induced map $f:J(C)\to E$ and hence, a non-standard
NC-thickening $\wt{J^{NC}}$ of $J=J(C)$, given by $\und{\End}(\LL)$, where $\LL$ is a line bundle
on the standard NC-smooth thickening $J^{NC}$, extending $f^*\OO_E(p)$ for some point $p\in E$
(see Corollary \ref{nonstandard-NC-cor}). If the Poincar\'e line bundle $\PP$ on $C\times J$ 
extended to a line bundle on
$C\times \wt{J^{NC}}$, trivialized over $p_0$, then arguing as in Corollary \ref{J-identification-cor},
we would get an isomorphism of $\wt{J^{NC}}$ with the
standard NC-smooth thickening. Since we know that $\wt{J^{NC}}$ is non-standard, it follows that
$\PP$ cannot be extended to a line bundle on $C\times \wt{J^{NC}}$, trivialized over $p_0$.
Let $C=U_1\cup U_2$ be an affine open covering of $C$, equipped with liftings $\si_i:U_i\to J^\natural$,
where $p_0\not\in U_2$, $p_0\in U_1$ and $\si_1(p_0)=0$. 
Then as in the above Construction we get relative flat connections $\nabla_i$ on the restrictions 
$\PP_C|_{U_i\times J}$ for $i=1,2$, and then line bundles $\bK_{U_i,J}(\PP_C|_{U_i\times J},\nabla_i)$ on
$U_i\times \wt{J^{NC}}$, where $\bK_{U_1,J}(\PP_C|_{U_1\times J},\nabla_1)$ is trivialized over $p_0$.
Since these do not glue into a global line bundle on $C\times\wt{J^{NC}}$,
this implies that $\bK_{U_{12},J}(\PP_C|_{U_{12}\times J},\nabla_1)$ and 
$\bK_{U_{12},J}(\PP_C|_{U_{12}\times J},\nabla_2)$
are not isomorphic. Thus, the functor $\bK_{U_{12},J}$ 
does depend on a choice of a connection for our non-standard
thickening $\wt{J^{NC}}$ (for the standard thickening it does not, by Proposition \ref{D-module-prop}(i)).
\end{rem}

We now turn to the study of the NC-Fourier-Mukai transform
$$F^{NC}: D(\Qcoh(C))\to D(\mod-\OO^{NC}_A):E\mapsto Rp_{2*}(\PP(f)^{NC}\ot_{p_1^{-1}\OO_C} p_1^{-1}E)$$
associated with a morphism $f:C\to \hat{A}$, where $\PP(f)^{NC}=\PP(f)^{NC}_{\si_1,\si_2,\varphi}$ is
the NC-extension of the line bundle $\PP(f)$ on $C\times A$,
where for a smooth variety we set $D^b(X)=D^b(\Coh(X))$.
Let also $F:D^b(C)\to D^b(A)$ denote the usual Fourier-Mukai functor associated with $\PP(f)$.
Apriori the image of $F^{NC}$ is
just in the derived category of $\OO_{A^{NC}}$-modules.
Let us define the subcategory of {\it globally perfect} objects 
$$\Per^{gl}(A^{NC})\sub  D(\mod-\OO^{NC}_A)$$
as the full subcategory consisting of objects represented globally by a finite complex of
vector bundles on $A^{NC}$.

\begin{thm}\label{curve-Fourier-thm} 
For every $E\in D^b(C)$ one has $F^{NC}(E)\in\Per^{gl}(A^{NC})$.
Furthermore, if $E$ is a coherent sheaf on $C$ such that $F(E)$ is a vector bundle on $A$,
where $F:D^b(C)\to D^b(A)$ is the usual Fourier-Mukai transform associated with $\PP(f)$, then 
$F^{NC}(E)$ is a vector bundle on $A^{NC}$.
We have a commutative diagram of functors
\begin{equation}\label{F-res-diagram}
\begin{diagram}
D^b(C)&\rTo{F^{NC}}& \Per^{gl}(A^{NC})\\
&\rdTo{F}&\dTo_{}\\
&&D^b(A)
\end{diagram}
\end{equation}
where the vertical arrow is the natural restriction functor.
\end{thm}

\Pf . First, to prove that $F^{NC}(E)$ is globally perfect it is enough to check this 
for a set of objects $E$ generating $D^b(C)$
as a triangulated category. Since $F^{NC}(\OO_p)$, where $p\in C$, is a line bundle on $A^{NC}$,
it is enough to consider the case when $E$ is a coherent sheaf on $C$ such that $F(E)$
is a vector bundle on $A$.
Let $j_i:U_i\hra C$ and $j_{12}:U_{12}\hra C$ be the natural open embeddings.
Recall that the restriction $\PP(f)^{NC}|_{U_i\times A}$ is the $\OO^{NC}_{U_i\times A/U_i}$-module
associated with the relative flat connection on $\PP(f)|_{U_i\times A}$ coming from
the isomorphism 
$$\PP(f)|_{U_i\times A}\simeq (\si_i\times\id_A)^*\PP^\natural$$ 
and the natural relative flat connection $\nabla^\natural$ on $\PP^\natural$. Hence,
by Lemma \ref{rel-conn-cd-lem}(ii), we get an isomorphism
$$\PP(f)^{NC}|_{U_i\times A}\simeq (\si_i\times\id_A)^*\bK_{A^\natural,A}(\PP^\natural,\nabla^\natural).$$ 
Thus, for a quasicoherent sheaf $\FF$ on $U_i$ we have an isomoprhism
$$F^{NC}(j_{i*}(\FF))\simeq Rp_{2*}(\PP(f)^{NC}|_{U_i\times A}\ot_{p_1^{-1}\OO_{U_i}} p^{-1}\FF)\simeq
 F^{NC,\natural}(\si_{i*}\FF),$$
 where $F^{NC,\natural}$ is the functor \eqref{FNC-natural-eq}.
Let us write for brevity $\bK=\bK_{pt,A}$.
Combining this isomorphism with Proposition \ref{F-NC-natural-prop} we get
\begin{equation}\label{FNC-K-isom1}
F^{NC}(j_{i*}E|_{U_i})\simeq F^{NC,\natural}(\si_{i_*}E|_{U_i})\simeq
\bK(F^D(\si_{i*}E|_{U_i}))
\end{equation}
for $i=1,2$, where $F^D(\si_{i*}E|_{U_i})\simeq F(j_{i*}E|_{U_i})$ as an $\OO_A$-module, but it is equipped with
a connection coming from $\si_i$.
Similarly, we have two isomorphisms
\begin{equation}\label{FNC-K-isom2}
F^{NC}(j_{12*}E|_{U_{12}})\simeq \bK(F^D(\ov{\si}_{i*}E|_{U_{12}})),
\end{equation}
where for $i=1,2$, $\ov{\si}_i=\si_i|_{U_{12}}:U_{12}\to A^\natural$.
The isomorphisms \eqref{FNC-K-isom1}, \eqref{FNC-K-isom2} fit into the commutative diagrams (for $i=1,2$) 
\begin{diagram}
F^{NC}(j_{i*}(E|_{U_i}))&\rTo{\sim}& \bK(F^D(\si_{i*}E|_{U_i}))\\
\dTo{}&&\dTo{}\\
F^{NC}(j_{12*}E|_{U_{12}})&\rTo{\sim}&\bK(F^D(\ov{\si}_{i*}E|_{U_{12}}))
\end{diagram}
where the vertical arrows are induced by restrictions from $U_i$ to $U_{12}$.
Note that we have isomorphisms of $\OO_A$-modules
$$F^D(\ov{\si}_{1*}E|_{U_{12}}))\simeq F^D(\ov{\si}_{2*}E|_{U_{12}}))\simeq F(j_{12*}E|_{U_{12}}),$$
however, the two $D$-module structures are different: one connection, $\nabla_1$, is induced by $\si_1$
and another, $\nabla_2$, by $\si_2$.
Applying $F^{NC}$ to the Cech exact sequence
$$0\to E\to j_{1*}E|_{U_1}\oplus j_{2*}E|_{U_2}\to j_{12*}E|_{U_{12}}\to 0
$$
we get an exact triangle
$$F^{NC}(E)\to \bK(\VV_0)\rTo{\a} \bK(\VV_1)\to F^{NC}(E)[1],$$
where $\VV_0=F(j_{1*}E|_{U_1})\oplus F(j_{2*}E|_{U_2}))$ and $\VV_1=F(j_{12*}E|_{U_{12}})$
are viewed as $D$-modules on $A$ (for $\VV_1$ we take the connection induced by $\si_1$), 
and $\a$ is the morphism coming from the above identifications \eqref{FNC-K-isom1}, \eqref{FNC-K-isom2}.
We can represent $\a$ as the composition 
\begin{align*}
&\bK(F(j_{1*}E|_{U_1}))\oplus \bK(F(j_{2*}E|_{U_2})))\to\\
&\bK(F(j_{12*}E|_{U_{12}}),\nabla_1)\oplus \bK(F(j_{12*}E|_{U_{12}}),\nabla_2)\rTo{(-\id,\varphi)}
\bK(F(j_{12*}E|_{U_{12}}),\nabla_1)
\end{align*}
where $\varphi$ is induced by some isomorphism of $\AA$-dg-modules
$$K(F(j_{12*}E|_{U_{12}}),\nabla_2)\rTo{\sim} K(F(j_{12*}E|_{U_{12}}),\nabla_1).$$
Thus, $\a$ is the map on $\und{H}^0$ induced by a closed morphism of degree zero of
$\AA$-dg-modules 
$$f:\VV_0\hat{\ot}_\OO\AA=K(F(j_{1*}E|_{U_1})\oplus F(j_{2*}E|_{U_2}))\to
K(F(j_{12*}E|_{U_{12}}),\nabla_1)=\VV_1\hat{\ot}_\OO\AA.$$
Hence, $F^{NC}(E)$ is realized by the complex 
$$\und{H}^0(\VV_0\hat{\ot}_\OO \AA,d)\rTo{\und{H}^0(f)} \und{H}^0(\VV_1\hat{\ot} \AA,d)$$
Note that the complex
$$\VV_0 \rTo{f|_{\VV_0}\mod F^1} \VV_1$$
represents $F(E)$.
In other words, the map $\VV_0\to \VV_1$ induced by $f\mod F^1+\AA^1$ is surjective,
and its kernel is $F(E)$. This immediately implies that $f$ is surjective. Next, we claim that
$\ker(f)$ is locally isomorphic to $F(E)\hat{\ot}_\OO \AA$, as an $\AA$-module.
Indeed, locally we can choose a splitting $s:\VV_1\to\VV_0$, so that the corresponding projection 
$\pi:\VV_0\to F(E)$ fits into an exact sequence
$$0\to \VV_1\rTo{s} \VV_0\rTo{\pi} F(E)\to 0.$$
Hence, we have an induced exact sequence of $\AA$-modules
$$0\to \VV_1\hat{\ot}_\OO \AA\rTo{s\ot \id} \VV_0\hat{\ot}_\OO \AA\rTo{\pi\ot\id} F(E)\hat{\ot}_\OO\AA\to 0.$$
The composed morphism of $\AA$-modules
$$\VV_1\hat{\ot}_\OO\AA\rTo{s\ot\id} \VV_0\hat{\ot}_\OO\AA\rTo{f} \VV_1\hat{\ot}_\OO\AA$$
induces the identity modulo $F^1+\AA^1$, hence, it is an isomorphism.
By diagram chasing this implies that the composed morphism
$$\ker(f)\to \VV_0\hat{\ot}_\OO\AA\rTo{\pi\ot\id}F(E)\hat{\ot}_\OO\AA$$
is an isomorphism. This proves our claim.
Note that the differential $d$ on $\VV_0\hat{\ot}_\OO\AA$ induces a differential $d$ on
$\ker(f)$, so that we have an exact sequence of dg-modules over $\AA$
\begin{equation}\label{dg-mod-ex-seq}
0\to (\ker(f),d)\to (\VV_0\hat{\ot}_\OO\AA,d)\to \VV_1\hat{\ot}_\OO\AA\to 0.
\end{equation}
By Theorem \ref{module-thm}(ii), 
for any dg-module structure on $F(E)\hat{\ot}_\OO\AA$
the higher cohomology of the differential vanish, while $\und{H}^0$ is a vector bundle on $A^{NC}$.
Hence, the same is true for the dg-module $(\ker(f),d)$,
and the exact sequence \eqref{dg-mod-ex-seq} leads to an exact sequence of 
$\OO_A^{NC}$-modules
$$0\to \und{H}^0(\ker(f),d)\to 
\und{H}^0(\VV_0\hat{\ot}_\OO \AA,d)\rTo{\und{H}^0(f)} \und{H}^0(\VV_1\hat{\ot} \AA,d)\to 0.$$
Therefore, $F^{NC}(E)\simeq \und{H}^0(\ker(f),d)$ which is a vector bundle on $A^{NC}$.

Note that we have a natural isomorphism 
$$\PP(f)^{NC}\ot_{\OO^{NC}_{C\times A/C}}\OO_{C\times A}\rTo{\sim} \PP(f)$$
of $\OO_{C\times A}$-modules.
Hence, we obtain a natural map 
\begin{equation}\label{F-res-map}
\begin{array}{l}
F^{NC}(E)\ot^{\Bbb L}_{\OO^{NC}_A}\OO\to 
Rp_{2*}(p_1^{-1}E\ot_{p_1^{-1}\OO_C} \PP(f)^{NC}\ot_{\OO^{NC}}\OO)
\rTo{\sim} Rp_{2*}(p_1^{-1}E\ot_{p_1^{-1}\OO_C}\PP(f))\\
=F(E).
\end{array}
\end{equation}
If $E$ is a coherent sheaf such that $F(E)$ is a vector bundle, i.e., $E$ is an {\it $\IT_0$-sheaf} in
the terminology of \cite{Mukai}, then the above argument shows that
the map \eqref{F-res-map} is an isomorphism. Since every coherent sheaf can be embedded into an
$\IT_0$-sheaf, we can calculate $F^{NC}$ using resolutions of $\IT_0$-sheaves, and
the commutativity of the diagram \eqref{F-res-diagram} follows.
\ed

The above Theorem immediately implies that the Picard bundles on the Jacobian extend to
the standard NC-thickening.

\begin{cor} Let $V$ be the vector bundle on the Jacobian $J=J(C)$ obtained as the Fourier
transform of a coherent sheaf on $C$. Then $V$ extends to a vector bundle on 
the standard NC-thickening $J^{NC}$.
\end{cor}




\section{Analytic NC manifolds}\label{NC-analytic-sec}
\newcommand{\ord}{\mathsf{ord}}
\newcommand{\Jet}{\mathsf{Jet}}
\newcommand{\obs}{\mathsf{obs}}
\newcommand{\Ad}{\sf Ad}

\subsection{Analytic NC manifolds.}~\label{analytic-nc-manifolds}

Let $X$ be a complex manifold. We define the notion of a \emph{holomorphic NC-connection}
in the same way as an algebraic NC-connection (see Definition \ref{alge-nc-defi}), replacing
algebraic forms by holomorphic ones. The constructions of Theorem~\ref{dg-constr-thm}
go through in the holomorphic setting. Thus, starting with a holomorphic
torsion free connection $\nabla$, we obtain a holomorphic NC-connection $D$, which is a differential on
\[\AA^{an}_X:=(\Omega_X^\bullet)^{an}\otimes_{\OO_X^{an}} \hat{T}_{\OO_X^{an}}(\Om^{1})^{an},\]
the analytic version of the sheaf $\AA_X$ studied in Theorem~\ref{dg-constr-thm}.

Let $U$ be an open subset of $\C^n$.
Complex analytic coordinates $x_1,\ldots,x_n$ induce a trivialization 
\[ (\Om^1_U)^{an} \cong \OO_U^{an}dx_1\oplus\ldots\oplus \OO_U^{an}dx_n.\] 
of the cotangent bundle.
Let us set $e_i:= dx_i$. The above trivialization induces a \emph{flat} torsion free connection 
$\nabla$ on $(\Om^1_U)^{an}$ such that $\nabla(e_i)=0$. 
We have the corresponding NC-connection $D=d_{\nabla,\id}$ (see Lemma~\ref{torsion-free-lem}).

\medskip
In view of Theorem~\ref{smoothness-thm} we define local models of an analytic NC-smooth thickening of $U$ by 
\[ \OO_U^{nc}:= \underline{H}^0(D)=\ker(D:(\AA^0_U)^{an}\to (\AA^1_U)^{an}).\]
This local model $\OO_U^{nc}$ can also be described explicitly in the manner of Theorem~\ref{local-constr-thm}. Indeed, a local section of $\OO_U^{nc}$ can be uniquely written as an infinite series $\sum_\lambda [[f_\la(e)]] M_\la$ where $M_\la$ are monomials in the Lie words of degree $\ge 2$ in $e_i$, 
$f_\la$ are holomorphic functions on $U$, and $[[f_\la(e)]]$ is an infinite series of the form
\[ [[f_\la(e)]] = \sum (-1)^{i_1+\ldots i_n} (\partial_1)^{i_1}\ldots(\partial_n)^{i_n}f_\la e_1^{i_1}\cdots e_n^{i_n}\]
where $\partial_i=\partial/\partial x_i$. As before, we use the grading on $U(\LLie_+(e_1,\ldots,e_n))$
induced by the grading on the free Lie algebra, where $\deg(e_i)=1$.
We define $\deg(M_\la)$ by viewing a monomial $M_\la$ as an
element of $U(\LLie_+(e_1,\ldots,e_n))$ and using the above grading.
The $I$-filtration on $\OO_U^{nc}$  (see \eqref{I-n-def-eq}) can be described explicitly as follows:
$I_d\OO_U^{nc}$ consists of formal series $\sum_\lambda [[f_\la(e)]] M_\la$ such that $\deg(M_\la)\ge d$. Locally the corresponding quotient admits a trivialization
\[ \OO_U^{nc}/I_{d+1}\OO_U^{nc} \cong \OO_U^{an}\otimes_\C\big(U(\LLie_+(e_1,\ldots,e_n))_0\oplus \ldots \oplus U(\LLie_+(e_1,\ldots,e_n))_d\big)\]
as sheaves of abelian groups over $U$. The product structure, written in this trivialization, is described in~\cite[Proposition (3.4.3)]{Kapranov}. Again we let $\Aut(\OO^{nc}/I_{d+1}\OO^{nc})$ be automorphisms that are compatible with the projection map $\pi: \OO^{nc}/I_{d+1}\OO^{nc} \rightarrow \OO^{an}$.

\begin{defi}\label{analytic-def}
An \emph{analytic NC manifold} is a ringed space $(X,\OO^{nc})$ which is locally of the form $(U,\OO_U^{nc})$.
\end{defi}

\begin{lem}\label{surj-lem}
The natural homomorphism $\underline{\Aut}(\OO^{nc}/I_{d+1}\OO^{nc}) \rightarrow \underline{\Aut}(\OO^{nc}/I_{d}\OO^{nc})$ is surjective. 
\end{lem}

\Pf. The question is local, so we can work with the explicit local model described above. In this case, for any automorphism $\phi\in \Aut(\OO^{nc}/I_{d}\OO^{nc})$, let us write
\[ \phi([[x_i]])=[[x_i]]+\sum_{d-1\geq \deg(M_\la)\geq 2} [[f^i_\la(e)]]M_\la.\]
We define an $\OO^{an}$-linear automorphism $\tilde{\phi}$ of $\hat{T}_\OO^{an}((\Om^1)^{an})=\hat{T}_\OO^{an}(\OO^{an}e_1\oplus\ldots\oplus\OO^{an}e_n)$ which acts on the generators by
\[ \tilde{\phi}: e_i \mapsto e_i-\sum_{d-1\geq \deg(M_\la)\geq 2} [[f^i_\la(e)]]M_\la.\]
Since the element $\sum_{d-1\geq \deg(M_\la)\geq 2} [[f^i_\la(e)]]M_\la$ is $D$-closed, this automorphism commutes with the differential $D$. It is also straightforward to verify that $\tilde{\phi}$ induces an automorphism on $\ker(d)$ that extends the automorphism $\phi$. This shows that the composition
\[ \Aut(\OO^{nc}) \rightarrow \Aut(\OO^{nc}/I_{d+1}\OO^{nc}) \rightarrow \Aut(\OO^{nc}/I_d\OO^{nc})\]
is surjective. Hence, the second arrow is also surjective.\ed

For the proof of Theorem \ref{main-an-thm} 
we will need the following analytic analog of the notion of a $(d)_I$-smooth thickening.

\begin{defi}\label{analytic-d-def}
An \emph{analytic $(d)_I$-manifold} is a ringed space $(X,\OO_X^{(d)})$ which is locally of the form 
$(U,\OO^{nc}/I_{d+1}\OO^{nc})$. 
\end{defi}

By an \emph{algebraic $(d)_I$-manifold} we mean a $(d)_I$-smooth thickening of a smooth algebraic variety.
Let $(X,\OO_X^{(d)})$ be an algebraic (resp., analytic) $(d)_I$-manifold. 
For an open subset $U\subset X$ we define a category $Th^{(d+1)}(\OO_X^{(d)},U)$ whose objects are 
$(d+1)_I$-manifold thickenings of $\OO_X^{(d)}|_U$, and morphisms are morphisms of ringed spaces that are identical on $\OO_X^{(d)}|_U$.

\begin{lem}~\label{classification-lem}
Let $(X,\OO_X^{(d)})$ be an algebraic (resp., analytic) $(d)_I$-manifold. 
The association $U\mapsto Th^{(d+1)}(\OO_X^{(d)},U)$ defines a gerbe in the Zariski (resp., classical) 
topology of $X$ banded by  
$$\und{\Hom}(\Omega_X^1, (U\LLie_+(\Omega_X^1))_{d+1}).$$
\end{lem}

\Pf. The fact that all morphisms are isomorphisms is checked by induction in $Q$ using standard
facts about central extensions (see \cite[Prop.\ 1.2.6]{Kapranov}). The local triviality in the algebraic
case was checked in Proposition~\ref{trunc-twisted-prop} (in the analytic case the local triviality is a part
of the definition).
It remains to observe that by \cite[Prop.\ 1.2.5]{Kapranov}, an automorphism
of $\OO_U^{(d+1)}\in \CC(\OO_X^{(d)},U)$ trivial on $\OO_U^{(d)}$ corresponds to a derivation of
$\OO_U$ with values in 
$$I_{d+1}\OO_U^{(d+1)}\simeq (U\LLie_+(\Omega_U^1))_{d+1}.$$
\ed

\subsection{Analytification functors.} Let $(X,\OO^{NC})$ be an algebraic NC manifold, i.e., an NC-smooth thickening of a smooth algebraic variety over $\C$. We can associate with $(X,\OO^{NC})$
an analytic NC manifold $(X^{an},\OO^{nc})$ as follows. 
The underlying analytic manifold $X^{an}$ is simply the analytification of $X$. To construct $\OO^{nc}$
we observe that by Theorem~\ref{twisted-dg-thm}, we have a global dg resolution $(\AA_\TT,D)$ of $\OO^{NC}$ that is locally (in the Zariski topology) of the form $(\AA_U,D)$ considered above. Now we simply apply the analytification functor to $\AA_\TT$ and set
$$\OO^{nc}:=\und{H}^0(\AA_\TT^{an},D).$$
We refer to the resulting analytic NC manifold $(X^{an},\OO^{nc})$ as the \emph{analytification of} $(X,\OO^{NC})$. 

The analytification of an algebraic $(d)_I$-manifold is defined similarly using Proposition~\ref{trunc-twisted-prop}.

\medskip

\noindent
{\it Proof of Theorem \ref{main-an-thm}}.
It suffices to prove that the analytification functor 
\begin{equation}\label{d-I-gerbes-map}
Th^{(d)_I}_X\to Th_{X^{an}}^{(d)_I}
\end{equation}
between the gerbes of algebraic and analytic $(d)_I$-manifold thickenings of $X$ and $X^{an}$,
is an equivalence. 
Let us prove by induction in $Q$ that for an algebraic $(d)_I$-manifold $\OO_X^{(d)}$ the map
\[\Aut(X, \OO_X^{(d)}) \to \Aut(X^{an}, (\OO_X^{(d)})^{an}).\]
is an isomorphism (where we consider automorphisms trivial on the abelianization). 
We use the short exact sequence
\[  0 \ra \underline{\Hom}\big(\Omega_X^1, U\LLie_+(\Omega_X^1)_{d}\big) \ra \underline{\Aut}( \OO^{NC}/I_{d+1}\OO^{NC})\ra \underline{\Aut}(\OO^{NC}/I_{d}\OO^{NC})\ra 0,\]
which gives the corresponding long exact sequence
\[ \mbox{ \footnotesize $0\ra \Hom\big(\Omega_X^1, U\LLie_+(\Omega_X^1)_{d}\big) \ra 
\Aut( \OO^{NC}/I_{d+1}\OO^{NC})\ra \Aut(\OO^{NC}/I_{d}\OO^{NC}) \ra \Ext^1\big(\Omega_X^1, U\LLie_+(\Omega_X^1)_{d+1}\big)\ra\ldots$}.\]
Lemma \ref{surj-lem} implies that
we also have a similar long exact sequence in the analytic setting. Now our claim follows by the GAGA theorem and the five lemma.
The fact that \eqref{d-I-gerbes-map} is essentially surjective is also deduced by induction in $Q$,
using Lemma \ref{classification-lem} and the GAGA theorem for
the cohomology of the coherent sheaf $\und{\Hom}(\Omega_X^1, (U\LLie_+(\Omega_X^1))_{d+1})$ and
its analytification.
\ed

In the remaining subsections, we will construct examples of NC manifolds using analytic methods.

\subsection{$L_\infty$ spaces in K\"ahler geometry.} To explain how K\"ahler geometry might enter into constructions of analytic NC manifolds, we first recall an analogous construction of Costello (see~\cite{Costello}) in the commutative case: the construction of an \emph{$L_\infty$ space} associated to a K\"ahler manifold $(X,g)$. In fact Costello's construction works for any complex manifold. However, in the K\"ahler situation we can write down explicit formulas of the $L_\infty$ structure in terms the curvature tensor of a K\"ahler metric $g$, following a beautiful paper of Kapranov~\cite{Kapranov2}.

Consider the infinite dimensional vector bundle $\Jet_X$, the infinite holomorphic jet bundle of $X$. The bundle $\Jet_X$ admits a decreasing filtration whose associated graded are symmetric powers of $\Omega^1_X$. This filtration almost never splits, thus in general we do not have an isomorphism
\[\Jet_X\cong \hat{S}_{\OO_X} \Omega^1_X.\]
However, Kapranov~\cite{Kapranov} observed that we can write the Dolbeault complex of $\Jet_X$ in terms of the Dolbeault complex of $\hat{S}_{\OO_X} \Omega^1_X$ endowed with a perturbed differential. This perturbed differential can be explicitly written down using the Levi-Civita connection $\nabla$ of $g$ and its associated curvature tensor $R$, which we briefly recall in the following paragraph.

To avoid confusions, we will denote by $\Omega^{1,0}_X$ the sheaf of smooth sections of the holomorphic vector bundle $\Omega^1_X$. The curvature operator maybe written as a morphism
\[ R: T_X^{1,0}\otimes T_X^{1,0} \ra \Omega^{0,1}_X\otimes T_X^{1,0}.\]
Due to the K\"ahler condition, the operator $R$ factor through
\[ R: \Sym^2 T_X^{1,0} \ra \Omega^{0,1}_X\otimes T_X^{1,0}.\]
For each $n\geq 2$ we define an operator $R_n\in \Omega^{0,1}_X\big(\Hom(\Sym^2T_X^{1,0}\ot (T_X^{1,0})^{\ot (n-2)}, T_X^{1,0})\big)$ inductively by putting
\[ R_2:= R, \mbox{\;\; and \;\;} R_{n+1}:=\nabla R_{n}.\]
Kapranov showed that $R_n$ is a totally symmetric operator, i.e. it is an operator
\[ \Sym^n T_X^{1,0} \ra \Omega^{0,1}_X \ot T_X^{1,0}.\]
Dualizing, we obtain a morphism
\[ R_n^\vee: \Omega^{1,0}_X \ra \Omega^{0,1}_X\otimes \Sym^n \Omega^{1,0}_X.\]
Now the perturbed differential $D$ on the Dolbeault complex $\Omega^{0,*}_X(\hat{S} \Omega^1_X)$ is defined to be the unique derivation of this algebra extending the $\overline{\partial}$ operator, and acting on $\Omega^{1,0}$ by the infinite sum $\overline{\partial}+\sum_{n\geq 2} R_n^\vee$.
The following Theorem is obtained by Kapranov~\cite[Theorem 2.8.2]{Kapranov2}.
\begin{thm}
We have $(\overline{\partial}+\sum_{n\geq 2} R_n^\vee)^2=0$. Furthermore, there is an isomorphism of commutative differential graded algebras
\[ \Omega^{0,*}(\Jet_X) \cong \big(\Omega^{0,*}(\hat{S}(\Omega^1_X)), \overline{\partial}+\sum_{n\geq 2} R_n^\vee\big).\]
\end{thm}

It is well-known that $\Jet_X$ has a flat holomorphic connection. Furthermore the canonical embedding
\[ i: \OO_X\hookrightarrow \Omega^*(\Jet_X)=\Omega^*\ot_{\OO_X} \Jet_X,\]
sending a holomorphic function to the associated constant jets, is a quasi-isomorphism of differential graded commutative algebras. Applying Kapranov's Theorem to replace $\Jet_X$ by $\big(\Omega^{0,*}(\hat{S}(\Omega^1_X)), \overline{\partial}+\sum_{n\geq 2} R_n^\vee\big)$, we obtain an explicit resolution of $\OO_X$ by a quasi-free commutative differential graded algebra over the smooth de Rham algebra $\Lambda_X:=\oplus_{p,q} \Omega_X^{p,q}$. 

More precisely the underlying algebra is simply $\Lambda_X\otimes_{C^\infty_X} \hat{S}_{C^\infty}(\Omega_X^{1,0})$. The differential on this algebra is the unique derivation which extends the de Rham differential on $\Lambda_X$, and which acts on the space $\Omega_X^{1,0}$ by the infinite sum $d_\tau+ \overline{\partial}+\nabla+ \sum_{n\geq 2} R_n^\vee$ (as depicted below).
\begin{diagram}
\Lambda^1 &&  \Lambda^1\otimes \Omega_X^{1,0}  &&\Lambda^1\otimes \Sym^2\Omega_X^{1,0} && \cdots \\
&\luTo(2,2)~{d_\tau} & \uTo~{\overline{\partial}+\nabla} & \ruTo(2,2)~{R^\vee_2}  &\cdots \\
 & & \Omega_X^{1,0} & 
\end{diagram}
Here $d_\tau(1\otimes \alpha):=\alpha\otimes 1$ using the fact that $\Omega_X^{1,0}$ is a summand of $\Lambda^1$, and $\nabla$ is the $(1,0)$-type connection associated to the K\"ahler metric $g$. There is a quasi-isomorphism of commutative differential graded algebras
\begin{equation}\label{resolution-eq}
\OO_X \cong \big(\Lambda_X\otimes_{C^\infty_X} \hat{S}_{C^\infty}(\Omega_X^{1,0}),d_\tau+ \overline{\partial}+\nabla+ \sum_{n\geq 2} R_n^\vee\big).
\end{equation}
This resolution, referred to as the \emph{$L_\infty$ space} associated to $X$, plays an important role in the recent work of Costello on Witten genus~\cite{Costello}.

\subsection{$A_\infty$ spaces and NC thickenings.} Now we are going to 
establish an NC analog of \eqref{resolution-eq} (see Proposition \ref{if-part-prop} below).
Let $X$ be a complex manifold. The NC analog of the differential appearing in \eqref{resolution-eq}
is given by the notion of an \emph{NC-connection} $\DD$, which is a differential on the algebra
\[\Lambda_X\otimes_{C^\infty_X} \hat{T}_{C_X^\infty}(\Omega_X^{1,0})\] 
extending the de Rham differential on $\Lambda_X$, and acting on the space of generators $\Omega_X^{1,0}$ by an operator of the form 
\[d_\tau+\overline{\partial}+\nabla+A^\vee_2+A^\vee_3+\ldots,\]
where $\nabla$ is a $(1,0)$-type connection on the bundle $\Omega_X^{1,0}$
(see Definition \ref{nc-connection-def}). For $k\ge 2$ we split
the $C^\infty_X$-linear morphism $A^\vee_k: \Omega_X^{1,0} \ra \Lambda^1_X\ot (\Omega_X^{1,0})^{\ot k}$ 
as a sum
\[ A^\vee_k:= E_k+ F_k\]
where $E_k: \Omega_X^{1,0} \ra \Lambda^{1,0}_X\ot (\Omega_X^{1,0})^{\ot k}$ and $F_k:\Omega_X^{1,0} \ra \Lambda^{0,1}_X\ot (\Omega_X^{1,0})^{\ot k}$. Thus, the operator $\DD$ is of the form
\[ \DD^{1,0}+\DD^{0,1}:=(d_\tau+\nabla+E_2 +E_3+\ldots)+(\overline{\partial} + F_2 +F_3+\ldots).\]
The integrability condition $\DD^2=0$ then implies that 
\[(\DD^{1,0})^2=[\DD^{1,0},\DD^{0,1}]=(\DD^{0,1})^2=0.\]
Thus, the complex $\big(\Lambda_X\otimes_{C^\infty_X} \hat{T}_{C^\infty}(\Omega_X^{1,0}),\DD\big)$ is in fact a bicomplex whose degree $(p,q)$-component is  \[\Lambda^{p,q}\ot_{C^\infty_X} \hat{T}_{C^\infty}(\Omega_X^{1,0}).\]
The following proposition describes the local structure of the cohomology of $\DD^{0,1}$.
\medskip
\begin{prop}~\label{nc-jet-prop}
Let $U$ be a Stein open subset of $X$. Then there exists an algebra automorphism $\phi$ of $\hat{T}_{C^\infty}(\Omega_U^{1,0})$ of the form
\[ \phi(\alpha)=\alpha+\phi_2(\alpha)+\phi_3(\alpha)+\ldots\]
for $\alpha\in \Omega_U^{1,0}$, and $\phi_n: \Omega_U^{1,0}\ra (\Omega_U^{1,0})^{\ot n}$ are $C^\infty_U$-linear morphisms, such that on the Dolbeault complex $\Lambda_U^{0,*}\ot_{C^\infty_U}\hat{T}_{C_U^\infty}(\Omega_U^{1,0})$ we have the identity
\[(\id\ot \phi)^{-1}(\overline{\partial})(\id\ot\phi)= \DD^{0,1}.\]
This implies that $\underline{H}^i\big(\Lambda_U^{0,*}\ot_{C^\infty_U}\hat{T}_{C^\infty}(\Omega_U^{1,0}),\DD^{0,1}\big) =0$ if $i\geq 1$, and $\underline{H}^0\big(\Lambda_U^{0,*}\ot_{C^\infty_U}\hat{T}_{C^\infty}(\Omega_U^{1,0}),\DD^{0,1}\big)$ is isomorphic to $\hat{T}_{\OO_U}(\Omega_U)$ where $\Omega_U$ denotes the space of holomorphic one forms.
\end{prop}

\Pf. We inductively construct $\phi_n$'s as follows. For the construction of $\phi_2$, we use the fact that $\DD^{0,1}=\overline{\partial} + F_2 +F_3+\ldots$ squares to zero. This implies in particular that $\overline{\partial} F_2=0$. Over a Stein manifold the $\overline{\partial}$-closedness implies the exactness. We set $\phi_2$ so that $\overline{\partial} \phi_2=F_2$. Assuming there exists $\phi_k$ for $k\leq n$ such that 
\[ (\id\ot \phi_{\leq n})\circ \DD^{0,1}\circ (\id\ot \phi_{\leq n})^{-1}=\overline{\partial} \pmod{\hat{T}^{\geq n+1}}.\]
This implies that
\[ (\id\ot \phi_{\leq n})\circ \DD^{0,1}\circ (\id\ot \phi_{\leq n})^{-1}=\overline{\partial}+G_{n+1}\pmod{\hat{T}^{\geq n+2}}\]
for some morphism $G_{n+1}: \Omega_U^{1,0}\ra \Lambda_X^{0,1}\ot (\Omega_U^{1,0})^{\ot n+1}$. The integrability $(\DD^{0,1})^2=0$ implies that $\overline{\partial} G_{n+1}=0$. Since we are over a Stein manifold there exists $\phi_{n+1}: \Omega_U^{1,0}\ra (\Omega_U^{1,0})^{\ot n+1}$ such that $\overline{\partial}(\phi_{n+1})= G_{n+1}$. We then check that
\[ (\id\ot \phi_{\leq n+1})\circ \DD^{0,1}\circ (\id\ot \phi_{\leq n+1})^{-1}=\overline{\partial} \pmod{\hat{T}^{\geq n+2}}.\]
The proposition is proved.\ed

\begin{rem}
In view of the commutative situation, the cohomology sheaf $\underline{H}^0\big(\Lambda_X^{0,*}\ot_{C^\infty_X}\hat{T}_{C^\infty}(\Omega_X^{1,0}),\DD^{0,1}\big)$ should be thought of as a generalization of the notion of jet bundle in the NC setting.
\end{rem}

\medskip

Observe that the natural projection $\underline{H}^0\big(\Lambda_X\otimes_{C^\infty_X} \hat{T}_{C^\infty}(\Omega_X^{1,0}),\DD\big) \ra C^\infty_X$ to smooth functions lands in fact inside holomorphic functions. So there is a morphism of algebras
\[\underline{H}^0\big(\Lambda_X\otimes_{C^\infty_X} \hat{T}_{C^\infty}(\Omega_X^{1,0}),\DD\big) \ra \OO_X.\] 

\medskip

\begin{prop}~\label{if-part-prop}
Let $\DD$ be an NC-connection on a complex manifold $X$. Then we have
$\underline{H}^i\big(\Lambda_X\otimes_{C^\infty_X} \hat{T}_{C^\infty}(\Omega_X^{1,0}),\DD\big)=0$ for $i\geq 1$. Furthermore, the sheaf of algebras $\underline{H}^0\big(\Lambda_X\otimes_{C^\infty_X} \hat{T}_{C^\infty}(\Omega_X^{1,0}),\DD\big)$ defines an analytic NC manifold structure on $X$ with respect to the projection $\underline{H}^0\big(\Lambda_X\otimes_{C^\infty_X} \hat{T}_{C^\infty}(\Omega_X^{1,0}),\DD\big) \ra \OO_X$.
\end{prop}

\Pf. Both statements are local, so it is enough to prove them over a Stein open subset $U$ of $X$. 
As we observed before, the complex $\big(\Lambda_X\otimes_{C^\infty_X} \hat{T}_{C^\infty}(\Omega_X^{1,0}),\DD\big)$ is a bicomplex graded by the $(p,q)$-degree of forms in $\Lambda_X$. In particular, this is a bounded complex in the first quadrant whose cohomology can be computed via the spectral sequence with the $E^1$-page given by
the cohomology of $\DD^{0,1}=\overline{\partial}+F_2+\ldots$. 
By Proposition~\ref{nc-jet-prop}, there exists an automorphism $\phi$ of $\hat{T}_{C^\infty}(\Om^{1,0}_U)$
with $\phi(\a)=\a+\phi_2(\a)+\phi_3(\a)+\ldots$ such that
$(\id\ot\phi)\DD^{0,1}(\id\ot\phi)^{-1}=\dbar$. Hence, $\id\ot\phi$ induces an isomorphism of
the cohomology of $\DD^{0,1}$ with $\Omega_U^*\ot_{\OO_U}\hat{T}(\Omega^1_U)$. 
In particular, the spectral sequence collapses, and 
the total cohomology is isomorphic to the cohomology of
$D=(\id\ot\phi)\DD^{1,0}(\id\ot\phi)^{-1}$ acting on $\Omega_U^*\ot_{\OO_U}\hat{T}(\Omega^1_U)$. 
Note that $D$ is a derivation. We claim that in fact $D$ is a holomorphic NC-connection.
Indeed, for $f\in\OO_U$ we have
$$D(f)=(\id\ot\phi)\DD^{1,0}(f)=(\id\ot\phi)(df\ot 1)=df\ot 1.$$
Also, for $\alpha\in \Omega^1_U$ we have
$$D(1\ot\a)=(\id\ot\phi)(\a\ot 1+\ldots)=\a\ot 1+\ldots,$$
where the dots denote the terms in $\Omega_U^*\ot_{\OO_U}\hat{T}^{\ge 1}(\Omega^1_U)$.
Hence, $D$ is a holomorphic NC-connection, and
the assertion now follows from the holomorphic version of Theorem \ref{main-nc-thm}.
\ed

\begin{prop}
Let $X$ be a complex manifold. Then $X$ admits an analytic NC thickening if and only if there exists an NC-connection on $X$.
\end{prop}

\Pf. The ``if" part follows from Proposition~\ref{if-part-prop}. Conversely, assume that we are given an analytic NC thickening. Then by Theorem~\ref{twisted-dg-thm}~\footnote{More precisely, we use the analytic version of Theorem~\ref{twisted-dg-thm} which is proved similarly.} there exists a quadruple $(\TT,\JJ,\varphi, D)$ in the sense of Definition~\ref{twisting-def}. Since in the smooth category exact sequences of vector bundles always split, there exists an isomorphism between $\TT$ and its associated graded:
\[ \psi:C^\infty(\TT) \cong \prod_{n\geq 0} C^\infty(\JJ^n/\JJ^{n+1}).\]
Followed by the isomorphism $\varphi$ we get an isomorphism
\[ \varphi\circ\psi: C^\infty(\TT) \cong \hat{T} (\Omega^{1,0}_X).\]
Observe that the first short exact sequence associated to the $\JJ$-adic filtration on $\TT$
\[ 0\ra \JJ/\JJ^2 \ra \TT/\JJ^2 \ra \TT/\JJ \ra 0\]
always splits as $\TT/\JJ\cong \OO$ and $\TT$ is an $\OO$-algebra. Thus, we may choose $\psi$ so that $\varphi\circ\psi$ is a \emph{holomorphic} isomorphism
\[ (\varphi\circ\psi)^{\leq 1}: C^\infty(\TT/\JJ^2) \cong \hat{T}^{\leq 1} (\Omega^{1,0}_X).\]
Now the connection operator $D$ on $\TT$ induces a connection $\DD$ on $\hat{T}(\Omega^{1,0}_X)$ via the isomorphism $\varphi\circ \psi$. The component \[\Omega^{1,0}_X \ra \Lambda_X^{0,1}\ot \Omega^{1,0}_X\]
of $\DD$ is simply $\overline{\partial}$ due to the holomorphicity of $(\varphi\circ\psi)^{\leq 1}$. Thus, we conclude that $\DD$ is an NC-connection as described in Definition~\ref{nc-connection-def}.\ed

\subsection{Analytic NC manifolds from K\" ahler manifolds with constant holomorphic sectional curvature.} 
We are going to show that the above construction of NC-smooth thickenings is applicable to
K\"ahler manifolds with constant holomorphic sectional curvature. First, let us recall the definition of such manifolds. Since the K\"ahler form $g$ is nondegenerate $(1,1)$-form, it induces an operator (which is in fact an isomorphism) 
\[ \rho: T^{1,0}_X \ra \Omega_X^{0,1}.\]
Then $X$ is said to have constant holomorphic sectional curvature if the curvature operator $R: \Sym^2T^{1,0}_X \ra \Omega_X^{0,1}\ot T^{1,0}_X$ can be written as
\[ R(V,W) = c\rho(V)\otimes W +c\rho(W)\ot V\] 
for some constant $c\in \R$. If we define an operator $A: T^{1,0}_X\otimes T^{1,0}_X \ra \Omega^{0,1}_X\ot T^{1,0}_X$ by 
\[ A(V,W):= c\rho(V)\otimes W,\]
the curvature operator $R$ is then the symmetrization of the operator $A$. 

We define an operator associated to the metric $g$ by
\[\DD_g:=d_\tau+\overline{\partial}+\nabla+A^\vee\] 
on the algebra $\Lambda_X\otimes_{C^\infty_X} \hat{T}_{C^\infty}(\Omega_X^{1,0})$ as the unique derivation extending the de Rham differential and acting on $\Omega_X^{1,0}$ by the diagram
\begin{diagram}
\Lambda^1 &&  \Lambda^1\otimes \Omega_X^{1,0}  &&\Lambda^1\otimes T^2\Omega_X^{1,0} & \\
&\luTo(2,2)~{d_\tau} & \uTo~{\overline{\partial}+\nabla} & \ruTo(2,2)~{A^\vee} \\
 & & \Omega_X^{1,0}.& 
\end{diagram}
Here $d_\tau(1\otimes \alpha):=\alpha\otimes 1$ is the same as in the previous subsection. The second operator $\overline{\partial}+\nabla$ is the Levi-Civita connection associated to the metric $g$. The last operator 
\[A^\vee: \Omega_X^{1,0} \ra \Lambda^1\otimes T^2\Omega_X^{1,0},\] 
which is dual to $A$ (over the de Rham algebra $\Lambda$), is explicitly given by left multiplication by $c\tau(g)$ where $\tau(g)$ is the K\"ahler form $g$ considered as an element in $\Lambda^1\otimes \Omega_X^{1,0}$. Observe that the operator $D_g$ defined here is of the form defined in Definition~\ref{nc-connection-def}. Next we prove that it squares to zero, which shows that $D_g$ is an NC-connection.

\medskip
\begin{prop}~\label{nc-kahler-prop}
Let $(X,g)$ be a K\"ahler manifold with constant holomorphic sectional curvature. Then we have  $D_g^2=0$.
\end{prop}

\Pf. We need to check that $(d_\tau+\overline{\partial}+\nabla+A^\vee)^2=0$. First it is clear that $d_\tau^2=0$. We then note that $[d_\tau,\overline{\partial}+\nabla]=0$ due to the K\"ahler condition which implies the torsion freeness of $\overline{\partial}+\nabla$. The next equation that 
\[   (\overline{\partial}+\nabla)^2=-[d_\tau, A^\vee] \]
is due to the fact that the curvature $R$ is the symmetrization of $A$. Next the equation
\[ [\overline{\partial}+\nabla, A^\vee ]=0\]
holds since the connection $\overline{\partial}+\nabla$ preserves the K\"ahler metric, and $A^\vee$ is defined to multiply with $c\tau( g )$. Finally we need to show that
\[ [A^\vee, A^\vee]=0,\]
which is equivalent to its dual $A$ being associative. The operator $A$ is associative since we have
\[ A(A(U,V),W)=c\rho(A(U,V))\otimes W=c^2\rho(U)\rho(V)\otimes W;\]
and
\[ A(U,A(V,W))= c\rho(U)\otimes A(V,W)= c^2\rho(U)\rho(V)\otimes W.\]
The proposition is proved.\ed

\begin{rem}
In the definition of $A^\vee$ we can also define it to be the right multiplication by $\tau( g )$.
\end{rem}

\begin{cor}
Let $(X,g)$ be a K\"ahler manifold with constant holomorphic sectional curvature. Then $X$ admits an analytic NC thickening.
\end{cor}

\Pf. This follows from the previous proposition and Proposition~\ref{if-part-prop}.\ed

\begin{exs} K\"ahler manifolds with constant sectional curvature are classified:

\noindent\textbf{(A.)} In the case $c=0$ such a manifold corresponds to the quotient of $(\C^n,\omega_0=\sqrt{-1}/2\sum dz_id\overline{z}_j)$ by a discrete subgroup of its isometry group. For example, complex tori are such
quotients.

\noindent\textbf{(B.)} If $c>0$, these are complex projective spaces.

\noindent\textbf{(C.)} If $c<0$, we get complex hyperbolic manifolds which are quotient of upper half space $H^n$ endowed with its hyperbolic metric by a discrete subgroup of its isometry group. An interesing example is the moduli space of cubic surfaces (see~\cite{ACT}).

\end{exs}

\section{Analytic vector bundles over NC manifolds}\label{analytic-vb-sec}

\subsection{NC algebraic geometry and analytic geometry}\label{GAGA-sec} 
We first define the analytification of NC vector bundles. Let $X$ be a smooth variety over $\C$, equipped
with an algebraic NC-thickening $\OO_X^{NC}$.
Let $E$ be an algebraic vector bundle over $X$ of rank $n$, i.e. a locally free $\OO_X$-module of rank $n$ in the Zariski topology of $X$. Let $E^{NC}$ be an algebraic NC thickening of $E$, i.e. a locally free $\OO_X^{NC}$-module of rank $n$ whose abelianization is $E$. 
The embedding of algebraic sections to analytic sections is a morphism of sheaves in the Zariski topology of $X$
\[ \OO_X^{NC} \hookrightarrow \OO_X^{nc}.\]
We define the analytification of $E^{NC}$ by
\[ E^{nc}=(E^{NC})^{an}:=\OO_X^{nc}\ot_{\OO_X^{NC}} E^{NC}. \]
Note that $E^{nc}$ is a locally free $\OO^{nc}$-module with respect to Zariski topology (and hence
with respect to the classical topology). It is clear that the abelianization of $E^{nc}$ is $E^{an}$.

Next we compare the algebraic $\Ext$ groups with the analytic ones, where in the analytic setting we
consider $\Ext$'s in the category of sheaves of $\OO^{nc}$-modules with respect to the classical topology.  

\begin{thm}
\label{ext-thm}
Let $X$ be a smooth projective scheme, and let $X^{NC}$ be a smooth NC thickening of $X$. Let $E^{NC}$ and $F^{NC}$ be two vector bundles on $X^{NC}$ of finite rank. Then there is a natural isomorphism
\begin{equation}\label{Ext-map-eq}
\Ext^*(E^{NC},F^{NC}) \ra \Ext^*(E^{nc}, F^{nc}),
\end{equation}
where the $\Ext$-groups are computed in the category of $\OO_X^{NC}$- and 
$\OO_X^{nc}$-modules, respectively.
\end{thm}

\Pf. We have natural isomorphisms (see \cite[Thm.\ 4.2.1]{G-HA})
\[ \Ext^i(E^{NC},F^{NC})\simeq H^i(X, \underline{\Hom}_{\OO^{NC}} (E^{NC},F^{NC})),\]
\[ \Ext^i(E^{nc},F^{nc})\simeq H^i(X^{an}, \underline{\Hom}_{\OO^{nc}} (E^{nc},F^{nc}))\]
Since the classical topology is finer than the Zariski topology, we have a natural map
$$H^i(X, \underline{\Hom}_{\OO^{NC}} (E^{NC},F^{NC}))\to H^i(X, \underline{\Hom}_{\OO^{nc}} (E^{nc},F^{nc}))
\to H^i(X^{an}, \underline{\Hom}_{\OO^{nc}} (E^{nc},F^{nc})),$$
where in the middle we consider cohomology with respect to the Zariski topology, hence we get a map
\eqref{Ext-map-eq}.

Let $d\in \Z_{\geq 0}$ be a nonnegative integer. We first prove the statement for vector bundles over $X^{(d)}$ by induction on $Q$. The case $d=0$ is the usual GAGA theorem (see \cite{GAGA}). In the following we set $G_{d+1}:=(U\LLie_+(\Omega_X^1))_{d+1}$. By Theorem~\ref{module-thm}(iii), we have a short exact sequence
\[ 0 \ra G_{d+1}\ot \underline{\Hom}(E,F) \ra \underline{\Hom}(E^{(d+1)},F^{(d+1)}) \ra \underline{\Hom}(E^{(d)},F^{(d)}) \ra 0,\]
and in the analytic setting we have
\[ 0 \ra G^{an}_{d+1}\ot \underline{\Hom}(E^{an},F^{an}) \ra \underline{\Hom}((E^{(d+1)})^{an},(F^{(d+1)})^{an}) \ra \underline{\Hom}((E^{(d)})^{an},(F^{(d)})^{an}) \ra 0.\]
The analytification map induces a morphism of the associated long exact sequences of cohomology
\begin{diagram}
\Ext^{i-1,(d)} &\rTo & H^i(G_{d+1}\ot \underline{\Hom}(E,F)) &\rTo& \Ext^{i,(d+1)} & \rTo&\Ext^{i,(d)} & \rTo& H^{i+1}(G_{d+1}\ot \underline{\Hom}(E,F))\\
\dTo{\cong} & & \dTo{\cong} && \dTo{} && \dTo{\cong} && \dTo{\cong} \\
\Ext_{an}^{i-1,(d)} &\rTo & H^i(G^{an}_{d+1}\ot \underline{\Hom}(E^{an},F^{an})) &\rTo& \Ext_{an}^{i,(d+1)} & \rTo&\Ext_{an}^{i,(d)} & \rTo& H^{i+1}(G^{an}_{d+1}\ot \underline{\Hom}(E^{an},F^{an}))
\end{diagram}
where $\Ext^{i,(d)}:=\Ext^i(E^{(d)},F^{(d)})$ (resp., $\Ext_{an}^{i,(d)}=\Ext^i((E^{(d)})^{an},(F^{(d)})^{an})$). 
The induction assumption together with the usual GAGA implies that 
all the vertical arrows except for the middle one are isomorphisms.
Hence, by the five lemma, the middle vertical arrow is also an isomorphism. 
To prove the theorem, we observe that 
\[ \underline{\Hom}(E^{NC},F^{NC})= \liminv \underline{\Hom}(E^{(d)},F^{(d)}).\]
Since, the projective system $(H^0(U,\underline{\Hom}((E^{(d)})^{an},(F^{(d)})^{an})))$ satisfies the Mittag-Leffler condition for every open affine subset $U$, 
as in the proof of~\cite[Lemma 1.1.6]{KS}, we get the short exact sequences
\[0\ra R^1\liminv \Ext^{i-1}(E^{(d)},F^{(d)}) \ra \Ext^i(E^{NC},F^{NC}) \ra \liminv \Ext^i(E^{(d)},F^{(d)}) \ra 0.\]
Similarly, in the analytic case we have
\[0\ra R^1\liminv \Ext^{i-1}((E^{(d)})^{an},(F^{(d)})^{an}) \ra \Ext^i(E^{nc},F^{nc}) \ra \liminv \Ext^i((E^{(d)})^{an},(F^{(d)})^{an}) \ra 0.\]
As was shown above, the analytification morphism induces isomorphisms on the leftmost and the rightmost groups. Thus, the groups in the middle are also be isomorphic.\ed

Let us denote by $\Per^{gl}(X^{NC})$ the global perfect derived category of finite complexes of locally free right $\OO_X^{NC}$-modules of finite rank, and let $\Per^{gl}(X^{nc})$ be the similar category for $\OO_X^{nc}$ (and the classical topology of $X^{an}$). Theorem~\ref{ext-thm} above shows that the analytification functor
\[\An: \Per^{gl}(X^{NC}) \ra \Per^{gl}(X^{nc})\]
is fully-faithful. In the next theorem we will prove that $\An$ is also essentially surjective. 
Hence, the functor $\An$ is an equivalence of categories as claimed in Theorem \ref{nc-gaga-thm}.

\begin{thm}\label{nc-bundle-thm}
Let $(X,\OO^{NC})$ be an algebraic NC-smooth thickening of a smooth projective variety $X$, and let $V$ be a locally free $\OO^{nc}$-module of finite rank over its analytification $X^{an}$. Then there exists a
locally free $\OO^{NC}$-module $E^{NC}$ over $X$ such that its analytification $E^{nc}$ is isomorphic to $V$.
\end{thm}

\Pf. Since $V=\liminv V^{(d)}$, it is enough to construct a projective system of locally free $\OO^{(d)}$-modules
$E^{(d)}$ whose analytification gives $(V^{(d)})$ (then we can set $E=\liminv E^{(d)}$).
We will construct the $\OO^{(d)}$-modules $E^{(d)}$
inductively using the $I$-filtration. The existence of $E^{(0)}$ follows from the usual GAGA theorem. Assume that there exists a locally free $\OO^{(d)}$-module $E^{(d)}$ whose analytification is isomorphic to $V^{(d)}$. We would like to construct a locally free $\OO^{(d+1)}$-module $E^{(d+1)}$ whose analytification is isomorphic to $V^{(d+1)}$. For this consider the following exact sequence of sheaves of (nonabelian) groups in the Zariski topology of $X$:
\[ 0\rightarrow \Mat_{r} \big((U\LLie_+(\Omega_X^1))_{d+1}\big)\rightarrow GL_r(\OO^{(d+1)}) \rightarrow GL_r(\OO^{(d)})\rightarrow 0.\]
According to the general formalism of nonabelian cohomology (see e.g., \cite{Serre}),
the locally free $\OO^{(d)}$-module $E^{(d)}$ corresponds to a class in $c\in H^1(X,GL_r(\OO^{(d)}))$,
and there is a naturally defined obstruction
$$\de(c)\in H^2\big(X,(U\LLie_+(\Omega_X^1))_{d+1}\ot \und{\End}(E)\big),$$ 
such that $\de(c)=0$ if and only if $c$ lifts to a class in $H^1(X,GL_r(\OO^{(d+1)}))$ (note that here
the sheaf $\und{\End}(E)$ appears as the twisted version of the sheaf $\Mat_{r}(\OO_X)$).
Furthermore, such liftings of $c$ form a homogeneous space for the natural action of
$H^1\big(X,(U\LLie_+(\Omega_X^1))_{d+1}\ot \und{\End}(E)\big)$.
Now we observe that by the usual GAGA, the analytification map induces isomorphisms
$$H^i\big(X,(U\LLie_+(\Omega_X^1))_{d+1}\ot \und{\End}(E)\big)\to 
H^i\big(X^{an},(U\LLie_+(\Omega_X^1)^{an})_{d+1}\ot \und{\End}(E^{an})\big)$$
for $i=1,2$. Since the analytification map is compatible with the above cohomological constructions,
we derive the existence of $E^{(d+1)}$ with desired properties.
\ed

\begin{cor} (of Theorem \ref{nc-gaga-thm}). In the context of Theorem \ref{nc-bundle-thm}, 
the pair consisting of a locally free $\OO^{NC}$-module $E^{NC}$ together with an isomorphism $E^{nc}\to V$, is unique up to a unique isomoprhism.
\end{cor}

\section{Analytic NC-Jacobian and $A_\infty$-structures}\label{ainf-sec}

\subsection{NC thickening of moduli space of line bundles.} Let $Z$ be a compact complex manifold,
and let $L$ be a holomorphic line bundle over $Z\times B$, where $B$ is another 
complex manifold. We assume that

\medskip
\noindent
($\star$) the Kodaira-Spencer map $\kappa: T_B \to H^1(Z,\OO)\ot\OO_B$ associated with the family $L$ is an isomorphism.

\medskip

Set $E=\bigoplus_{i\ge 0}E^i:=H^\bullet(Z,\OO)$.
We choose a Hermitian structure on $Z$, and use Hodge theory to obtain harmonic representatives of 
$E$ in the differential graded algebra $(\Omega^{0,*},\dbar)$. Moreover, we have maps
\begin{diagram}
E &\pile{\rTo{i}\\ \lTo_p} & \Omega^{0,*}
\end{diagram}
such that $p\circ i=\id$ and $i\circ p= \id+\dbar h+h \dbar$,
where $h: \Omega^{0,*} \to \Omega^{0,*}$ is a homotopy operator. 
Using the Kontsevich-Soibelman's tree formula (see \cite{KoSo})
we get a minimal $A_\infty$-algebra structure on $E$.

Let us consider the (trivial) graded vector bundle over $B$, 
$$\EE=\bigoplus_{i\ge 0}\EE^i:=E\otimes_\C \OO_B.$$ 
The $A_\infty$-structure on $E$ induces by extension of scalars
an $A_\infty$-algebra structure on $\EE$, i.e., we have $\OO_B$-linear morphisms
\[ m_k: \EE^{\otimes k} \ra \EE\]
that satisfy the $A_\infty$-identities (with $m_1=0$).

The bundle $\EE$, being canonically trivialized, has a flat connection $\nabla$, and
the maps $m_k$ are $\nabla$-horizontal. Let us denote by $m_i\ot \id$ the operations on
the holomorphic de Rham complex of $\EE$,
$$\EE_\Omega:=\EE\otimes_{\OO_B}\Omega^*_B,$$ 
obtained from $m_i$ by extending scalars from $\OO_B$ to $\Omega^*_B$.
Let us set
\[ \tilde{m}_1=\nabla,\;\;\;\; \tilde{m}_k= m_k\otimes \id \mbox{\; for $k\geq 2$},\]
where in the first formula we extend $\nabla$ to a differential on the de Rham complex of
$\EE$ in the standard way.
This is almost the $A_\infty$-structure on $\EE_\Omega$ we want, but not quite: we are going to add an $m_0$ term. 

For this we observe that there is a canonical element $\omega\in \EE^1\otimes \Omega^1_B$. Indeed, by Assumption 
($\star$), the bundle $\EE^1=H^1(Z,\OO)\otimes_\C \OO_B$ is canonically isomorphic to the tangent bundle of $B$. We let $\omega$ be the element corresponding to $\id\in T_B\ot \Omega^1_B=\End(T_B)$. 
In local Kuranishi coordinates (see for example~\cite[Section 1.5, 1.6]{Fukaya}) $y_1,\ldots, y_n$ on $B$ corresponding to a basis $f_1,\ldots, f_n$ of $H^1(Z,\OO)$, the Kodaira-Spencer map is given by 
\[ f_i \mapsto \partial/\partial y_i.\]
Thus, in these coordinates we have $\omega=\sum_{i=1}^n f_i\ot dy_i$.
Note that since $f_i$ are $\nabla$-horizontal, we have
\begin{equation}\label{nabla-omega-eq}
\nabla(\omega)=0.
\end{equation}

Now we add a curvature term 
to the $A_\infty$-structure $\tilde{m}_k$ on $\EE_\Omega$ by putting $\tilde{m}_0=\omega$. 

\begin{lem}
The maps $\left\{\tilde{m}_k\right\}_{k=0}^{\infty}$ form a curved $A_\infty$-structure on the de Rham complex $\EE_\Omega$, with the respect to the total grading on $\EE_\Omega$.
\end{lem}

\Pf. For each $N\geq 0$, we need to prove the $A_\infty$-identity
\begin{equation}\label{A-infty}
\sum_{i+l+j=N, k=i+j, l\geq 0} (-1)^{i+lj}
\tilde{m}_{k+1}(\id^i\otimes \tilde{m}_l \otimes \id^j)=0.
\end{equation}
Since the maps $\tilde{m}_n$ for $n\ge 2$ are $\nabla$-horizontal, using the $A_\infty$-identity
for $(m_n)$ we obtain
\[ \sum_{i+l+j=N, k=i+j, l\geq 1} (-1)^{i+lj}
\tilde{m}_{k+1}(\id^i\otimes \tilde{m}_l \otimes \id^j)=0.\]
To prove \eqref{A-infty} it remains to check the identity
\[\sum_{i+j=k} (-1)^i \tilde{m}_{k+1}(\id^i\ot \omega\ot \id^j)=0.\]
In the case $k=0$ this identity follows from the equation \eqref{nabla-omega-eq}.
To deal with the cases $k\geq 1$ we use the fact that the differential graded algebra $\Omega^{0,*}$ is (super)commutative. Thus, by Cheng-Getzler's $C_\infty$ transfer theorem~\cite[Theorem 12]{CG} the structure maps $m_k$ obtained by Kontsevich-Soibelman's tree formula in fact form a $C_\infty$-structure, 
i.e., an $A_\infty$-structure that vanishes on the image of the shuffle product. 
Now we observe that the operator $\sum \tilde{m}_{k+1}(\id^i\ot \omega\ot \id^j)$, 
when applied to an element $\alpha_1\otimes\cdots\ot \alpha_k$, is equal (up to a sign) to
\[ \tilde{m}_{k+1}\big((\alpha_1\ot\cdots\ot\alpha_k)\bullet\omega\big),\]
where $\bullet$ denotes the shuffle product map. So the above sum vanishes by Cheng-Getzler's theorem.\ed

\begin{rem} In the case when $Z$ is K\"ahler, by \cite{DGMS}, the dg-algebra $(\Om^{0,*},\dbar)$ is formal.
Hence, in this case we can assume that $m_k=0$ with $k>2$, so $\EE_\Omega$ is a curved dg-algebra.
\end{rem}

Thus, $\EE_\Omega$ is a (curved) 
$A_\infty$-algebra with structure maps $\left\{\tilde{m}_k\right\}_{k=0}^{\infty}$. We view $\EE_\Omega$
as a strictly unital $A_\infty$-algebra over $\Omega^*_B$ via the canonical embedding
\[ \Omega^*_B\to \EE_\Omega=H^*(Z,\OO)\otimes_\C \Omega^*_B\]
defined by $\alpha \mapsto \id_\OO\otimes \alpha$.

\begin{lem}
The canonical projection $\EE_\Omega \to H^0(Z,\OO)\otimes \Omega_B^*=\Omega_B^*$ is an augmentation of $\EE_\Omega$ as an $A_\infty$-algebra over $\Omega^*_B$.
\end{lem}

\Pf. This follows from the fact that the $A_\infty$-algebra $E=H^*(Z,\OO)$ is augmented, and that the curvature $\omega$ maps to zero under this projection map. \ed

Let us denote by $\overline{E}:=\oplus_{k\geq 1}H^k(Z,\OO)$ the augmentation ideal. 
Similarly, we set $\overline{\EE}:= \overline{E}\otimes_\C \OO$. 

\begin{defi}
Let $\BB\EE_\Omega$ be the Bar construction of the augmented $A_\infty$-algebra $\EE_\Omega$ over
$\Omega^*_B$. 
This is a differential graded coalgebra over $\Omega_B^*$. 
We define $(\EE_\Omega)^!$, the Koszul dual algebra of $\EE_\Omega$ to be the dual of $\BB\EE_\Omega$, i.e. we have
\[ (\EE_\Omega)^!:= (\BB \EE_\Omega)^\vee=\hat{T}_\Omega (\Omega\otimes_\C \overline{E}^\vee[-1]),\]
where $\Omega=\Omega^*_B$.
\end{defi}

We are going to describe the sheaf of algebras $\underline{H}^0\big( (\EE_\Omega)^!\big)$ over $B$. In particular, we will see that if $H^2(Z,\OO)=0$ then the algebra $\underline{H}^0\big( (\EE_\Omega)^!\big)$ is a smooth NC thickening of $B$. We first exhibit the differential graded algebra $(\EE_\Omega)^!$ as a double complex in the second quadrant. For this observe that the complex $(\EE_\Omega)^!$ is bigraded by $\Z^{\leq 0}\times \N$ where the first grading is induced from that of $E$, and the second grading is simply the degree of differential forms on $B$. More precisely, an element $\alpha\in H^i(Z,\OO)^\vee\otimes \Omega_B^j$ has bidegree $(1-i,j)$, so that the bigraded components of $(\EE_\Omega)^!$ are
\[ \AA^{i,j}:= \bigoplus_{k-(r_1+\ldots+r_k)=i, \;\; r_1,\ldots,r_k\geq 2} \big(\hat{T} \ot (\EE^{r_1})^\vee \ot \cdots \ot (\EE^{r_k})^\vee\ot \hat{T}\big)\otimes \Omega^{j},\]
where $\hat{T}=\hat{T}_\OO((\EE^1)^\vee)$. The complex $(\EE_\Omega)^!$ may be depicted in the $(i,j)$-plane as follows.

\[\minCDarrowwidth15pt\begin{CD}
@AAA                    @A D_\nabla AA @.   @A D_\nabla AA  @A D_\nabla AA  @.\\
\cdots  @>d_m>> \AA^{i,j}  @>>>\cdots @>>> \AA^{-1,j} @>d_m>>  \AA^{0,j} @>>> 0\\
@AAA                    @A D_\nabla AA @.   @A D_\nabla AA  @A D_\nabla AA  @.\\
\vdots   @>d_m>> \vdots  @>>> \cdots @>>> \vdots @>d_m>> \vdots @>>> 0\\
@AAA                    @A D_\nabla AA @.   @A D_\nabla AA  @A D_\nabla AA  @.\\
\cdots   @>d_m>> \AA^{i,0}  @>>> \cdots @>>> \AA^{-1,0} @>d_m>> \AA^{0,0} @>>> 0\\
\end{CD}\]
Here the differential $\nabla$ is the connection operator, $d_\omega$ is the Koszul dual differential defined by the curvature term $\omega$, and $d_m$ is the operator defined by the $A_\infty$-structure on $E$.

Since the double complex $\AA^{*,*}$ is bounded vertically by the dimension of $B$, we can calculate the cohomology of the total complex using 
the spectral sequence associated to the vertical filtration of $\AA^{*,*}$. Let us first consider the rightmost vertical complex
\[ \big( \AA^{0,*}, D_\nabla \big).\]
We claim that its cohomology is concentrated in degree $0$.
Indeed, note that $\AA^{0,*}=\Omega^*\otimes \hat{T}_\OO((\EE^1)^\vee)$. By Assumption 
($\star$), the dual of the Kodaira-Spencer map is an isomorphism $(\EE^1)^\vee\cong \Omega^1_B$. Moreover, under this identification the operator $d_\omega$ is the $\Omega^*$-linear derivation acting on generators by $1\otimes \alpha \mapsto \alpha \ot 1$. By results in Section~\ref{connection-constr-sec}, we know this complex has cohomology concentrated in degree $0$, which is by definition a smooth NC thickening of $B$. As before we denote this $0$th cohomology by $\OO^{nc}$.

Now consider the vertical complex $\AA^{i,*}$ with $i<0$. 
Note that the operator $d_\omega$ acts on $(\EE^{r})^\vee$ by zero if $r\geq 2$. 
Hence, the vertical cohomology of $\nabla+d_\omega$ is also concentrated  in degree zero, and is equal to
\[ \bigoplus_{k-(r_1+\ldots+r_k)=i, \;\; r_1,\ldots,r_k\geq 2} \OO^{nc} \ot_\C (E^{r_1})^\vee \ot_\C \cdots \ot_\C (E^{r_k})^\vee\ot_\C \OO^{nc}\]
at degree $(i,0)$. This implies that the first page of the spectral sequence is of the form
\[ \minCDarrowwidth15pt\begin{CD}
\cdots \bigoplus \OO^{nc} \ot_\C (E^{r_1})^\vee \ot_\C \cdots \ot_\C (E^{r_k})^\vee\ot_\C \OO^{nc} @> d_m>> \cdots @>d_m>> \OO^{nc}\ot (E^2)^\vee \ot \OO^{nc} @>>> \OO^{nc}.
\end{CD}\]
Therefore, the spectral sequence degenerates on the next page.
In particular, we get the following result.

\begin{prop}
Let $R$ be the image of $1\otimes E^2 \otimes 1$ under the morphism $d_m$ in $\OO^{nc}$. Then we have
\[ \underline{H}^0\big((\EE_\Omega)^!\big) \cong \OO^{nc}/\langle R \rangle\]
where $\langle R \rangle$ is the two-sided ideal generated by $R$.
\end{prop}

\begin{exs}

\textbf{A.} Assume the algebra $H^*(Z,\OO)$ is the exterior algebra in $H^1(Z,\OO)$,
and the higher products $m_k$ vanish (e.g., this is true if $Z$ is a complex torus).
Then $R$ is generated by elements of the form $e_i\otimes e_j-e_j\otimes e_i=[e_i,e_j]$,
where $e_1,\ldots, e_n$ is a basis of $(H^1(Z,\OO))^\vee$. Hence, we found that
\[ \underline{H}^0\big((\EE_\Omega)^!\big) \cong \OO^{nc}_{ab}=\OO_B.\]

\medskip
\noindent \textbf{B.} Assume that $H^2(Z,\OO)=0$. Then $R=0$, so we get 
\[\underline{H}^0\big((\EE_\Omega)^!\big)=\OO^{nc}.\]
\end{exs}

In general the sheaf of algebras $\underline{H}^0\big((\EE_\Omega)^!\big)$ might be quite complicated. Nevertheless, we record the following elementary property.

\begin{prop}
The ideal $\langle R\rangle$ is contained in $F^1\OO^{nc}$. Thus, 
$\und{H}^0\big((\EE_\Omega)^!\big)$ is an NC-thickening of $B$, where the morphism
$\underline{H}^0\big((\EE_\Omega)^!\big)\twoheadrightarrow \OO_B$ is induced by
the projection $\OO^{nc}\twoheadrightarrow \OO_B$.
\end{prop}

\Pf. Let $r=d_m(t)$ for some $t\in (E^2)^\vee$. Then we have $r\in\hat{T}^{\geq 2} \big((E^1)^\vee\big)$ since $d_m$ is dual to the $A_\infty$-structure maps $m_2, m_3,\ldots$ which all have two or more inputs. Now the assertion follows from Proposition~\ref{filtration-prop}(i). \ed

\subsection{Analytic construction of the NC Poincar\'e bundle.} In this subsection we specialize to the case when $Z=C$ is a smooth curve, and $B=J(C)$ (or simply $J$) is the Jacobian variety of $C$. Let $\PP$ be a universal line bundle over $C\times J$. We will give an analytic construction of an NC-thickening of $\PP$, and compare it with the algebraic thickening constructed in Section \ref{NC-Poincare-sec}.

For a point $L\in J(C)$, we first give an analytic description of the universal line bundle $\PP|_{C\times V}$ for some analytic neighborhood $V$ of $L$. The line bundle $L$ over $C$ has a resolution by a two term complex of the form
\[ R_L:=\Omega_C^{0,0} \stackrel{\overline{\partial}+ a_L}{\longrightarrow} \Omega_C^{0,1}\]
for some harmonic element $a_L\in \Om_C^{0,1}$.

A local family of deformations of $L$ can be described as follows. We choose a basis $f_1,\ldots, f_n$ of of the vector space $H^1(C,\OO)$, and let $y_1,\ldots, y_n$ be the dual linear coordinates on this vector space. We consider $H^1(C,\OO)$ as a subspace of $\Omega_C^{0,1}$ consisting of harmonic forms. Let $V$ be a small neighborhood of the origin in $H^1(C,\OO)$. 
Let us consider the following family of two term complexes parametrized by $(y_1,\ldots, y_n)\in V$:
\[ \begin{CD}
 R_{L(y)}:=\Omega_C^{0,0} @>\overline{\partial}+ a_L+\sum_i y_i f_i>> \Omega_C^{0,1}.
\end{CD}\]
Note that $R_{L(y)}$ is a resolution of a deformation $L(y)$ of $L$. 
To describe the bundle $\PP|_{C\times V}$ we use the relative Dolbeault complex
\[ \begin{CD}
R_{\PP|_{C\times V}}:= \Omega_\pi^{0,0} @>\overline{\partial}+ a_L+\sum_i y_i f_i>> \Omega_\pi^{0,1}
\end{CD}\]
where $\Om^{0,*}_\pi$ are differential forms in $C$-direction, holomorphic in $V$-direction:
\[ \Omega^{0,*}_\pi:= \left\{ s\in q^*\Omega_C^{0,*} | \overline{\partial}_{y_i} s=0, \forall i=1,\ldots, n. \right\} \]
where $\pi: C\times V \ra V$ and $q: C\times V\ra C$ are the projections. Note that $\Omega^{0,*}_\pi$ is a sheaf of $\OO_{C\times V}$-modules, and that it has a relative holomorphic flat connection $\nabla$ in the $V$-direction, i.e., an operator
\[ \nabla: \Omega^{0,*}_\pi \ra \Omega^{0,*}_\pi\otimes \pi^*\Omega^1_V\]
satisfying the Leibniz rule and such that $\nabla^2=0$. However, $\nabla$ does not commute with the differential $\overline{\partial}+ a_L+\sum_i y_i f_i$ on $\Omega^{0,*}_\pi$, so it does not induce a relative flat connection on $\PP|_{C\times V}$. 

Note that the connection $\nabla$ does commute with $\overline{\partial}+a_L$,
and that the complex
\[ \begin{CD}
R_{q^*L}:= \Omega_\pi^{0,0} @>\overline{\partial}+a_L>> \Omega_\pi^{0,1}
\end{CD}\]
is the Dolbeault resolution of the line bundle $q^*L$. The induced connection in $V$-direction on $q^*L$ 
coincides with the obvious one.

The relative flat connection $\nabla$, by the construction in subsection~\ref{D-mod-sec}, induces a partial analytic $NC$-smooth thickening of $C\times V$ in the $V$-direction. By definition this partial thickening is $\underline{H}^0$ of a differential graded algebra $\pi^*\AA_V$ where $\AA_V=\Omega^*_V\otimes_{\OO_V}\hat{T}(\Omega^1_V)$. Its differential is the (analytic) NC-connection $D_\nabla$ associated to the flat connection $\nabla$.  Viewing $R_{q^*L}$ is a complex of $\pi^{-1}D_V$-modules, 
by the construction of Definition~\ref{D-mod-defi}, we get a structure of a
right dg-module over $\pi^*\AA_V$ on the completed tensor product 
\[R_{q^*L}\hat{\otimes} \pi^*\AA_V=R_{q*L} \hat{\otimes}_{\OO_{C\times V}} \pi^*\big(\Omega^*_V\otimes_{\OO_V}\hat{T}(\Omega^1_V)\big)\]
giving a thickening of $q^*L$. To obtain a thickening of $\PP|_{C\times V}$, we replace the tensor product differential 
$\overline{\partial}+a_L+D_\nabla$  with
\begin{equation}\label{deformed-differential-eq}
\overline{\partial}+ a_L+D_\nabla+\sum_i (y_i-e_i) f_i,
\end{equation}
where $e_i:=dy_i$ are horizontal sections of $\Omega^1_V\sub\hat{T}(\Omega^1_V)$, and the sum $\sum_i (y_i-e_i) f_i$ acts by left multiplication. Let us denote by 
$$R_{\PP|_{C\times V}}\hat{\ot}\pi^*\AA_V$$
this complex with the differential given by \eqref{deformed-differential-eq}. 

\begin{lem}
The complex $R_{\PP|_{C\times V}}\hat{\ot}\pi^*\AA_V$ is a right dg-module over $\Omega^{0,*}_\pi\hat{\ot}\pi^*\AA_V$.
\end{lem}

\Pf. This follows from the fact that 
\begin{equation}\label{MC-element-eq}
\a=a_L+\sum_i (y_i-e_i) f_i
\end{equation} 
is a Maurer-Cartan element  for the dg-algebra 
$$(\Omega^{0,*}_\pi\hat{\ot}\pi^*\AA_V, \ov{\pa}+D_\nabla).$$
\ed

\begin{rem}~\label{vanishing-rem}
Theorem~\ref{coh-vanishing} in the relative setting implies that $\underline{H}^i(\pi^*\AA_V)=0$ for $i>0$. Let us set $\OO_{C\times V^{nc}}:= \underline{H}^0(\pi^*\AA_V)$. This is an NC-thickening of $\OO_{C\times V}$ in the $V$-direction.
\end{rem}

\vspace{1mm}

\begin{prop}~\label{analytic-constr-prop}
The dg-module $R_{\PP\mid_{C\times V}}\hat{\ot}\pi^*\AA_V$ is locally trivial of rank one over the dg-algebra $\Omega^{0,*}_\pi\hat{\ot}\pi^*\AA_V$.
\end{prop}

\Pf. By definition, this dg-module is the rank one twisted complex associated with the Maurer-Cartan element
$\a$ given by \eqref{MC-element-eq}
of $A:=\Omega^{0,*}_\pi\hat{\ot}\pi^*\AA_V$ whose differential is $\overline{\partial}+D_\nabla$. Thus,
it is enough to prove that the Maurer-Cartan element $\alpha$ is locally gauge equivalent to zero. Let $U\subset C$ be a Stein open subset of $C$. We will prove that the Maurer-Cartan moduli space of the dg algebra $\Gamma(U\times V, A)$ is just a one point set. Thus every Maurer-Cartan element of it is gauge equivalent to zero.

For this, we make use of the following well-known result. See for example~\cite[Theorem 2.1]{Getzler}.

\begin{lem}~\label{mc-lem}
Let $h=F^0 h\supset F^1 h \supset F^2 h \supset \ldots$ and $g=F^0g \supset F^1g\supset F^2g\supset \ldots$ be filtered dg Lie algebras (that is, $dF^ih\subset F^ih$ and $[F^ih,F^jh]\subset F^{i+j+1}h$, and likewise for $g$) such that $h$ and $g$ are complete with respect to the filtrations. Let $f: h\ra g$ be a morphism of filtered dg Lie algebras which induces weak equivalences of the associated chain complexes
\[ \gr^i f: F^ih/F^{i+1}h \ra F^ig/F^{i+1}g.\]
Then $f$ induces a bijection between the associated Maurer-Cartan moduli spaces.
\end{lem}

Note that for a dg algebra $A$, by definition, its Maurer-Cartan moduli space is the same as that of the dg Lie algebra $A^{\Lie}$. We apply the above Lemma to the case $h=\Gamma(U\times V, \OO_{C\times V^{nc}})^{\Lie}$ and $g=\Gamma(U\times V, A)^{\Lie}$, both endowed with the $F$-filtration. One checks easily that the $F$-filtration define a filtered dg Lie algebra structure on $h$ and $g$. The morphism $f: h\ra g$ is the canonical embedding. Moreover, by Remark~\ref{vanishing-rem}, $f$ satisfies the property required in Lemma~\ref{mc-lem}. It follows that the Maurer-Cartan space of $g$ is a point since the degree one part of $h$ is only zero.\ed

\begin{cor}
The sheaf $\underline{H}^0(R_{\PP|_{C\times V}}\hat{\ot}\pi^*\AA_V)$ is a locally free $\OO_{C\times V^{nc}}$-module of rank one. Moreover, it is an NC thickening of the line bundle $\PP|_{C\times V}$.
\end{cor}

\Pf. The first assertion follows immediately from the previous proposition. For the second, we observe that the Maurer-Cartan element $\a=a_L+\sum_i (y_i-e_i) f_i$
maps down to $a_L+\sum_i y_i f_i$ under the canonical projection map 
$$\Omega^{0,*}_\pi\hat{\ot}\pi^*\AA_V\ra \Omega^{0,*}_\pi.$$
Twisting $\Omega^{0,*}_\pi$ by the latter Maurer-Cartan element defines the line bundle $\PP|_{C\times V}$. \ed

The construction of an NC thickening of $\PP|_{C\times V}$ can be made global over $C\times J$. 

\medskip

\noindent
{\bf Construction}.
Recall that the Jacobian $J$ is isomorphic to $H^1(\OO)/H^1(C,\Z)$. Since we have chosen coordinates $y_1,\ldots,y_n$ on $H^1(\OO)$, we can cover $J$ by open subsets 
$V_\a$ with coordinates on them that only differ by translations.

Let us first recall the gluing data of the usual Poincar\'e bundle over $C\times J$. Locally over $C\times V_{\alpha}$, with affine coordinates $y_i^\alpha$ centered at the point $L_\alpha$, the Poincar\'e bundle is isomorphic to $\underline{H}^0$ of the complex
\[ \big(\Omega^{0,*}_\pi, \overline{\partial}+a_{L_\alpha}+ \sum_i y^\alpha_i f_i\big).\]
Over the intersection $C\times (V_{\alpha}\cap V_\beta)$ there is a transition isomorphism
\[ \phi_{\alpha\beta}: \big(\Omega^{0,*}_\pi, \overline{\partial}+a_{L_\alpha}+ \sum_i y^\alpha_i f_i\big) \ra \big(\Omega^{0,*}_\pi, \overline{\partial}+a_{L_\beta}+ \sum_i y^\beta_i f_i\big).\]
Since the two coordinates $y^\alpha$ and $y^\beta$ only differ by translation, the function $\phi_{\alpha\beta}$ is a function pulled back from $C$. The collection of transition isomorphisms $\left\{ \phi_{\alpha\beta} \right\}$ satisfies the cocycle condition.

In the NC case, recall that over $C\times V_\alpha$ we have the complex
\[\big(\Omega_\pi^{0,*}\hat{\ot}\pi^*\AA_{V_\alpha}, \overline{\partial}+ a_{L_\alpha}+\nabla+d_\omega+\sum_i (y^\alpha_i-e_i) f_i\big).\]
Similarly on $C\times V_\beta$, the differential is
\[\overline{\partial}+ a_{L_\beta}+\nabla+d_\omega+\sum_i (y^\beta_i-e_i) f_i.\]
The transition functions $\phi_{\alpha\beta}$ of the Poincar\'e bundle lift to $\phi_{\alpha\beta}\ot \id$ which still satisfies the cocycle condition. Thus, the local dg-modules $R_{\PP|_{C\times V_\alpha}}\hat{\ot}\pi^*\AA_{V_\alpha}$ can be pasted together to obtain a global dg $\pi^*\AA_J$-module which we denote by $R_\PP\hat{\ot}\pi^*\AA_J$. Taking $\underline{H}^0$ gives an NC thickening of the Poincar\'e bundle $\PP$ over $C\times J$,
which we denote by $\PP^{nc}$.

\begin{prop}
Let $p\in C$ be a fixed point on $C$, and assume that $\PP$ is trivialized over $p\times J$. Then $\PP^{nc}$ is also naturally trivialized over $p\times J^{nc}$.
\end{prop}

\Pf. We defined the transition functions of $\PP^{nc}$ to be simply $\left\{\phi_{\alpha\beta}\ot \id\right\}$. It follows that if we start with the Poincar\'e bundle $\PP$ which is trivialized over $p\times J$, then so will
be $\PP^{nc}$.\ed

\begin{prop}\label{analytic-Poincare-id-prop} There exists
an automorphism $\phi: J^{nc} \ra J^{nc}$, trivial on the abelianization,
such that the analytic Poincar\'e bundle $\PP^{nc}$ on $C\times J^{nc}$ constructed above is isomorphic
to the analytification $(\id_C\times\phi)^*\PP^{NC, an}$, where $\PP^{NC}$ is constructed in Section~\ref{NC-Poincare-sec}.
\end{prop}

\Pf . Indeed, this follows from the universality of $\PP^{NC}$ combined with the results of Section \ref{GAGA-sec}.
\ed

\subsection{NC integral transforms.} By Proposition \ref{analytic-Poincare-id-prop}, we can rewrite (up
to an automorphism of $J^{nc}$)
the NC-Fourier-Mukai transform $F^{NC}$ from Section~\ref{NC-Poincare-sec} as the functor associated with
the kernel $P^{nc}$ over $C\times J^{nc}$,
\[ \Phi_{P^{nc}}: D^b(C) \ra \D(J^{nc}),\]
where for a ringed space we set $\D(X)=D(\mod-\OO_X)$. Explicitly, we have
\[ \Phi_{P^{nc}}(E):=   R\pi_*\big(q^*E\otimes \PP^{nc}\big).\]

Let us fix a line bundle $L\in J$. Consider the formal neighborhood of $L$ in $J$, which can be
identified with the formal polydisc $\De=\Spf \C[[y_1,\ldots, y_n]]$, using
 the affine coordinates $(y_i)$. The construction of the functor $\Phi_{P^{nc}}$ can be modified 
 by replacing the NC manifold $J^{nc}$ with a formal one, $\Spf \C\dlb{y_1,\ldots, y_n}\drb$ (where we take the 
 formal spectrum of an NC-complete algebra in the sense introduced by Kapranov~\cite[Def.\ 2.2.2]{Kapranov}). 
 Namely, the formal scheme $\De=\Spf \C[[y_1,\ldots,y_n]]$ has a torsion-free flat connection such that
 $\nabla(\partial/\partial y_i)=0$. This gives a differential graded resolution $\AA_{\De}$ of the algebra $\C\dlb{y_1,\ldots, y_n}\drb$ with differential $D_\nabla$. It remains to construct a kernel $\PP^{nc}_L$ over $C\times \Spf\C\dlb{y_1,\ldots,y_n}\drb$.  As before, we denote by
 $q: C\times \Spf \C\dlb{y_1,\ldots, y_n}\drb\ra C$ and $\pi:C\times \Spf 
\C\dlb{y_1,\ldots, y_n}\drb\ra \Spf \C\dlb{y_1,\ldots, y_n}\drb$ the natural projections. 
 Similarly to the construction of $\PP^{nc}$, we define $\PP_L^{nc}$ as 
 $\underline{H}^0$ of the complex
\[ \Omega^{0,*}_\pi\hat{\ot}\pi^*\AA_{\De}\]
endowed with the differential $\overline{\partial}+a_L+D_\nabla+\sum_i(y_i-e_i)f_i$.  The kernel $\PP_L^{nc}$ induces a Fourier-Mukai functor
\[ \Phi_{L}: D^b(C) \ra \D(\Spf \C\dlb{y_1,\ldots, y_n}\drb).\]

In~\cite[Sec.\ 2.2-2.4]{P-BN} the first named author constructed another functor
\[ \FF_L: D^b(C) \ra \D(\Spf \C\dlb{y_1,\ldots, y_n}\drb)\]
using the $A_\infty$-structure on $D^b(C)$. 
More precisely, $\FF_L(E)=(H^*(C,E\ot L)\ot \C\dlb{y_1,\ldots, y_n}\drb, d)$ where the differential 
$d$ is constructed using the structure of an $A_{\infty}$-module on
$H^*(C,E\ot L)=\Ext^*(E^\vee,L)$ 
over the $A_\infty$-algebra $\Ext^*_C(L,L)$, by passing to the (completed) Koszul dual (see \cite[(2.1.2)]{P-BN}).
Since $\Ext^2_C(L,L)=0$, this $A_\infty$-algebra has trivial products\footnote{Even though $\Ext^*(L,L)$
has trivial $A_\infty$-algebra structure, $A_\infty$-modules over it are not the same
as ordinary modules over an associative algebra.},
so its Koszul dual is just $\C\dlb{y_1,\ldots, y_n}\drb$.

\begin{prop}
Let $E\in D^b(C)$ be a bounded complex of locally free sheaves. Then there is a quasi-isomorphism 
\[\FF_L(E) \ra \Phi_L(E).\]
\end{prop}

\Pf.  
To compute $R\pi_*$ we use the relative Dolbeault resolution of $q^*E\otimes \PP_L^{nc}$. Thus, 
$\Phi_L(E)$ is represented by the complex
\[ \pi_*(q^*E\ot \Om^{0,*}_\pi\hat{\ot} \pi^*\AA_{\De})=\pi_*(q^*E\otimes\Omega^{0,*}_\pi) \hat{\ot} \AA_{\De}.\]
This complex is endowed with the differential $d_E+\overline{\partial}+a_L+D_\nabla+\sum_i(y_i-e_i)f_i$. We separate this differential into two parts by putting
\[Q:=d_E+\overline{\partial}+a_L+D_\nabla, \mbox{\;\; and \;\;} \delta:=\sum_i(y_i-e_i)f_i.\]
We view $\delta$ as a perturbation added to the differential $Q$, and we are going to apply the Homological Perturbation Lemma (see~\cite[(1.1)]{HK}) to a certain homotopy retraction of $Q$. 
Note that the differential $Q$ consists of two parts: the first part $d_E+\overline{\partial}+a_L$ acting on $\pi_*(q^*E\otimes\Omega^{0,*}_\pi)$, and the second part $D_\nabla$ acting on $\AA_{\De}$. The complex $\pi_*(q^*E\otimes\Omega^{0,*}_\pi)$, as a vector space, is equal to
\[ \Omega^{0,*}(E)\ot_\C \C[[y_1,\ldots,y_n]].\]
The differential $d_E+\overline{\partial}+a_L$ only acts on $\Omega^{0,*}(E)$ and the corresponding cohomology is $H^*(C,E\ot L)$. Thus, the cohomology of the complex 
$$\big( \Omega^{0,*}(E)\ot_\C \C[[y_1,\ldots,y_n]], d_E+\overline{\partial}+a_L\big)$$ 
is equal to $H^*(C,E\ot L)\otimes_{\C} \C[[ y_1,\ldots, y_n]]$. Choosing Hermitian metrics and using harmonic representatives we obtain a homotopy retraction
\begin{diagram}
H^*(C,E\ot L) &\pile{\rTo{i}\\ \lTo_p} & \big(\Omega^{0,*}(E), d_E+\overline{\partial}+a_L\big)
\end{diagram}
such that $p\circ i=\id$ and $i\circ p= \id+(d_E+\overline{\partial}+a_L)h+h(d_E+\overline{\partial}+a_L)$ where $h:\Omega^{0,*}(E)\ra \Omega^{0,*}(E)$ is a homotopy. Extending the coefficients to $\C[[ y_1,\ldots, y_n]]$ yields a homotopy retraction whose morphisms we still denote by $i$, $p$, and $h$:
\begin{diagram}
H^*(C,E\ot L)\otimes_\C \C[[ y_1,\ldots, y_n]] &\pile{\rTo{i}\\ \lTo_p} & \Omega^{0,*}(E)\ot_\C \C[[y_1,\ldots,y_n]]=\pi_*(q^*E\otimes\Omega^{0,*}_\pi)
\end{diagram}
Moreover, since the maps $i$, $p$ and $h$ are obtained by extending coefficients to $ \C[[ y_1,\ldots, y_n]]$, these are morphisms of $D$-modules (in variables $y_1,\ldots,y_n$),
where the left-hand side is equipped with the flat connection so that elements in $H^*(C,E\ot L)\otimes 1$ are 
horizontal, and the $D$-module structure on the right-hand side is defined similarly. Thus, applying the construction in Theorem~\ref{module-thm}, the maps $i$, $p$ and $h$ further induce a homotopy retraction between the corresponding dg-modules over $\AA_{\De}$,
\begin{equation}\label{filtered-hom-retract}
\begin{diagram}
 H^*(C,E\ot L)\otimes_\C \AA_{\De} &\pile{\rTo{i}\\ \lTo_p} & \pi_*(q^*E\ot \Om^{0,*}_\pi)\hat{\ot} \AA_{\De}.
\end{diagram}
\end{equation}
The differential on the left-hand side is $D_\nabla$, while the differential on the right-hand side is 
exactly $Q=d_E+\overline{\partial}+a_L+D_\nabla$. Let us equip 
$\AA_{\De}=\Om^\bullet_{\C[[y_1,\ldots,y_n]]/\C}\dlb e_1,\ldots,e_n\drb$ with the decreasing
filtration associated with the grading given by $\deg(y_i)=\deg(dy_i)=\deg(e_i)=1$.
Then \eqref{filtered-hom-retract} is a filtered homotopy retraction to which we can
apply the Homological Perturbation Lemma:
adding the perturbation $\delta=\sum_i(y_i-e_i)f_i$ to $Q$ yields a perturbed differential on $H^*(C,E\ot L)\otimes_\C \AA_{\De}$ given by
\[ D_\nabla+   p \delta i +p\delta h\delta i +p\delta h\delta h\delta i+\ldots.\]
Moreover, the Homological Perturbation Lemma also yields a perturbation $i'$ of the morphism $i$ so that it is still a quasi-isomorphism
\[\big(H^*(C,E\ot L)\otimes_\C \AA_{\De}, D_\nabla+   p \delta i +\ldots\big) \stackrel{i'}{\ra} \big(\pi_*(q^*E\ot \Om^{0,*}_\pi)\hat{\ot} \AA_{\De},Q+\delta\big)\]
between the perturbed complexes. To prove our assertion we rewrite the infinite series $p \delta i +p\delta h\delta i +p\delta h\delta h\delta i+\ldots$ in terms of
the structure of an $A_\infty$-module over $\Ext^*(L,L)$ on $H^*(C,E\ot L)=\Ext^*(E^\vee,L)$,
as in the definition of $\FF_L(E)$ (see ~\cite[Section 2.4]{P-BN}). 
Here it is crucial that we use Kontsevich-Soibelman's tree formula to form these $A_\infty$-structures,
since as was shown in~\cite{Markl}, this tree formula agrees with the ordinary homological perturbation formula.  This implies an identification of operators acting on $H^*(C,E\ot L)\otimes_\C \AA_{\De}$:
\begin{equation}\label{series-b-nc-eq}
p \delta i +p\delta h\delta i +p\delta h\delta h\delta i+\ldots = \rho_1(b^{nc})+\rho_2(b^{nc},b^{nc})+\ldots
\end{equation}
where 
$$\rho_k: \Ext^*(L,L)^{\ot k}\to \End(H^*(C,E\ot L))$$
are the operators giving the $A_\infty$-module structure on $H^*(C,E\ot L)$, extended to $H^*(C,E\ot L)\otimes_\C \AA_{\De}$ by $\AA_{\De}$-linearity, and $b^{nc}:=\sum_i(y_i-e_i) f_i$. 

Recall (see Sec.\ \ref{NC-affine-space-sec}) that the cohomology of $(\AA_\De, D_\nabla)$ is concentrated in degree zero and is given by
the subalgebra 
$$K:=\C\dlb y_1-e_1,\ldots,y_n-e_n\drb\sub \C[[y_1,\ldots,y_n]]\dlb e_1,\ldots,e_n \drb.$$
It is easy to see (using the filtration introduced above) that the embedding of a subcomplex
$$\big(H^*(C,E\ot L)\ot K, p \delta i +p\delta h\delta i +p\delta h\delta h\delta i+\ldots\big) \hra
\big(H^*(C,E\ot L)\otimes_\C \AA_{\De}, D_\nabla+   p \delta i +\ldots\big)
$$
is a quasiisomorphism. It remains to observe that the right-hand side of \eqref{series-b-nc-eq} gives
the same differential on $H^*(C,E\ot L)\ot \C\dlb y_1-e_1,\ldots,y_n-e_n\drb$ as in the definition
of $\FF_L(E)$.
\ed

\end{document}